\documentclass[11pt]{article}
\usepackage{amsfonts}
\usepackage{amsgen,amsmath,amstext,amsbsy,amsopn,amssymb,amsthm,subfigure}
\usepackage[dvips]{graphicx}
\usepackage{comment}
\usepackage{enumitem}
\usepackage{booktabs}
\usepackage{blkarray}
\usepackage{algorithm}
\usepackage{algpseudocode}
\usepackage{url}
\usepackage{longtable}
\usepackage{bm}
\usepackage[usenames]{color}

\allowdisplaybreaks

\newtheorem{Theorem}{Theorem}

\newtheorem{Lemma}{Lemma}
\newtheorem{Remark}{Remark}
\newtheorem{Corollary}{Corollary}

\newtheorem{Proposition}{Proposition}

\newcommand{\be}{\begin{equation}}
\newcommand{\ee}{\end{equation}}
\newcommand{\bea}{\begin{eqnarray}}
\newcommand{\eea}{\end{eqnarray}}
\newcommand{\beas}{\begin{eqnarray*}}
	\newcommand{\eeas}{\end{eqnarray*}}

\newcommand{\bbR}{\mathbb{R}}

\newcommand{\Var}{{\rm Var}}

\newcommand{\bbP}{{\mathbb{P}}}
\newcommand{\bbE}{{\mathbb{E}}}

\newcommand{\argmin}{\mathop{\rm arg\min}}

\usepackage{prettyref}

\newrefformat{eq}{(\ref{#1})}
\newrefformat{chap}{Chapter~\ref{#1}}
\newrefformat{sec}{Section~\ref{#1}}
\newrefformat{alg}{Algorithm~\ref{#1}}
\newrefformat{fig}{Fig.~\ref{#1}}
\newrefformat{tab}{Table~\ref{#1}}
\newrefformat{rmk}{Remark~\ref{#1}}
\newrefformat{clm}{Claim~\ref{#1}}
\newrefformat{def}{Definition~\ref{#1}}
\newrefformat{cor}{Corollary~\ref{#1}}
\newrefformat{lmm}{Lemma~\ref{#1}}
\newrefformat{prop}{Proposition~\ref{#1}}
\newrefformat{app}{Appendix~\ref{#1}}
\newrefformat{hyp}{Hypothesis~\ref{#1}}
\newrefformat{th}{Theorem~\ref{#1}}
\newrefformat{ass}{Assumption~\ref{#1}}

\providecommand{\keywords}[1]{\textbf{\textit{Keywords---}} #1}

\allowdisplaybreaks

\begin{document}

\title{On the Non-asymptotic and Sharp Lower Tail Bounds of Random Variables}

\author{Anru R. Zhang$^1$ ~ and ~ Yuchen Zhou$^2$}

\date{(\today)}

\maketitle

\bigskip

\footnotetext[1]{Department of Statistics, University of Wisconsin-Madison and Department of Biostatistics and Bioinformatics, Duke University}

\footnotetext[2]{Department of Statistics, University of Wisconsin-Madison}

\begin{abstract}
	The non-asymptotic tail bounds of random variables play crucial roles in probability, statistics, and machine learning. Despite much success in developing upper bounds on tail probability in literature, the lower bounds on tail probabilities are relatively fewer. In this paper, we introduce systematic and user-friendly schemes for developing non-asymptotic lower bounds of tail probabilities. In addition, we develop sharp lower tail bounds for the sum of independent sub-Gaussian and sub-exponential random variables, which match the classic Hoeffding-type and Bernstein-type concentration inequalities, respectively. We also provide non-asymptotic matching upper and lower tail bounds for a suite of distributions, including gamma, beta, (regular, weighted, and noncentral) chi-square, binomial, Poisson, Irwin-Hall, etc. We apply the result to establish the matching upper and lower bounds for extreme value expectation of the sum of independent sub-Gaussian and sub-exponential random variables. A statistical application of signal identification from sparse heterogeneous mixtures is finally considered.
\end{abstract}


\section{Introduction}\label{sec:intro}

The tail bounds and concentration inequalities, which study how a random variable deviates from some specific value (such as the expectation), are ubiquitous in enormous fields, such as random matrix theory, high-dimensional statistics, and machine learning \cite{boucheron2013concentration,vershynin2010introduction,vershynin2018high,wainwright2019high}. Many upper bounds on tail probabilities, such as the Markov's inequality, Chebyshev's inequality, Hoeffding's inequality \cite{hoeffding1963probability}, Bernstein's inequality \cite{bernstein1924modification,bernstein1937certain}, and Bennett's inequality \cite{bennett1962probability}, have been well-regarded and extensively studied in literature. The Chernoff-Cram\`er bound (\cite[Theorem 1]{chernoff1952measure}, also see, e.g., \cite{boucheron2013concentration} and \cite{wainwright2019high}), with a generic statement given below, has been a basic tool to develop upper bounds of tail probabilities. 
\begin{Proposition}[Chernoff-Cram\`er Bound]\label{pr:chernoff} If $X$ is a real-valued random variable with moment generating function $\phi_X(t) = \mathbb{E}\exp(tX)$ defined for $t \in T\subseteq \mathbb{R}$. Then for any $x\in \mathbb{R}$, 
	\begin{equation*}
	\begin{split}
	\bbP\left(X\geq x\right) \leq \inf_{t\geq 0, t \in T}\phi_X(t) \exp(-tx);\quad	\bbP\left(X\leq x\right) \leq \inf_{t\leq 0, t \in T}\phi_X(t) \exp(-tx).
	\end{split}
	\end{equation*}
\end{Proposition}
The Chernoff-Cram\`er bound has been the beginning step for deriving many probability inequalities, among which the Hoeffding's inequality \cite{hoeffding1963probability}, Bernstein's inequality \cite{bernstein1924modification,bernstein1937certain,uspensky1937introduction}, Bennett's inequality \cite{bennett1962probability}, Azuma's inequality \cite{azuma1967weighted}, and McDiarmid's inequality \cite{mcdiarmid_1989} are well-regarded and widely used. 
Motivated by applications in high-dimensional statistics and machine learning, the Chernoff-Cram\`er bound has been widely used to develop various concentration inequalities for sub-Gaussian, sub-exponential random variables. 
These results and applications have been collected in recent textbooks (see, e.g., \cite{boucheron2013concentration}, \cite{lugosi2009concentration}, \cite{steele1997probability}, \cite{vershynin2010introduction}, \cite{wainwright2019high}) and class notes (see, e.g.,  \cite{pollard2015good}, \cite{roch2015moments}, \cite{sridharan2002gentle}).

Despite enormous achievements on upper tail bounds in literature, there are relatively fewer results on the corresponding lower bounds: how to find a sharp and appropriate $L>0$, such that $\bbP(X\geq x) \geq L$ (or $\bbP(X\leq x) \geq L$) holds?  Among existing literature, the Cram\`er-Chernoff theorem characterized the asymptotic tail probability for the sum of i.i.d. random variables (\cite{lyons1992random} and \cite{chernoff1952measure}, Theorem 1; also see \cite{van2000asymptotic}, Proposition 14.23): suppose $Z_1,\ldots, Z_k$ are i.i.d. copies of $Z$. Then, 
\begin{equation*}
\begin{split}
k^{-1}\log\left(\mathbb{P}(Z_1+\cdots +Z_k \geq ak)\right) \to \log\left(\inf_{t\geq 0} \mathbb{E}( e^{t(Z - a)})\right)
\end{split}
\end{equation*}
as $k\to\infty$. Berry-Esseen central limit theorem \cite{berry1941accuracy,esseen1942liapunov} provided a non-asymptotic quantification of the normal approximation residual for the sum of independent random variables: let $Z_1,\ldots, Z_k$ be i.i.d. copies of $Z$, where $\mathbb{E}Z = 0$ and $\mathbb{E}|Z|^3<\infty$; then for all $x\in \mathbb{R}$,
\begin{equation}\label{ineq:berry-esseen}
	\begin{split}
	\Phi(-x) - \frac{C\mathbb{E}|Z|^3}{\sqrt{k}(\Var(Z))^{3/2}} \leq \bbP\left(Z_1+\cdots + Z_k \geq \sqrt{k\Var(Z)}x\right) \leq \Phi(-x)+\frac{C\mathbb{E}|Z|^3}{\sqrt{k}(\Var(Z))^{3/2}},
	\end{split}
\end{equation}
where $\Phi(\cdot)$ is the cumulative distribution function of standard normal distribution. This lower bound is not universally sharp, as the left hand side of \eqref{ineq:berry-esseen} can be negative for $x \geq C\sqrt{\log(k)}$. \cite{slud1977distribution,csiszar1998method,cover2012elements} established upper and lower tail bounds for binomial distribution based on its probability mass function. Kolmogorov introduced the Bernstein-type lower bound for the sum of independent bounded random variables \cite[Lemma 8.1]{ledoux2013probability}. \cite{gluskin1995tail} established tail and moment estimates for the sums of independent random variables with logarithmically concave tails. 
\cite[Chapter 3]{de2012decoupling} considered decoupling by developing 
lower tail bounds of hypercontractive variables. 
\cite{bagdasarov1996reversion} studied a reversion of Chebyshev's inequality based on the moment generating function and properties of the convex conjugate. \cite{theodosopoulos2007reversion} proved a tight reversion of the Chernoff bound using the tilting procedure.  
	There is still a lack of easy-to-use and sharp lower bounds on tail probabilities for generic random variables in the finite sample setting.

Departing from existing results, in this paper, we introduce systematic and user-friendly schemes to develop non-asymptotic and sharp lower bounds on tail probabilities for generic random variables. The proofs can be conveniently applied to various settings in statistics and machine learning. We also discuss the following implications of the developed results. 
\begin{itemize}[leftmargin=*]
	\item 
	In Section \ref{sec:sum-of-independent}, we establish the lower bounds on tail probabilities for the weighted sums of independent sub-Gaussian and sub-exponential random variables. These new results match the classic Hoeffding-type and Bernstein-type upper bounds in the literature.

	\item In Section \ref{sec:example}, we study the matching upper and lower tail bounds for a suite of commonly used distributions, including gamma, beta, (regular, weighted, and noncentral) chi-square, binomial, Poisson, Irwin-Hall distributions. To the best of our knowledge, these are the first sharp lower bounds on tail probabilities for these distributions. Especially, we establish a reverse Chernoff-Cram\`er bound (in the forthcoming Lemma \ref{th:reverse-chernoff}) to better cope with binomial and Poisson distributions, which may be of independent interest. 
	\item 
	In Section \ref{sec:application}, we consider the applications of the established results. We derive the matching upper and lower bounds for extreme values of the sums of independent sub-Gaussian and sub-exponential random variables. A statistical problem of signal identification from heterogeneous mixtures is finally studied.
\end{itemize}

\section{Generic Lower Bounds on Tail Probabilities}\label{sec:generic}

We use uppercase letters, e.g., $X, Z$, to denote random variables and lowercase letters, e.g., $x, t$, to denote deterministic scalars or vectors. $x\vee y$ and $x \wedge y$ respectively represent the maximum and minimum of $x$ and $y$. We say a random variable $X$ is centered if $\mathbb{E} X = 0$. 
For any random variable $X$, let $\phi_X(t) = \mathbb{E} \exp(tX)$ be the moment generating function.
For any vector $u\in \mathbb{R}^k$ and $q \geq 0$, let $\|u\|_{q} = (\sum_i |u_i|^q)^{1/q}$ be the $\ell_q$ norm. In particular, $\|u\|_\infty = \max_i |u_i|$. We use $C, C_1, \ldots, c, c_1, \cdots $ to respectively represent generic large and small constants, whose actual values may differ from time to time. Throughout the paper, we present the previous results as propositions and our new results as theorems/corollaries.

For any centered random variable $X$, if $\phi_X(t) = \mathbb{E}\exp(tX)$ is finite for a range of $t \in T\subseteq \mathbb{R}$, by Taylor's expansion, we have $\phi_X(t) = 1 + t^2\mathbb{E}X^2 + o(t^2)$ for $t$ in a neighborhood of 0. Thus, there exist constants $c_1, C_1>0$, such that 
\begin{equation}\label{ineq:thm2-pre}
\begin{split}
\exp\left(c_1(\mathbb{E}X^2) t^2\right)\leq \phi_X(t) \leq \exp\left(C_1(\mathbb{E}X^2) t^2\right)
\end{split}
\end{equation}
holds for $t$ in a neighborhood of 0. The following Theorem \ref{thm: general} provides the matching upper and lower bounds on tail probabilities for any random variable $X$ satisfying Condition \eqref{ineq:thm2-pre} in a certain regime. 
\begin{Theorem}[Tail probability bound: Small $x$]\label{thm: general}
	Suppose $X$ is a centered random variable satisfying
	\begin{equation}\label{eq32}
	\begin{split}
	c_2\exp(c_1\alpha t^2)\leq \phi_{X}(t) \leq C_2\exp(C_1\alpha t^2), \quad \forall 0 \leq t \leq M,
	\end{split}
	\end{equation}
	where $C_1 \geq c_1 > 0, C_2, c_2 > 0$ are constants. 
	Suppose either of the following statements holds: (1) $\alpha M^2 \geq \frac{16\left(1 + \log\left(1/c_2\right)\right)}{c_1}$; (2) $c_2 = 1$ and $\alpha M^2 \geq c''$ for some constant $c'' > 0$. Then there exist constants $c, c', C > 0$ such that
	\begin{equation}\label{eq33}
	\begin{split}
	c\exp\left(-C x^2/\alpha\right) \leq \bbP\left(X \geq x\right) \leq C\exp\left(-cx^2/\alpha\right),\quad \forall 0 \leq x \leq c'M\alpha.
	\end{split}
	\end{equation}
	Moreover, if $\alpha>0, M = +\infty$, there exist constants $c, C > 0$, such that 
	\begin{equation}\label{eq34}
	\begin{split}
	c\exp\left(-Cx^2/\alpha\right) \leq \bbP\left(X \geq x\right) \leq C\exp\left(-cx^2/\alpha\right),\quad \forall x\geq 0.
	\end{split}
	\end{equation}
\end{Theorem}
The proof of Theorem \ref{thm: general} relies on the appropriate choice of values $(t, \theta)$ in Paley-Zygmund inequality \cite[Chapter 1.6]{kahane1993some}, \cite{paley1932note}: 
\begin{equation*}
	\begin{split}
	\bbP\left(\exp(tX) \geq \theta\bbE \exp(tX)\right) \geq (1 - \theta)_+^2\frac{(\bbE \exp(tX))^2}{\bbE \exp(2tX)}.
	\end{split}
\end{equation*}

\begin{Remark}
	Previously, \cite{gluskin1995tail} studied the tail and moment estimates for the sums of independent random variables. Different from their assumptions that $X$ can be written as the sum of symmetric, independent, and identically distributed random variables with log-concave tails, Theorem \ref{thm: general} studies the generic random variables with bounded moment generating function $\phi_X(t)$ for a range of $t$. In addition, \cite{gluskin1995tail} focused on the tail probability $\bbP(X\geq x)$ for large $x$, while Theorem \ref{thm: general} focuses on bounded $x$.
\end{Remark}
Next, we consider the tail probability for large $x$. 
\begin{Theorem}[Tail probability bound: Large $x$]\label{thm:general 2}
	Suppose $Y$ and $Z$ are independent random variables. $\phi_{Z}$ is the moment generating function of $Z$ and $\phi_{Z}(t) \leq e^{C_1\alpha t^2}$ for all $-M \leq t \leq 0$. $Y$ satisfies the tail inequality
	$\bbP\left(Y \geq x \right) \geq T(x)$ for all $x \geq w\alpha.$
	Then $X = Y + Z$ satisfies
	\begin{equation}
	\begin{split}
	\bbP(X \geq x) \geq \left(1 - \exp\left(-\min\left\{w^2/(16C_1^2), Mw/4\right\}\alpha\right)\right)T(2x), \quad \forall x \geq w\alpha/2.
	\end{split}
	\end{equation}
\end{Theorem}

Theorem \ref{thm:general 2} immediately yields the lower tail bound for the sum of independent random variables. 
\begin{Corollary}\label{cor:general 2}
	Suppose $Z_1, \dots, Z_k$ are centered and independent random variables. Assume $\phi_{Z_i}(t) \leq \exp(C_1t^2)$ for $-M \leq t \leq 0$, where $\phi_{Z_i}$ is the moment generating function of $Z_i$, $1 \leq i \leq k$. $Z_1$ satisfies the tail inequality
	\begin{equation*}
	\begin{split}
	\bbP\left(Z_1 \geq x \right) \geq T(x), \quad \forall x \geq wk.
	\end{split}
	\end{equation*}
	Then, $X = Z_1 + \dots + Z_k$ satisfies
	\begin{equation}
	\begin{split}
	\bbP(X \geq x) \geq \left(1 - \exp\left(-\min\left\{w^2/(16C_1^2), Mw/4\right\}k\right)\right)T(2x), \quad \forall x \geq wk/2.
	\end{split}
	\end{equation}
\end{Corollary}

\section{Tail Bounds for the Sums of Independent Random Variables}\label{sec:sum-of-independent}

With the generic lower bounds on tail probability developed in the previous section, we are in position to study the tail probability bounds for the sum of independent random variables. 
We first consider the tail probability bounds for weighted sums of independent sub-Gaussian random variables. The upper tail bound, which is referred to as the Hoeffding-type concentration inequality, has been widely used in high-dimensional statistics and machine learning literature (see, e.g., \cite{vershynin2010introduction}). 
\begin{Proposition}[Hoeffding-type inequality for the sum of sub-Gaussians]\label{pr:hoeffding-type}
	Suppose $Z_1,\ldots, Z_k$ are centered and independently sub-Gaussian distributed, in the sense that either of the following holds: (1) $\phi_{Z_i}(t) \leq \exp(Ct^2)$; (2) $\bbP(|Z_i|\geq x) \leq C\exp(-cx^2)$. Then $X = u_1Z_1 + \cdots + u_kZ_k$ satisfies
	\begin{equation*}
	\begin{split}
	\bbP\left(|X|\geq x\right) \leq \exp(-c'x^2/\|u\|_2^2), \quad \forall x>0.
	\end{split}
	\end{equation*}
\end{Proposition}
The proof of Proposition \ref{pr:hoeffding-type} can be found in \cite[Proposition 5.10]{vershynin2010introduction}. With the additional condition on the lower tail bound of each summand, we can prove the following lower bound on tail probabilities for the sums of independent sub-Gaussians. This result matches the classic Hoeffding-type inequality in Proposition \ref{pr:hoeffding-type}. 
\begin{Theorem}[Hoeffding-type inequality for sub-Gaussians: Matching lower bounds]\label{cor: sub-Gaussian}
	Suppose $Z_1, \dots, Z_k$ are independent, $\mathbb{E}Z_i = 0$, and $\phi_{Z_i}$ is the moment generating function of $Z_i$.  Suppose either of the following statements holds: (1) there exist constants $c_1, C_1 > 0$, such that 
	\begin{equation}\label{ineq:sub-gaussian-state1}
	\begin{split}
	\exp(c_1t^2) \leq \phi_{Z_i}(t) \leq \exp\left(C_1t^2\right), \quad \forall t \geq 0, 1 \leq i \leq k;
	\end{split}
	\end{equation}
	(2) there exist constants $c_2, C_2, c_3, C_3 > 0$ , such that 
	\begin{equation}\label{ineq:sub-gaussian-state2}
	\begin{split}
	\bbP(Z_i\geq x) \geq c_2\exp(-C_2x^2), \quad \bbP\left(|Z_i| \geq x\right)\leq C_3\exp(-c_3x^2),\quad \forall x\geq 0, 1\leq i \leq k.
	\end{split}
	\end{equation}
	Then there exist constants $c, c', C > 0$, such that for any fixed values $u_1,\ldots, u_k \geq 0$, $X = u_1Z_1 + \cdots + u_kZ_k$ satisfies
	\begin{equation}
	\begin{split}
	c\exp\left(-Cx^2/\|u\|_2^2\right) \leq \bbP\left(X \geq x\right) \leq \exp\left(-c'x^2/\|u\|_2^2\right), \quad \forall x \geq 0.
	\end{split}
	\end{equation} 
\end{Theorem}

Although Theorem \ref{thm:general 2} focuses on the right tail bound, i.e., $\bbP(X\geq x)$, similar results hold for the left tail, i.e., $\bbP(X\leq -x)$, by symmetry. We can further prove the following lower bound on tail probability for the sum of random variables with two-sided sub-Gaussian tails.
\begin{Corollary}
	Suppose $Z_1, \dots, Z_k$ are independent, $\mathbb{E}Z_i = 0$, $\phi_{Z_i}$ is the moment generating function of $Z_i$.  Suppose either of the following statements holds: (1) $\exp(c_1t^2) \leq \phi_{Z_i}(t) \leq \exp\left(C_1t^2\right)$ for constants $c_1, C_1 > 0$, $t \in \bbR$, $1 \leq i \leq k$; (2) $c_2\exp(-C_2x^2) \leq \bbP(Z_i\geq x)\wedge \bbP\left(Z_i \leq -x\right) \leq \bbP(Z_i\geq x)\vee \bbP\left(Z_i \leq -x\right)\leq C_3\exp(-c_3x^2)$ for constants $c_2, C_2, c_3, C_3 > 0$ and all $x \geq 0$, $1 \leq i \leq k$.
	Then there exist constants $c, c', C > 0$, for any fixed real values $u_1,\ldots, u_k$, $X = u_1Z_1 + \cdots + u_kZ_k$ satisfies
	\begin{equation*}
	\begin{split}
	& c\exp\left(-Cx^2/\|u\|_2^2\right) \leq \bbP\left(X \geq x\right) \leq C\exp(-cx^2/\|u\|_2^2),\quad \forall x>0,\\ 
	& c\exp(-Cx^2/\|u\|_2^2) \leq \bbP\left(X \leq -x\right) \leq C\exp(-cx^2/\|u\|_2^2), \quad \forall x > 0.
	\end{split}
	\end{equation*}
\end{Corollary}

On the other hand, the class of sub-Gaussian random variables considered in Theorem \ref{cor: sub-Gaussian} may 
fail to cover many useful random variables with heavier tails. 
To cover broader cases, sub-exponential distributions were introduced and widely used in literature (see the forthcoming Proposition \ref{pr:bernstein-type} for definition of sub-exponential distribution). The following Bernstein-type inequality (see, e.g., \cite[Proposition 5.16]{vershynin2010introduction}, \cite[Theorem 2.3]{boucheron2013concentration} and \cite[Proposition 2.2]{wainwright2019high}) is a classic result on the tail bound for the sum of sub-exponential random variables. 
\begin{Proposition}[Bernstein-type inequality for the Sums of independent Sub-exponentials (\cite{vershynin2010introduction}, Proposition 5.16)]\label{pr:bernstein-type}
	Let $Z_1,\ldots, Z_k$ be independent centered sub-exponential random variables in the sense that $\mathbb{E}\exp(tZ_i) \leq \exp(Ct^2)$ for all $|t|\leq c$. Then for every $u_1,\ldots, u_k$, $X = u_1Z_1+ \cdots  + u_k Z_k$ satisfies
	\begin{equation*}
	\begin{split}
	& \bbP\left(\left|X\right| \geq x\right) \leq 2\exp\left(-cx^2/\|u\|_2^2\right), \quad \forall 0\leq x \leq \|u\|_2^2/\|u\|_\infty;\\
	& \bbP\left(\left|X\right| \geq x\right) \leq 2\exp\left(- cx/\|u\|_\infty\right), \quad \forall x \geq \|u\|_2^2/\|u\|_\infty.
	\end{split}
	\end{equation*}
\end{Proposition}
With the additional lower bound on the moment generating function of each summand, we prove the following matching upper and lower tail bounds for the sum of independent sub-exponential random variables. 
\begin{Theorem}[Bernstein-type inequality for the sums of independent sub-exponentials: Matching upper and lower bounds]\label{cor: sub-exponential}
	Suppose $Z_1,\ldots, Z_k$ are centered independent random variables and $\phi_{Z_i}$ is the moment generating function of $Z_i$. Suppose $Z_i$'s are sub-exponential in the sense that there exist two constants $c_1, C_1 > 0$,
	\begin{equation*}
		\begin{split}
		\exp(c_1\alpha t^2) \leq \phi_{Z_i}(t) \leq \exp(C_1\alpha t^2),\quad \forall |t| \leq M, 1\leq i \leq k.
		\end{split}
	\end{equation*}
	Suppose $\alpha M^2 \geq c''$, where $c'' >0$ is a constant. If $u_1,\ldots, u_k$ are non-negative values, then $X = u_1Z_1 + \dots + u_kZ_k$ satisfies
    \begin{equation*}
    \begin{split}
    c\exp\left(-C\frac{x^2}{\alpha\|u\|_2^2}\right) \leq \bbP\left(X \geq x\right) \leq \exp\left(-\tilde{c}\frac{x^2}{\alpha\|u\|_2^2}\right)
    \end{split}
    \end{equation*}
    for any $0 \leq x \leq \frac{c'M\alpha\|u\|_2^2}{\|u\|_{\infty}}$. Here $c, c', C, \tilde{c} > 0$ are constants.
	
	In addition, if there exists one $Z_i$ ($1\leq i \leq k$) satisfying
\begin{equation*}
	    \begin{split}
	    \bbP\left(Z_i \geq x \right) \geq c_0\exp(-C_0Mx), \quad \forall x \geq 2c'M\alpha,
	    \end{split}
	\end{equation*}
and $u_i \geq c_2\|u\|_{\infty}$ for constant $c_2 > 0$,
then we have
	\begin{equation}\label{eq37}
	\begin{split}
	\bar{c}\exp\left(-\bar{C}\frac{Mx}{\|u\|_{\infty}}\right) \leq \bbP(X \geq x) \leq \exp\left(-\tilde{c}\frac{Mx}{\|u\|_{\infty}}\right), \quad \forall x \geq \frac{c'M\alpha\|u\|_2^2}{\|u\|_{\infty}},
	\end{split}
	\end{equation}
	where $\bar{c}, \bar{C} > 0$ are constants.	
\end{Theorem}

\section{Sharp Tail Bounds of Specific Distributions}\label{sec:example}

In this section, we establish matching upper and lower tail bounds for a number of commonly used distributions. 
\subsection{Gamma Distribution}

Suppose $Y$ is gamma distributed with shape parameter $\alpha$, i.e., $Y\sim \text{Gamma}(\alpha)$, where the density is 
\begin{equation}\label{eq:gamma-density}
\begin{split}
f(y; \alpha) = \frac{1}{\Gamma(\alpha)}y^{\alpha-1}e^{-y}, \quad y > 0.
\end{split}
\end{equation}
Here $\Gamma(\alpha)$ is the gamma function. Although the density of gamma distribution is available, it is highly non-trivial to develop the sharp tail probability bound in a closed form. Previously, \cite[Pages 27-29]{boucheron2013concentration} established an upper bound on the tail probability of gamma distribution:
	Suppose $Y\sim \text{Gamma}(\alpha)$ and $X = Y - \alpha$. Then 
	\begin{equation}\label{ineq:previous-gamma}
	\begin{split}
	\bbP\left(X \geq \sqrt{2\alpha t} + t\right)  \leq e^{-t},
	\quad \bbP\left(X \leq -\sqrt{2\alpha t}\right) \leq e^{-t}, \quad \forall t \geq 0.
	\end{split}
	\end{equation}
	\begin{equation*}
	\begin{split}
	\text{Or equivalently,} \quad \bbP\left(X \geq x\right) \leq \exp\left(-\frac{x^2}{x+\alpha + \sqrt{\alpha^2+2x\alpha}}\right),\quad 
	\bbP\left(X \leq -x\right) \leq \exp\left(-\frac{x^2}{2\alpha}\right), \quad \forall x \geq 0.
	\end{split}
	\end{equation*}	

We can prove the following lower bound on the tail probability of gamma distribution that matches the upper bound. Since the density of $\text{Gamma}(\alpha)$ has distinct shapes for $\alpha \geq 1$ and $\alpha<1$: 
$\lim_{y \to 0} f(y;\alpha) =\infty$ if $\alpha<1$; $\lim_{y \to 0} f(y;\alpha) <\infty$ if $\alpha  \geq 1$,
the tail bound behaves differently in these two cases and we discuss them separately.
\begin{Theorem}[Gamma tail bound]\label{th:Gamma tail lower bound}
	Suppose $Y\sim \text{Gamma}(\alpha)$ and $X = Y - \alpha$. 
	\begin{itemize}[leftmargin=*]
		\item There exist two uniform constants $c, C > 0$, such that for all $\alpha \geq 1$ and $x \geq 0$, 
		\begin{equation*}
		\begin{split}
		c\exp\left(-Cx\wedge \frac{x^2}{\alpha}\right) \leq \bbP\left(X \geq x\right) \leq \exp\left(-cx\wedge \frac{x^2}{\alpha}\right).
		\end{split}
		\end{equation*}
		For any $\beta > 1$, there exists $C_{\beta} > 0$ only relying on $\beta$, such that for all $\alpha \geq1$ and $0 \leq x \leq \frac{\alpha}{\beta}$, we have
		\begin{equation*}
		\begin{split}
		c\cdot\exp\left(-C_{\beta}\frac{x^2}{\alpha}\right) \leq \bbP\left(X \leq -x\right) \leq \exp\left(-\frac{x^2}{2\alpha}\right).
		\end{split}
		\end{equation*}
		For any $x > \alpha$, $\bbP(X\leq -x) = 0$.
		\item For all $0<\alpha < 1$,
		\begin{equation}\label{ineq:alpha<1-tail}
		\begin{split}
		\frac{1}{e}\cdot\frac{(\alpha + x + 1)^\alpha - (\alpha + x)^\alpha}{e^{\alpha + x}\Gamma(\alpha + 1)} \leq \bbP\left(X \geq x\right) \leq \frac{e}{e - 1}\cdot\frac{(\alpha + x + 1)^\alpha - (\alpha + x)^\alpha}{e^{\alpha + x}\Gamma(\alpha + 1)}, \quad \forall x\geq 0;
		\end{split}
		\end{equation}
		\begin{equation}\label{ineq:alpha<1-tail2}
		\begin{split}
		\frac{\left\{(\alpha - x)\vee 0\right\}^\alpha}{e\Gamma(\alpha + 1)} \leq \bbP\left(X \leq -x\right) \leq \frac{\left\{(\alpha - x)\vee 0 \right\}^\alpha}{\Gamma(\alpha + 1)},\quad \forall x \geq 0.
		\end{split}
		\end{equation}
	\end{itemize}
\end{Theorem}
The proof for $\alpha\geq 1$ case is based on Theorem \ref{thm:general 2}; the proof of $\alpha<1$ is via the direct integration and approximation of gamma density.

\subsection{Chi-square Distribution}

The Chi-square distributions form a special class of gamma distributions and are widely used in practice. Suppose $Y\sim \chi^2_k, X = Y - k$. \cite[Lemma 1]{laurent2000adaptive} introduced the following upper tail bound for Chi-square distribution: for any $x\geq 0$,
\begin{equation*}
	\begin{split}
	\bbP(X \geq x) \leq \frac{x^2}{2(k+x)+2\sqrt{k^2+2kx}},\quad \bbP(X\leq -x) \leq \frac{x^2}{4k}.
	\end{split}
\end{equation*}
Theorem \ref{th:Gamma tail lower bound} implies the following lower bound on Chi-square distribution tail probability that matches the upper bound.
\begin{Corollary}[$\chi^2$ tail bound]\label{th:chi^2}
	Suppose $Y \sim\chi^2_k$ and $X = Y - k$ for integer $k\geq 1$. There exist uniform constants $C, c>0$ and a constant $C_{\varepsilon} > 0$ that only relies on $\varepsilon$, such that
	\begin{equation}\label{ineq:chi-square-tail}
	\begin{split}
	\bbP\left(X \geq x\right) \geq c\exp\left(- Cx\wedge \frac{x^2}{k}\right), \quad \forall x>0;
	\end{split}
	\end{equation}
	\begin{equation}\label{eq44}
	\begin{split}
	\bbP\left(X \leq - x \right)\left\{\begin{array}{ll} 
	\geq c\exp\left(-\frac{C_\varepsilon x^2}{k}\right), & \forall 0< x<(1-\varepsilon)k;\\
	= 0, & x \geq k.
	\end{array}\right.
	\end{split}
	\end{equation}
\end{Corollary}

In addition to the regular Chi-square distributions, the \emph{weighted} and \emph{noncentral Chi-square} distributions (definitions are given in the forthcoming Theorems \ref{thm:weighted chi^2 distribution} and \ref{thm: noncentral chi-square}) are two important extensions. 
We establish the matching upper and lower tail bounds for weighted Chi-square distributions in Theorem \ref{thm:weighted chi^2 distribution} and noncentral Chi-square distributions in Theorem \ref{thm: noncentral chi-square}, respectively.

\begin{Theorem}[Tail bounds of weighted $\chi^2$ distribution]\label{thm:weighted chi^2 distribution}
	Suppose $Y$ is weighted Chi-square distributed, in the sense that $Y = \sum_{i =1}^k u_iZ_i^2$, where $u_1,\ldots, u_k$ are fixed non-negative values and $Z_1,\ldots, Z_k \overset{iid}{\sim}N(0, 1)$. Then the centralized random variable $X = Y - \sum_{i=1}^k u_i$ satisfies
	\begin{equation*}
	\begin{split}
	c\exp\left(-Cx^2/\|u\|_2^2\right) \leq \bbP\left(X \geq x\right) \leq \exp(-\bar{c}x^2/\|u\|_2^2), \quad \forall 0 \leq x \leq \frac{\|u\|_2^2}{\|u\|_{\infty}};
	\end{split}
	\end{equation*}
	\begin{equation*}
	\begin{split}
	\tilde{c}\exp\left(-\tilde{C}x/\|u\|_{\infty}\right) \leq \bbP\left(X \geq x\right) \leq \exp\left(-\bar{c}x/\|u\|_{\infty}\right), \quad \forall x > \frac{\|u\|_2^2}{\|u\|_{\infty}};
	\end{split}
	\end{equation*}
	\begin{equation*}
	\begin{split}
	c_1\exp\left(-C_1x^2/\|u\|_2^2\right) \leq \bbP\left(X \leq -x\right) \leq \exp(-\frac{1}{4}x^2/\|u\|_2^2), \quad \forall 0 \leq x \leq c''\frac{\|u\|_2^2}{\|u\|_{\infty}}.
	\end{split}
	\end{equation*}
\end{Theorem}
\begin{Theorem}[Noncentral $\chi^2$ tail bound]\label{thm: noncentral chi-square}
	Let $Z$ be noncentral $\chi^2$ distributed with $k$ degrees of freedom and noncentrality parameter $\lambda$, in the sense that
	\begin{equation*}
		\begin{split}
		Z = \sum_{i=1}^k Z_i^2,\quad Z_i \sim N(\mu_i, 1) \text{ independently and } \sum_{i=1}^k \mu_i^2 = \lambda.
		\end{split}
	\end{equation*}
	Then the centralized random variable $X = Z - (k + \lambda)$ satisfies
	\begin{equation}
    \begin{split}
    	c\exp\left(-C\frac{x^2}{k + 2\lambda}\right) \leq \bbP\left(X \geq x\right) \leq \exp\left(-\bar{c}\frac{x^2}{k + 2\lambda}\right), \quad \forall 0 \leq x \leq k + 2\lambda.
    \end{split}
	\end{equation}	
	\begin{equation}
	\begin{split}
	\widetilde{c}\exp\left(-\widetilde{C}x\right) \leq \bbP\left(X \geq x\right) \leq \exp\left(-\bar{c}x\right), \quad \forall x \geq k + 2\lambda.
	\end{split}
	\end{equation}
	For all $\beta > 1$, there exists $C_{\beta} > 0$ that only relies on $\beta$ and a constant $c_1 > 0$,
	\begin{equation}
	\begin{split}
	c_1\exp\left(-C_{\beta}\frac{x^2}{k + 2\lambda}\right) \leq \bbP\left(X \leq -x\right) \leq \exp\left(-\frac{1}{4}\frac{x^2}{k + 2\lambda}\right), \quad \forall 0 < x \leq \frac{k + \lambda}{\beta}.
	\end{split}
	\end{equation}
\end{Theorem}
Here, the upper bounds in Theorems \ref{thm:weighted chi^2 distribution} and \ref{thm: noncentral chi-square} were previously proved by \cite[Lemma 1]{laurent2000adaptive} and \cite[Lemma 8.1]{birge2001alternative}, respectively. We first prove the lower bounds in Theorems \ref{thm:weighted chi^2 distribution} and \ref{thm: noncentral chi-square} in Sections \ref{sec:proof-weighted-chi^2} and \ref{sec:proof-noncentral-chi^2}, respectively.

\subsection{Beta Distribution}

Beta distribution is a class of continuous distributions that commonly appear in applications. Since the beta distribution is the conjugate prior of Bernoulli, binomial, geometric, and negative binomial, it is often used as the prior distribution for proportions in Bayesian inference. 
Recently, \cite{marchal2017sub} proved that the beta distribution is $\left(\frac{1}{4(\alpha + \beta + 1)}\right)$-sub-Gaussian, in the sense that the moment generating function of $Z\sim {\rm Beta}(\alpha, \beta)$ satisfies $\bbE\left[\exp\left(t\left(Z - \frac{\alpha}{\alpha + \beta}\right)\right)\right] \leq \exp\left(\frac{t^2}{8(\alpha + \beta + 1)}\right)$ for all $t \in \bbR$. They further gave an upper bound on tail probability:
\begin{equation}\label{ineq:previous-beta}
	\begin{split}
	\bbP\left(Z \geq \frac{\alpha}{\alpha + \beta} + x\right)\vee \bbP\left(Z \leq \frac{\alpha}{\alpha + \beta} - x\right) \leq \exp\left(-2(\alpha + \beta + 1)x^2\right), \quad \forall x\geq 0.
	\end{split}
\end{equation}
Based on the tail bounds of gamma distribution, we can prove the following matching upper and lower bounds for ${\rm Beta}(\alpha, \beta)$.
\begin{Theorem}[Beta distribution tail bound]\label{thm:beta}
	Suppose $Z \sim$ Beta($\alpha, \beta$), $\alpha, \beta \geq 1$. There exists a uniform constant $c > 0$ such that 
	\begin{equation*}
	\begin{split}
	\forall 0 < x < \frac{\beta}{\alpha + \beta}, \quad \bbP\left(Z \geq \frac{\alpha}{\alpha + \beta} + x\right) \leq \left\{\begin{array}{ll}
	2\exp\left(-c\left(\frac{\beta^2x^2}{\alpha} \wedge \beta x\right)\right), & \text{if } \beta > \alpha;\\
	2\exp\left(-c\frac{\alpha^2x^2}{\beta}\right), & \text{if } \alpha \geq \beta,
	\end{array}\right.
	\end{split}
	\end{equation*}
	\begin{equation*}
	\begin{split}
	\forall 0 < x < \frac{\alpha}{\alpha + \beta}, \quad \bbP\left(Z \leq \frac{\alpha}{\alpha+\beta} - x\right)  \leq \left\{\begin{array}{ll}
	2\exp\left(-c\left(\frac{\alpha^2x^2}{\beta} \wedge \alpha x\right)\right), & \text{if } \alpha > \beta;\\
	2\exp\left(-c\frac{\beta^2x^2}{\alpha}\right), & \text{if } \beta \geq \alpha.
	\end{array}\right.
	\end{split}
	\end{equation*}
	In addition, for any $\eta>1$, there exists $C_{\eta} > 0$ only depending on $\eta$ and a uniform constant $c > 0$, such that
	\begin{equation}
	\begin{split}
	\forall 0 < x \leq \frac{\beta}{\eta(\alpha + \beta)}, \quad \bbP\left(Z \geq \frac{\alpha}{\alpha + \beta} + x\right) \geq \left\{\begin{array}{ll}
	c\exp\left(-C_{\eta}\left(\frac{\beta^2x^2}{\alpha} \wedge \beta x\right)\right), & \text{if } \beta > \alpha;\\
	c\exp\left(-C_{\eta}\frac{\alpha^2x^2}{\beta}\right), & \text{if } \alpha \geq \beta,
	\end{array}\right.
	\end{split}
	\end{equation}
	\begin{equation*}
	\begin{split}
	\forall 0 < x \leq \frac{\alpha}{\eta(\alpha + \beta)}, \quad \bbP\left(Z \leq \frac{\alpha}{\alpha+\beta} - x\right) \geq \left\{\begin{array}{ll}
	c\exp\left(-C_{\eta}\left(\frac{\alpha^2x^2}{\beta} \wedge \alpha x\right)\right), & \text{if } \alpha > \beta;\\
	c\exp\left(-C_{\eta}\frac{\beta^2x^2}{\alpha}\right), & \text{if } \beta \geq \alpha.
	\end{array}\right.
	\end{split}
	\end{equation*}
	Moreover, 
	\begin{equation*}
	\begin{split}
	& \forall x > \frac{\beta}{\alpha+\beta},\quad \bbP\left(Z \geq \frac{\alpha}{\alpha+\beta}+x\right) = 0;\quad \forall x > \frac{\alpha}{\alpha+\beta},\quad \bbP\left(Z \leq \frac{\alpha}{\alpha+\beta}-x\right) = 0.
	\end{split}
	\end{equation*}
\end{Theorem}
Theorem \eqref{thm:beta} implies that \eqref{ineq:previous-beta} may not be sharp when $\alpha \ll \beta$ or $\alpha \gg \beta$. 

\subsection{Binomial Distribution}

The following classic result of upper tail bound for binomial distribution has been introduced in \cite{arratia1989tutorial} and \cite[Pages 23-24]{boucheron2013concentration}: 
	Suppose $Y\sim \text{Bin}(k, p)$ and $X = Y - kp$ is the centralization. Then for all $0 < x < k(1 - p)$, 
	\begin{equation}\label{ineq:previous-binomial}
	\begin{split}
	\bbP\left(X \geq x\right) \leq \exp\left(-kh_p(x/k + p)\right),
	\end{split}
	\end{equation}
	where $h_u(v) = v\log(\frac{v}{u}) + (1 - v)\log(\frac{1 - v}{1 - u})$ is the Kullback-Leibler divergence between $\text{Bernoulli}(u)$ and $\text{Bernoulli}(v)$ for all $0 < u, v < 1$.
Due to the delicate form of the moment generating function of binomial distribution, it may be difficult to select appropriate values of $(\theta, t)$ in $\bbP\left(e^{tX} \geq \theta\bbE e^{tX}\right) \geq (1 - \theta)^2\frac{(\bbE e^{tX})^2}{\bbE e^{2tX}}$ if we plan to apply Paley-Zygmund inequality to prove the desired lower tail bound. 
We instead prove the following reverse Chernoff-Cram\`er bound as a key technical tool. 
\begin{Lemma}[A reverse Chernoff-Cram\`er bound]\label{th:reverse-chernoff}
	Suppose $X$ is a random variable with moment generating function $\phi(t) = \mathbb{E}\exp(tX)$ defined for $t\in T\subseteq \mathbb{R}$. Then for any $x>0$, we have
	\begin{equation*}
	\begin{split}
	& \mathbb{P}(X\geq x) \geq \sup_{\substack{\theta,\delta>1, t\geq t' \geq 0\\t\theta \in T}}\left\{\phi(t) \exp(-t\delta x) - \phi(t\theta) \exp(-t\delta\theta x) - \exp(-(t\delta - t')x)\phi(t - t')\right\}.
	\end{split}
	\end{equation*}
	\begin{equation*}
	\begin{split}
	& \mathbb{P}(X\leq -x) \geq  \sup_{\substack{\theta,\delta>1, t\geq t' \geq 0\\-t\theta \in T}}\left\{\phi(-t) \exp(-t\delta x) - \phi(-t\theta) \exp(-t\delta\theta x) - \exp(-(t\delta - t')x)\phi(t' - t)\right\}.
	\end{split}
	\end{equation*}
\end{Lemma}

Then we can prove the following lower bound on tail probability for binomial distribution.
\begin{Theorem}[Lower bounds on binomial tail probability]\label{th:Sum of i.i.d. Bernoulli tail lower bound}
	Suppose $X$ is centralized binomial distributed with parameters $(k, p)$. 
	For any $\beta>1$, there exist constants $c_\beta, C_\beta > 0$ that only rely on $\beta$, such that 
	\begin{equation}\label{ineq:binomial-upper-tail}
	\begin{split}
	\bbP(X \geq x) \left\{\begin{array}{ll}
	\geq c_\beta\exp\left(-C_\beta k h_p\left(p + \frac{x}{k}\right)\right), & \text{if } 0 \leq x \leq \frac{k(1 - p)}{\beta} \text{ and } x+kp\geq 1;\\
	= 1 - (1-p)^k, & \text{if } 0<kp+x<1.
	\end{array}\right.
	\end{split}
	\end{equation}
	\begin{equation}\label{ineq:binomial-lower-tail}
	\begin{split}
	\bbP(X \leq - x) \left\{\begin{array}{ll}
	\geq c_\beta\exp\left(-C_\beta k h_{p}\left(p - \frac{x}{k}\right)\right), & \text{if } 0 \leq x \leq \frac{kp}{\beta}, x+ k(1 - p) \geq 1;\\
	= 1- p^k, & \text{if } 0<k(1 - p)+ x<1.
	\end{array}\right.
	\end{split}
	\end{equation}
\end{Theorem}

\begin{Remark}[Comparison to previous results]\rm
	Previously, \cite{slud1977distribution} studied the normal approximation for $X\sim \text{Bin}(k, p) - kp$ and proved that $\mathbb{P}\left(X \geq x\right) \geq 1 - \Phi\left(x/\sqrt{kp(1-p)}\right)$ if either (a) $p\leq 1/4$ and $x\geq 0$, or (b) $0\leq x \leq k(1-2p)$. Here, $\Phi$ is the cumulative distribution function of standard normal distribution. However, it is unclear if the Slud's inequality universally provides the sharp lower bound: when $p$ is close to zero, say $p = 1/k$, the lower bound provided by Slud's inequality is approximately $\exp\left(-\frac{k^2}{2}\right)$ for $x = \frac{k - 1}{2}$, which does not match the classic upper bound $\left(\approx \exp\left(-k\log k\right)\right)$. By the method of types, \cite[Lemma II.2]{cover2012elements} showed that
	\begin{equation*}
		\begin{split}
		\bbP\left(X \geq x\right) \geq \frac{1}{k+1}\exp\left(-kh_p(p+\frac{x}{k})\right)
		\end{split}
	\end{equation*}
	when $0\leq x \leq k(1-p)$ and $pk + x$ is an integer. This bound is sharp up to a factor of $\frac{1}{k+1}$ in comparison to the upper bound \eqref{ineq:previous-binomial}. 
	Our Theorem \ref{th:Sum of i.i.d. Bernoulli tail lower bound} yields sharper rate than the one by the method of types for small values of $x$. In addition, our Theorem \ref{th:Sum of i.i.d. Bernoulli tail lower bound} allows $pk+x$ or $pk-x$ to be more general real values than integers.
\end{Remark}

The sum of i.i.d. Rademacher random variables is an important instance of binomial distributions that commonly appears in practice. 
The established binomial tail bounds immediately imply the following result.
\begin{Corollary}[Sum of i.i.d. Rademacher random variables]\label{cor: i.i.d. Rademacher}
	Suppose $Z_1,\ldots, Z_k$ are i.i.d. Rademacher random variables, i.e., $\bbP(Z_i = 1) = \bbP(Z_i = -1) = 1/2$. Suppose $X = Z_1+\cdots+Z_k$. Then, For any $\beta>1$, there exist constants $c_\beta, C_\beta > 0$ that only rely on $\beta$, such that 
	\begin{equation*}		
	\begin{split}
	c_\beta\exp\left(-C_\beta \frac{x^2}{k}\right) \leq \bbP(X \geq x) = \bbP(X \leq -x) \leq \exp\left(-\frac{x^2}{4k}\right), \quad \forall 0 \leq x \leq \frac{k}{\beta}.
	\end{split}
	\end{equation*}
	\begin{equation*}
    \begin{split}
    \bbP\left(X\geq x \right) = \bbP\left(X\leq -x\right) = 0, \quad \forall x > k.
    \end{split}
	\end{equation*}
\end{Corollary}

\subsection{Poisson Distribution}\label{sec:poisson-distribution}
It was shown in literature that the Poisson distribution has the Bennett-type tail upper bound (see, e.g.,\cite{boucheron2013concentration} and \cite{pollard2015good}; also see \cite{bennett1962probability} for Bennett Inequality): 
	Suppose $Y \sim \text{Poisson}(\lambda)$, $X = Y - \lambda$. Then 
	\begin{equation}\label{ineq:previous-poisson1}
	\begin{split}
	\bbP\left(X \geq x\right) \leq \exp\left(-\frac{x^2}{2\lambda}\psi_{Benn}\left(\frac{x}{\lambda}\right)\right), \quad \forall x \geq 0,
	\end{split}
	\end{equation}
	\begin{equation}\label{ineq:previous-poisson2}
	\begin{split}
	\bbP\left(X \leq -x\right) \leq \exp\left(-\frac{x^2}{2\lambda}\psi_{Benn}\left(-\frac{x}{\lambda}\right)\right),\quad \forall 1\leq x\leq \lambda.
	\end{split}
	\end{equation}
	Here, $\psi_{Benn}$ is the Bennett function defined as
	\begin{equation}\label{eq:Bennnett_function}
	\begin{split}
	\psi_{Benn}(t) = \left\{\begin{array}{ll}
	\frac{(1 + t)\log(1 + t) - t}{t^2/2}, & \text{if } t > -1, t \neq 0;\\
	1, & \text{if } t = 0.
	\end{array}\right.
	\end{split}
	\end{equation}
We can prove the following lower bound on the tail probability for Poisson distribution.
\begin{Theorem}[Lower bound on Poisson tail]\label{cor:Poisson tail lower bound}
	Suppose $X$ is centralized Poisson distributed with parameter $\lambda$. Equivalently, $X = Y - \lambda$ and $Y \sim \text{Poisson}(\lambda)$.
	Then there exist constants $c, C > 0$ such that 
	\begin{equation}\label{ineq:poisson-upper-tail}
	\begin{split}
	\bbP(X \geq x) \left\{\begin{array}{ll}
	\geq c\cdot\exp\left(-C\cdot\frac{x^2}{2\lambda}\psi_{Benn}(x/\lambda)\right), & \text{if } x \geq 0 \text{ and } x+\lambda \geq 1;\\
	= 1 - e^{-\lambda}, & \text{if } x + \lambda <1.
	\end{array}\right.
	\end{split}
	\end{equation}
	For all $\beta > 1$, there exist constants $c_{\beta}, C_{\beta} > 0$ that only rely on $\beta$, such that
	\begin{equation}\label{ineq:poisson-lower-tail}
	\begin{split}
	\bbP(X \leq - x) \geq c_{\beta}\cdot\exp\left(-C_{\beta}\frac{x^2}{2\lambda}\psi_{Benn}(x/\lambda)\cdot\right), \quad \text{if } 0 \leq x \leq \frac{\lambda}{\beta}.
	\end{split}
	\end{equation}
\end{Theorem}
\begin{Remark}\rm
	Since the Poisson distribution is discrete, similarly to the binomial distribution, we need to discuss the boundary case of $x+\lambda<1$ when analyzing the tail bound $\bbP(X\geq x)$. By comparing to the lower bound in Theorem \ref{cor:Poisson tail lower bound}, we can see the classic Bennett-type upper bound \eqref{ineq:previous-poisson1} is not sharp when $x+\lambda<1$. 
\end{Remark}

\subsection{Irwin-Hall Distribution}
Suppose $U_1,\ldots, U_k$ are i.i.d. uniformly distributed on $[0, 1]$. Then $Y = \sum_{i = 1}^{k}U_i$ satisfies \emph{Irwin-Hall distribution}. 
We have the following matching upper and lower bounds on tail probabilities of Irwin-Hall distribution. 
\begin{Corollary}[Irwin-Hall tail bound]\label{thm:Irwin-Hall Distribution}
	Suppose $Y$ follows the Irwin-Hall distribution with parameter $k$. Denote $X = Y - \frac{k}{2}$. Then for $0 \leq x \leq \frac{k}{2}$, 
	\begin{equation}
	\begin{split}
	\bbP\left(X \leq -x\right) = \bbP\left(X \geq x\right) \leq \exp\left(-k\cdot h_{\frac{1}{2}}\left(\frac{1}{2} + \frac{x}{k}\right)\right) \leq \exp\left(-\frac{x^2}{k}\right).
	\end{split}
	\end{equation}
	There also exist constants $c, c', C > 0$, such that for all $0 \leq x \leq c'k$, we have
	\begin{equation}
	\begin{split}
	\bbP\left(X \leq -x\right) = \bbP\left(X \geq x\right) \geq c\cdot \exp\left(-C\frac{x^2}{k}\right).
	\end{split}
	\end{equation}
\end{Corollary}

\section{Applications}\label{sec:application}

\subsection{Extreme Values of Random Variables}

The extreme value theory plays a crucial role in probability, statistics, and actuarial science \cite{gumbel2012statistics,smith2003statistics,philippe2001value}. The central goal is to study the distribution of extreme value of a sequence of random variables. The lower bounds on tail probabilities developed in previous sections have a direct implication to extreme value distribution for generic random variables. 
To be specific, if $X_1,\ldots, X_k$ are weighted sums of sub-Gaussian or sub-exponential random variables, we can prove the following matching upper and lower bounds for $\mathbb{E}\sup_{1\leq i \leq k}X_i$.
\begin{Theorem}[Extreme value of the sums of independent sub-Gaussians]\label{th:extreme-value-gaussian}	
	Suppose $Z_{ij} (1 \leq i \leq k, 1 \leq j \leq n)$ are centered and independent sub-Gaussian random variables satisfying
		\begin{equation*}
		\begin{split}
		\exists \text{constants } c_1, C_1 > 0, \quad \exp(c_1\alpha t^2) \leq \phi_{Z_{ij}}(t) \leq \exp(C_1\alpha t^2),\quad \forall t \in \mathbb{R}, 1\leq i \leq k, 1 \leq j \leq n.
		\end{split}
		\end{equation*}
		If $u_1,\ldots, u_n$ are non-negative values, $X_i = u_1Z_{i1} + \dots + u_nZ_{in}$ for $1 \leq i \leq k$, then there exist constants $C, c > 0$, 
		\begin{equation*}
		\begin{split}
		c\sqrt{\alpha\|u\|_2^2\log k} \leq \bbE \sup_{1 \leq i \leq k}X_i \leq C\sqrt{\alpha\|u\|_2^2\log k}.
		\end{split}
		\end{equation*}
\end{Theorem}

\begin{Theorem}[Extreme value of the sums of independent sub-exponentials]\label{th:extreme-value-exponential}
	Suppose $Z_{ij} (1 \leq i \leq k, 1 \leq j \leq n)$ are centered and independent random variables. Suppose $Z_{ij}$'s are sub-exponential in the sense that: 
	\begin{equation*}
		\begin{split}
		\exists \text{constants } c_1, C_1 > 0, \quad \exp(c_1\alpha t^2) \leq \phi_{Z_{ij}}(t) \leq \exp(C_1\alpha t^2),\quad \forall |t|<M, 1\leq i \leq k, 1 \leq j \leq n.
		\end{split}
	\end{equation*}
	Suppose $\alpha M^2 \geq c''$, where $c'' >0$ is a constant. In addition, for any $j$, there exists one $Z_{ij}$ ($1\leq i \leq k$) further satisfying 
	\begin{equation*}
	\begin{split}
	\bbP\left(Z_{ij} \geq x \right) \geq c_0\exp(-C_0Mx), \quad \forall x \geq 2c'M\alpha,
	\end{split}
	\end{equation*}
	and $u_i \geq c_2\|u\|_{\infty}$, where $c_2 > 0$ is a constant.
	If $u_1,\ldots, u_n$ are non-negative values, $X_i = u_1Z_{i1} + \dots + u_nZ_{in}$ for $1 \leq i \leq k$, then there exist constants $C, c > 0$, 
	\begin{equation*}
		\begin{split}
		c\left[\sqrt{\alpha\|u\|_2^2\log k} \vee \frac{\|u\|_{\infty}\log k}{M}\right] \leq \bbE \sup_{1 \leq i \leq k}X_i \leq C\left[\sqrt{\alpha\|u\|_2^2\log k} \vee \frac{\|u\|_{\infty}\log k}{M}\right].
		\end{split}
	\end{equation*}
\end{Theorem}

\begin{Remark}
	If $X_1,\ldots, X_k\overset{iid}{\sim} \chi^2_n$, Theorem \ref{th:extreme-value-exponential} immediately implies
	\begin{equation*}
		\begin{split}
		c\left(\sqrt{n\log k} \vee \log k\right)\leq \mathbb{E}\sup_{1\leq i \leq k} X_i - n \leq C\left(\sqrt{n\log k} \vee \log k\right).
		\end{split}
	\end{equation*}
\end{Remark}

\subsection{Signal Identification from Sparse Heterogeneous Mixtures}

In this subsection, we consider an application of signal identification from sparse heterogeneous Mixtures. Motivated by applications in astrophysical source and genomics signal detections, the sparse heterogeneous mixture model has been proposed and extensively studied in recent high-dimensional statistics literature (see, e.g. \cite{tony2011optimal,cai2007estimation,donoho2004higher,jin2008proportion}). Suppose one observes $Y_1, \dots, Y_k\in \mathbb{R}$. $Z_1, \dots, Z_k \in \{0, 1\}$ are hidden labels that indicate whether the observations $Y_i$ are signal or noise
\begin{equation}\label{eq:Y_i}
\begin{split}
Y_i \sim P_0 1_{\{Z_i = 0\}} + P_1 1_{\{Z_i=1\}}.
\end{split}
\end{equation}
Here, $P_0$ and $P_1$ are the noise and signal distributions, respectively. Suppose $0<\varepsilon<1$, $Z_1,\ldots, Z_k$ independently satisfy
\begin{equation}\label{eq101}
\begin{split}
\bbP(Z_i = 1) = \varepsilon, \quad \bbP(Z_i =0) = 1 - \varepsilon.
\end{split}
\end{equation}
We aim to identify the set of signals, i.e., $I = \{i: Z_i = 1\}$, based on observations $Y_1,\ldots, Y_k$. 
When the observations of $Y_1,\ldots, Y_k$ are discrete count valued, it is more natural to model the observations as Poisson random variables as opposed to the more-commonly-studied Gaussian distributions in the literature. In particular, we consider
\begin{equation}\label{eq96}
\begin{split}
& Y_i \overset{iid}{\sim} \text{Poisson}(\mu)1_{\{Z_i =0\}} + \text{Poisson}(\lambda)1_{\{Z_i = 1\}}, \quad i=1,\ldots, k,
\end{split}
\end{equation}
where $\mu$ and $\lambda$ are the Poisson intensities for noisy and signal observations, respectively. To quantify the performance of any identifier $\{\hat{Z}_i\}_{i=1}^k \in \{0, 1\}^k$, we consider the following misidentification rate in Hamming distance,
\begin{equation*}
\begin{split}
M(\hat{Z}) \triangleq \frac{1}{k}\sum_{i = 1}^{k}\left|\hat{Z}_i - Z_i\right|.
\end{split}
\end{equation*}
Based on the upper and lower tail bounds of Poisson distributions in Section \ref{sec:poisson-distribution}, we can establish the following sharp bounds on misidentification rate.
\begin{Theorem}\label{thm: signal}
	Suppose $(Y_1, Z_1), \dots, (Y_k, Z_k)$ are i.i.d. pairs of observations and hidden labels that satisfy \eqref{eq101} and \eqref{eq96}. Assume $1 \leq \mu < \lambda$, $C_1\mu \leq \lambda \leq C_2\mu$ for constants $C_2 \geq C_1 > 1$, and $0 < \varepsilon < 1$. Then for any identification procedure $\hat{Z} \in \{0, 1\}^k$ based on $\{Y_i\}_{i=1}^k$, there exist two constants $C, c > 0$ such that
	\begin{equation*}
	\begin{split}
	\bbE M(\hat{Z}) \geq \left\{\begin{array}{ll}
	c\varepsilon, & 0 < \varepsilon < \varepsilon^-;\\
	c\cdot e^{-Cg(\widetilde\theta)}, & \varepsilon^- \leq \varepsilon \leq \varepsilon^+;\\
	c(1 - \varepsilon), & \varepsilon^+ <\varepsilon < 1.
	\end{array}\right.
	\end{split}
	\end{equation*}
	Here,
	\begin{equation*}
		\begin{split}
		\varepsilon^+ = \frac{1}{\exp\left(\mu\log\left(\frac{\lambda}{\mu}\right) + \mu - \lambda\right) + 1}, \varepsilon^- = \frac{1}{\exp\left(\lambda\log\left(\frac{\lambda}{\mu}\right) + \mu - \lambda\right) + 1}, \widetilde{\theta} = \frac{\log\left(\frac{1 - \varepsilon}{\varepsilon}\right) + \lambda - \mu}{\log\left(\frac{\lambda}{\mu}\right)}, 
		\end{split}
	\end{equation*}
    \begin{equation*}
    	\begin{split}
    	g(\theta) = -\log\left[(1 - \varepsilon)\exp\left(\frac{(\widetilde\theta - \mu)^2}{2\mu}\psi_{Benn}\left(\frac{\widetilde\theta - \mu}{\mu}\right)\right) + \varepsilon\exp\left(-\frac{(\lambda - \widetilde\theta)^2}{2\lambda}\psi_{Benn}\left(\frac{\widetilde\theta - \lambda}{\lambda}\right)\right)\right].
    	\end{split}
    \end{equation*}
	$\psi_{Benn}$ is the Bennett function defined in \eqref{eq:Bennnett_function}. In particular, the classifier $\widetilde{Z} = \{\widetilde{Z}_i\}_{i=1}^k, \widetilde{Z}_i = 1_{\{Y_i > \widetilde{\theta}\}}$ achieves the following misclassification rate in Hamming distance:
	\begin{equation*}
	\begin{split}
	\bbE M(\widetilde{Z}) \leq \left\{\begin{array}{ll}
	\varepsilon, & 0 < \varepsilon < \varepsilon^-;\\
	e^{-g(\widetilde{\theta})}, & \varepsilon^- \leq \varepsilon \leq \varepsilon^+;\\
	(1 - \varepsilon), & \varepsilon^+ <\varepsilon < 1.
	\end{array}\right.
	\end{split}
	\end{equation*}
\end{Theorem}
\begin{Remark}\rm 
	If the noise and signal distributions in Model \eqref{eq:Y_i}, i.e., $P_0$ and $P_1$, are other than Poisson, similar results to Theorem \ref{thm: signal} can be established based on the corresponding upper and lower tail bounds in the previous sections.
\end{Remark}

\section{Additional Proofs of Main Results}\label{sec:proof}

In this supplement we complete the proofs of the main results in this paper.

\subsection{Proof of Lemma \ref{th:reverse-chernoff}}
First, we claim that if $0 < t_2 \in T$ and $0 < t_1 < t_2$, then $t_1 \in T$. Actually, by Jensen's inequality, 
	\begin{equation*}
	\begin{split}
	\phi(t_1) = \bbE e^{t_1X} \leq \left(\bbE e^{t_2 X}\right)^{t_1/t_2} < +\infty.
	\end{split}
	\end{equation*}
	Thus $t_1 \in T$. \\
	For any $x, t>0$, $\delta, \theta>1$, we have
	\begin{equation*}
	\begin{split}
	\bbP\left(X \geq x\right) \geq & \bbP\left(x \leq X \leq \delta x\right) \geq \frac{\mathbb{E} \exp(tX) 1_{\{x\leq X\leq \delta x\}}}{\exp(t\cdot \delta x)}.
	\end{split}
	\end{equation*}
	Here, for $0 < t' \leq t$ and $t\theta \in T$,
	\begin{equation*}
	\begin{split}
	& \mathbb{E} \exp(tX)1_{\{x \leq X \leq \delta x\}}  =  \mathbb{E} \exp(tX) - \mathbb{E}\exp(tX)1_{\{X > \delta x\}} - \mathbb{E}\exp(tX)1_{\{X< x\}}\\
	\geq & \phi(t) - \mathbb{E}\exp\left(t\theta X - t\theta \delta x + t\delta x\right)1_{\{X > \delta x\}} - \mathbb{E}\exp(tX)\frac{e^{t'x}}{e^{t'X}}\\
	\geq & \phi(t) - \mathbb{E}\exp\left(t\theta X\right) \cdot \exp\left(-t\delta x(\theta-1)\right) - \exp(t'x)\bbE \exp((t - t')X)\\
	= & \phi(t) - \phi(t\theta) \cdot \exp\left(-t\delta x(\theta-1)\right) - \exp(t'x)\phi(t - t').
	\end{split}
	\end{equation*}
	The second line is because $t\theta X - t\theta \delta x + t\delta x - tX = t(\theta - 1)(X - \delta x) \geq 0$ if $X > \delta x$.\\
	Therefore, for any $x, t > 0$, $\theta, \delta >1$, and $0 < t' \leq t, t\theta \in T$, we have
	\begin{equation*}
	\begin{split}
	\mathbb{P} (X \geq x) \geq \phi(t) \exp(-t\delta x) - \phi(t\theta) \exp(-t\delta\theta x) - \exp(-(t\delta - t')x)\phi(t - t').
	\end{split}
	\end{equation*}
	By taking supremum, we have
	\begin{equation*}
	\begin{split}
	& \mathbb{P}(X\geq x) \\
	\geq & \sup_{\substack{\theta,\delta>1, t\geq t' \geq 0\\t\theta \in T}}\left\{\phi(t) \exp(-t\delta x) - \phi(t\theta) \exp(-t\delta\theta x) - \exp(-(t\delta - t')x)\phi(t - t')\right\}.
	\end{split}
	\end{equation*}
	Particularly, set $t' = t$, we have
	\begin{equation}\label{eq71}
	\begin{split}
	\mathbb{P}(X\geq x) \geq \sup_{\substack{\theta,\delta>1, t\geq t' \geq 0\\t\theta \in T}}\left\{\phi(t) \exp(-t\delta x) - \phi(t\theta) \exp(-t\delta\theta x) - \exp(-t(\delta - 1)x)\right\}.
	\end{split}
	\end{equation}
	By symmetric argument, we can also show
	\begin{equation*}
	\begin{split}
	& \bbP\left(X \leq -x\right) \\
	\geq & \sup_{\substack{\theta,\delta>1, t\geq t' \geq 0\\-t\theta \in T}}\left\{\phi(-t) \exp(-t\delta x) - \phi(-t\theta) \exp(-t\delta\theta x) - \exp(-(t\delta - t')x)\phi(t' - t)\right\}.
	\end{split}
	\end{equation*}
$\square$

\subsection{Proof of Theorem \ref{thm: general}}\label{sec:thm-general}
First, the upper bounds follow from the Chernoff-Cram\`er bound,
\begin{equation*}
\begin{split}
\bbP\left(X\geq x\right) \leq \inf_{0 \leq t \leq M}\phi_{X}(t)\exp\left(-tx\right) \leq \inf_{0 \leq t \leq M}C_2\exp\left(C_1\alpha t^2 - tx\right).
\end{split}
\end{equation*}
For all $0 \leq x \leq 2C_1M\alpha$,
\begin{equation*}
\begin{split}
t^* = \argmin_{0 \leq t \leq M}\left(C_1\alpha t^2 - tx\right) = \frac{x}{2C_1\alpha}.
\end{split}
\end{equation*}
Thus
\begin{equation*}
\begin{split}
\bbP\left(X \geq x\right) \leq C_2\exp\left(-\frac{1}{4C_1}\frac{x^2}{\alpha}\right), \quad \forall 0 \leq x \leq 2C_1M\alpha.
\end{split}
\end{equation*}
To prove the lower bound, we discuss in two scenarios: $\alpha M^2 \geq \frac{16\left(1 + \log\left(\frac{1}{c_2}\right)\right)}{c_1}$ and $\alpha M^2 \geq c'', c_2 = 1$, respectively as follows.
\begin{itemize}[leftmargin=*]
	\item
	\bm{$\alpha M^2 \geq \frac{16\left(1 + \log\left(\frac{1}{c_2}\right)\right)}{c_1}$}. Set $t = 0$, we know that $c_2 \leq 1 \leq C_2$. By Paley-Zygmund inequality, for any $0 \leq t \leq \frac{M}{2}$,
	\begin{equation}\label{eq126}
		\begin{split}
		&\bbP\left(X \geq c_1\alpha t - \frac{1 - \log(c_2)}{t}\right) = \bbP\left(e^{tX} \geq e^{-1}\cdot c_2e^{c_1\alpha t^2}\right) \geq \bbP\left(e^{tX} \geq e^{-1}\bbE e^{tX}\right)\\ \geq& \left(1 - e^{-1}\right)^2\frac{\left(c_2\exp\left(c_1\alpha t^2\right)\right)^2}{C_2\exp\left(4C_1\alpha t^2\right)} = \left(1 - e^{-1}\right)^2\frac{c_2^2}{C_2}\exp\left(-2(2C_1 - c_1)\alpha t^2\right).
		\end{split}
	\end{equation}
	By solving the equation
	\begin{equation}\label{eq129}
		\begin{split}
		c_1\alpha t - \frac{1 - \log(c_2)}{t} = x,
		\end{split}
	\end{equation}
	we have
	\begin{equation*}
		\begin{split}
		t = \frac{x + \sqrt{x^2 + 4c_1\alpha\cdot\left(1 + \log\left(\frac{1}{c_2}\right)\right)}}{2c_1\alpha} \leq \frac{2x + \sqrt{4c_1\alpha\cdot\left(1 + \log\left(\frac{1}{c_2}\right)\right)}}{2c_1\alpha}.
		\end{split}
	\end{equation*}
	If $\alpha M^2 \geq \frac{16\left(1 + \log\left(\frac{1}{c_2}\right)\right)}{c_1}$, then for any $0 \leq x \leq \frac{c_1\alpha M}{4}$,
	\begin{equation}\label{eq130}
		\begin{split}
		t \leq \frac{M}{4} + \frac{2\sqrt{c_1\alpha\cdot\left(1 + \log\left(\frac{1}{c_2}\right)\right)}}{2c_1\alpha} \leq \frac{M}{2}.
		\end{split}
	\end{equation}
	Moreover, for any $\sqrt{c_1\alpha\cdot\left(1 + \log\left(\frac{1}{c_2}\right)\right)} \leq x \leq \frac{c_1\alpha M}{4}$,
	\begin{equation*}
		\begin{split}
		0 \leq t \leq \frac{4x}{2c_1\alpha} = 2\frac{x}{c_1\alpha}.
		\end{split}
	\end{equation*}
	\eqref{eq126}, \eqref{eq129}, \eqref{eq130} and the previous inequality together imply that for any\\ $\sqrt{c_1\alpha\cdot\left(1 + \log\left(\frac{1}{c_2}\right)\right)} \leq x \leq \frac{c_1\alpha M}{4}$,
	\begin{equation*}
		\begin{split}
		\bbP\left(X \geq x\right) \geq \widetilde{c}\exp\left(-C\frac{x^2}{\alpha}\right),
		\end{split}
	\end{equation*}
	where $\widetilde{c} = \left(1 - e^{-1}\right)^2\frac{c_2^2}{C_2}, C = \frac{8(2C_1 - c_1)}{c_1^2}$. \\
	For any $0 \leq x \leq \sqrt{c_1\alpha\cdot\left(1 + \log\left(\frac{1}{c_2}\right)\right)}$, by the previous inequality,
	\begin{equation*}
		\begin{split}
		\bbP\left(X \geq x\right) \geq \bbP\left(X \geq \sqrt{c_1\alpha\cdot\left(1 + \log\left(\frac{1}{c_2}\right)\right)}\right) \geq \widetilde{c}\exp\left(-Cc_1\left(1 + \log\left(\frac{1}{c_2}\right)\right)\right).
		\end{split}
	\end{equation*}
	Set $c = \widetilde{c}\exp\left(-Cc_1\left(1 + \log\left(\frac{1}{c_2}\right)\right)\right)$, then for any $0 \leq x \leq \frac{c_1\alpha M}{4}$,
	\begin{equation*}
		\begin{split}
		\bbP\left(X \geq x\right) \geq c\exp\left(-C\frac{x^2}{\alpha}\right).
		\end{split}
	\end{equation*}
	\item
	\bm{$\alpha M^2 \geq c'', c_2 = 1$}. Similarly to \eqref{eq126}, for fixed $\lambda > 0$ and any $0 \leq t \leq \frac{M}{2}$,
	\begin{equation*}
		\begin{split}
		&\bbP\left(X \geq c_1\alpha t - \frac{\lambda}{t}\right) \geq \left(1 - e^{-\lambda}\right)^2\frac{1}{C_2}\exp\left(-2(2C_1 - c_1)\alpha t^2\right).
		\end{split}
	\end{equation*}
	By solving the equation
	\begin{equation*}
		\begin{split}
		c_1\alpha t - \frac{\lambda}{t} = x,
		\end{split}
	\end{equation*}
	we have
	\begin{equation*}
		\begin{split}
		t = \frac{x + \sqrt{x^2 + 4c_1\alpha\lambda}}{2c_1\alpha} \leq \frac{2x + \sqrt{4c_1\alpha\lambda}}{2c_1\alpha}.
		\end{split}
	\end{equation*}
	Set $\lambda = \frac{c_1c''}{16}$, for any $0 \leq x \leq \frac{c_1\alpha M}{4}$, we can check that $0 \leq t \leq \frac{M}{2}$. By the same method as the proof of the first scenario, we know that there exist constants $C, c > 0$, for any $0 \leq x \leq \frac{c_1\alpha M}{4}$,
	\begin{equation*}
	\begin{split}
	\bbP\left(X \geq x\right) \geq c\exp\left(-C\frac{x^2}{\alpha}\right).
	\end{split}
	\end{equation*}
\end{itemize}

In summary, there exist constants $c, c', C > 0$, for all $\alpha \geq 1$ and $0 \leq x \leq c'M\alpha$,
\begin{equation}\label{eq25}
\bbP\left(X \geq x\right) \geq c\exp\left(-C\frac{x^2}{\alpha}\right),
\end{equation}

Particularly, if \eqref{eq32} holds for all $t \geq 0$ with constants $C_1 \geq c_1 > 0, C_2, c_2 > 0$, set $M = \infty$ in \eqref{eq25}, we know that \eqref{eq34} holds for all $\alpha > 0$ and $x \geq 0$.\quad $\square$

\subsection{Proof of Theorem \ref{thm:general 2}}\label{sec:thm-general_2}
Due to the independence of $Y$ and $Z$,
	\begin{equation*}
	\begin{split}
	\bbP\left(X \geq x\right) \geq& \bbP\left(Y \geq 2x\right)\cdot \bbP\left(Z \geq -x\right)
	\geq T\left(2x\right)\cdot\left(1 - \bbP\left(Z \leq -x\right)\right).
	\end{split}
	\end{equation*}
	For all $\lambda \in [0, M]$,
	\begin{equation}\label{eq31}
	\begin{split}
	\bbP\left(-Z \geq z\right) \leq e^{-\lambda z}\bbE e^{\lambda(-Z)} \leq e^{-\lambda z + C_1\alpha\lambda^2}.
	\end{split}
	\end{equation}
	For any $0 \leq z \leq 2MC_1\alpha$, set $0 \leq \lambda = \frac{z}{2C_1\alpha} \leq M$ in \eqref{eq31}, we have $\bbP\left(-Z \geq z\right) \leq e^{-\frac{z^2}{4C_1\alpha}}$.
	For any $z > 2MC_1\alpha$, set $\lambda = M$ in \eqref{eq31}, we can get
	\begin{equation*}
	\begin{split}
	\bbP\left(-Z \geq z\right) \leq e^{-Mz + M^2C_1\alpha} \leq e^{-Mz + \frac{Mz}{2}} = e^{-\frac{Mz}{2}}.
	\end{split}
	\end{equation*}
	In summary,
	\begin{equation}\label{eq74}
	\begin{split}
	\bbP\left(-Z \geq z\right) \leq \left\{\begin{array}{l}
	e^{-\frac{z^2}{4C_1\alpha}}, \quad 0 \leq z \leq 2MC_1\alpha,\\
	e^{-\frac{Mz}{2}}, \quad z > 2MC_1\alpha.
	\end{array}\right.  
	\end{split}
	\end{equation}
	Set $z = x \geq \frac{w\alpha}{2}$, we have
	\begin{equation*}
	\begin{split}
	\bbP\left(Z \leq -x\right) \leq \exp\left(-\min\left\{\frac{w^2}{16C_1}, \frac{Mw}{4}\right\}\alpha\right).
	\end{split}
	\end{equation*}
	Therefore
	\begin{equation*}
	\begin{split}
	\bbP(X \geq x) \geq \left(1 - \exp\left(-\min\left\{\frac{w^2}{16C_1}, \frac{Mw}{4}\right\}\alpha\right)\right)T(2x), \quad \forall x \geq \frac{w\alpha}{2}.
	\end{split}
	\end{equation*}
$\square$

\subsection{Proof of Corollary \ref{cor:general 2}}
	Since $\bbE e^{-t(Z_2 + \cdots + Z_k)} \leq e^{C_1(k - 1)t^2} \leq e^{C_1kt^2}$, by setting $Y = Z_1, Z = Z_2 + \cdots + Z_k, \alpha = k$, Theorem \ref{thm:general 2} implies our assertion.
\quad $\square$

\subsection{Proof of Theorem \ref{cor: sub-Gaussian}}\label{sec:cor: sub-Gaussian2}
The proof of this theorem relies on the following fact. 
\begin{Lemma}\label{lm:imply}
	In Theorem \ref{cor: sub-Gaussian}, the second statement \eqref{ineq:sub-gaussian-state2} implies the first one \eqref{ineq:sub-gaussian-state1}. 
\end{Lemma}
The proof of Lemma \ref{lm:imply} follows \cite{vershynin2007random} and is also provided in Section \ref{sec:proof-imply}. 
For the upper bound, if $\phi_{Z_i}(t) \leq \exp(C_1t^2)$ for all t, by applying Chernoff-Cram\`er bound, we have
\begin{equation*}
\begin{split}
\bbP\left(X \geq x \right) \leq \min_{t \geq 0}e^{-tx}\prod_{i = 1}^{k}\phi_{Z_i}(x) \leq \min_{t \geq 0}e^{-tx}e^{C_1\|u\|_2^2t^2} = e^{-\frac{x^2}{4C_1\|u\|_2^2}}.
\end{split}
\end{equation*}
For the lower bound, we have $\exp\left(c_1\|u\|_2^2t^2\right)\leq \phi_{X}(t) \leq \exp\left(C_1\|u\|_2^2t^2\right)$ for all $t \geq 0$ and $u_1, \dots, u_k \geq 0$.
Set $\alpha = \|u\|_2^2$ in Theorem \ref{thm: general}, we know that there exist constants $c, C > 0$, for all $x \geq 0$,
\begin{equation*}
\begin{split}
\bbP\left(X \geq x\right) \geq c\exp\left(-Cx^2/\|u\|_2^2\right).
\end{split}
\end{equation*}
$\square$

\subsection{Proof of Theorem \ref{cor: sub-exponential}}\label{sec:cor: sub-exponential}
	We first consider the lower bound. Without loss of generality, we can assume that $\|u\|_{\infty} = 1$. The moment generating function of $X$ satisfies
	\begin{equation}
	\begin{split}
	\exp\left(c_1\alpha\|u\|_2^2t^2\right) \leq \phi_{X}(t) \leq \exp\left(C_1\alpha\|u\|_2^2t^2\right), \quad \forall 0 < t \leq \frac{M}{\|u\|_{\infty}} = M.
	\end{split}
	\end{equation}
	Since $\alpha M^2 \geq c''$, $\alpha\|u\|_2^2\cdot M^2 \geq \alpha\|u\|_{\infty}^2\cdot M^2 = \alpha M^2 \geq c''$. By Theorem \ref{thm: general}, there exist constants $c, c', C > 0$ such that
	\begin{equation*}
	\begin{split}
	\bbP\left(X \geq x\right) \geq c\exp\left(-C\frac{x^2}{\alpha\|u\|_2^2}\right), \quad \forall 0 \leq x \leq c'M\alpha\|u\|_2^2.
	\end{split}
	\end{equation*}
	Note that $\frac{\|u\|_2^2}{u_i} \geq \frac{\|u\|_{\infty}^2}{u_i} \geq \|u\|_{\infty} = 1$, if 
	$\bbP\left(Z_i \geq x\right) \geq c_0\exp\left(-C_0Mx\right), \forall x \geq 2c'M\alpha,$
	then
	\begin{equation}\label{eq122}
	\begin{split}
	\bbP\left(Z_i \geq x\right) \geq c_0\exp\left(-C_0Mx\right), \quad \forall x \geq 2c'M\alpha\frac{\|u\|_2^2}{u_i}.	
	\end{split}
	\end{equation}
	Moreover, the moment generating function of $\sum_{1 \leq j \leq k, j \neq i}u_jZ_j$ satisfies
	\begin{equation}
	\begin{split}
	\bbE e^{t[\sum_{1 \leq j \leq k, j \neq i}u_jZ_j]} \leq e^{C_1\alpha\sum_{1 \leq j \leq k, j \neq i}^ku_j^2t^2} \leq e^{C_1\alpha\|u\|_2^2t^2}, \quad \forall -M = -\frac{M}{\|u\|_{\infty}} \leq t \leq 0.
	\end{split}
	\end{equation}
	Set $Y = u_iZ_i, Z = \sum_{1 \leq j \leq k, j \neq i}u_jZ_j$ in Theorem \ref{thm:general 2}, for all $x > c'M\alpha\|u\|_2^2$, 
	\begin{equation}
	\begin{split}
	\bbP\left(X \geq x\right) \geq& \left(1 - e^{-c_4M^2\alpha\|u\|_2^2}\right)\cdot \bbP\left(u_i Z_i \geq 2x\right) \geq \left(1 - e^{-c_4M^2\alpha}\right)\cdot \bbP\left(Z_1 \geq \frac{2}{u_i}x\right)\\ \overset{\eqref{eq122}}{\geq}& \left(1 - e^{-c_4c''}\right)c_0\cdot\exp\left(-\frac{2}{u_i}C_0Mx\right) \geq \bar{c}\exp(-\bar{C}Mx).
	\end{split}
	\end{equation}
	Here $c_4 > 0$ is a constant, $\bar{c} = 1 - e^{-c_4c''}, \bar{C} = \frac{2C_0}{c_2}$.
	In summary, we have proved the lower bound.
	\ \par
	For the upper bound, notice that
	\begin{equation*}
	\begin{split}
	\phi_{X}(t) \leq \exp\left(C_1\alpha\|u\|_2^2t^2\right), \quad \forall 0 \leq t \leq \frac{M}{\|u\|_{\infty}},
	\end{split}
	\end{equation*}
	Similarly to \eqref{eq74}, we have
	\begin{equation}\label{eq121}
	\begin{split}
	\bbP\left(X \geq x\right) \leq \left\{\begin{array}{l}
	e^{-\frac{x^2}{4C_1\alpha\|u\|_2^2}}, \quad 0 \leq x \leq \frac{2MC_1\alpha\|u\|_2^2}{\|u\|_{\infty}};\\
	e^{-\frac{Mx}{2\|u\|_{\infty}}}, \quad x > \frac{2MC_1\alpha\|u\|_2^2}{\|u\|_{\infty}}.
	\end{array}\right.  
	\end{split}
	\end{equation}

\subsection{Proof of Theorem \ref{th:Gamma tail lower bound}}\label{sec:proof-gamma}
The moment generating function of $X$ is
\begin{equation}\label{eq123}
\begin{split}
\phi_X(t) = (1 - t)^{-\alpha}\exp(-t\alpha) = \exp\left(-\alpha(t + \log(1 - t))\right).
\end{split}
\end{equation}
First, we consider the right tail. We discuss in two scenarios: $\alpha \geq 1$ or $\alpha < 1$.
\begin{itemize}[leftmargin = *]
	\item \bm{$\alpha \geq 1$}.
	By Taylor's expansion,
	\begin{equation*}
	\begin{split}
	\frac{-t^2}{2(1 - t)} = -\sum_{i = 2}^{\infty}t^i \leq \log(1 - t) + t = -\sum_{i = 2}^{\infty}\frac{t^i}{i} \leq -\frac{t^2}{2},
	\end{split}
	\end{equation*}
	for $0 \leq t < 1$. Therefore,
	\begin{equation*}
	\begin{split}
	\exp\left(\frac{\alpha}{2}t^2\right) \leq \phi_X(t) \leq \exp\left(5\alpha\cdot t^2\right), \quad, 0 \leq t \leq \frac{9}{10}.
	\end{split}
	\end{equation*}
	By Theorem \ref{thm: general}, there exist two constants $c, c', C > 0$, for all $x \leq c'\alpha$,
	\begin{equation}\label{eq29}
	\begin{split}
	\bbP\left(X \geq x\right) \geq c\cdot\exp\left(-C\frac{x^2}{\alpha}\right).
	\end{split}
	\end{equation}
	$\forall x > c'\alpha$, let $k = \lfloor \alpha \rfloor$, and suppose $W \sim \chi_1^2$, $W_1, \dots, W_{2k}$ are i.i.d. copies of $W$, immediately we have $\frac{\sum_{i = 1}^{2k}W_i}{2} \sim \Gamma(k, 1)$. Thus
	\begin{equation}\label{eq27}
	\begin{split}
	\bbP\left(X \geq x\right) \geq \bbP\left(\sum_{i = 1}^{2k}W_i \geq 2(\alpha + x)\right) \geq \bbP\left(W_1 \geq 2(\alpha + x)\right).
	\end{split}
	\end{equation}
	Suppose $Z$ follows the standard normal distribution, note that for all $t \geq \sqrt{2}$, $\bbP\left(Z \geq t\right) \geq \frac{1}{\sqrt{2\pi}}\left(\frac{1}{t} - \frac{1}{t^3}\right)e^{-\frac{t^2}{2}} \geq \frac{1}{\sqrt{2\pi}}\frac{1}{2t}e^{-\frac{t^2}{2}}$,
	\begin{equation*}
	\begin{split}
	& \bbP\left(W_1 \geq 2(\alpha + x)\right) = 2\bbP\left(Z \geq \sqrt{2(\alpha + x)}\right)
	\geq 2\cdot\frac{1}{\sqrt{2\pi}}\frac{1}{2\sqrt{2(\alpha + x)}}e^{-(\alpha + x)}.
	\end{split}
	\end{equation*}
	Since $u < e^u$ fot all $u > 0$ and $\alpha < \frac{1}{c'}x$,
	\begin{equation*}
	\begin{split}
	\bbP\left(W_1 \geq 2(\alpha + x)\right) \geq \frac{1}{2\sqrt{\pi}}e^{-\frac{3}{2}(\alpha + x)} \geq \frac{1}{2\sqrt{\pi}}e^{-\frac{3}{2}(\frac{1}{c'} + 1)x}.
	\end{split}
	\end{equation*}
	Combine \eqref{eq27} and the previous inequality together, we conclude that for all $x \geq c'\alpha$,
	\begin{equation}\label{eq28}
	\begin{split}
	\bbP\left(X \geq x\right) \geq \frac{1}{2\sqrt{\pi}}e^{-\frac{3}{2}(\frac{1}{c'} + 1)x}.
	\end{split}
	\end{equation}
	\eqref{eq29} an \eqref{eq28} imply that for all $x \geq 0$, there exist two constants $c, C > 0$,
	\begin{equation}
	\begin{split}
	\bbP\left(X \geq x\right) \geq c\exp\left(-Cx \wedge \frac{x^2}{\alpha}\right).
	\end{split}
	\end{equation}
	\item \bm{$\alpha < 1$}. For all $y \geq 0$, since $\alpha - 1 < 0$,
	\begin{equation*}
	\begin{split}
	\int_{y}^{y + 1}t^{\alpha - 1}e^{-t}dt \geq \int_{y}^{y + 1}(t + 1)^{\alpha - 1}e^{-t}dt = \frac{1}{e}\int_{y + 1}^{y + 2}t^{\alpha - 1}e^{-t}dt.
	\end{split}
	\end{equation*}
	By induction, for all $n \in \mathbb{N}$, $\int_{y}^{y + 1}t^{\alpha - 1}e^{-t}dt \geq \frac{1}{e^n}\int_{y + n}^{y + 1 + n}t^{\alpha - 1}e^{-t}dt$. Therefore
	\begin{equation}\label{eq30}
	\begin{split}
	& \int_{\alpha + x}^{\alpha + x + 1}\frac{1}{\Gamma(\alpha)}t^{\alpha - 1}e^{-t}dt \leq \bbP\left(X \geq x\right) \\
	=& \int_{\alpha + x}^{\infty}\frac{1}{\Gamma(\alpha)}t^{\alpha - 1}e^{-t}dt = \sum_{n = 0}^{\infty}\int_{\alpha + x + n}^{\alpha + x + 1 + n}\frac{1}{\Gamma(\alpha)}t^{\alpha - 1}e^{-t}dt\\
	\leq& \sum_{n = 0}^{\infty}\frac{1}{e^n}\int_{\alpha + x}^{\alpha + x + 1}\frac{1}{\Gamma(\alpha)}t^{\alpha - 1}e^{-t}dt = \frac{e}{e - 1}\int_{\alpha + x}^{\alpha + x + 1}\frac{1}{\Gamma(\alpha)}t^{\alpha - 1}e^{-t}dt.
	\end{split}
	\end{equation}
	In order to reach the conclusion, we need to bound $\int_{\alpha + x}^{\alpha + x + 1}t^{\alpha - 1}e^{-t}dt$. Actually, note that
	\begin{equation*}
	\begin{split}
	e^{-(\alpha + x + 1)}\int_{\alpha + x}^{\alpha + x + 1}t^{\alpha - 1}dt \leq \int_{\alpha + x}^{\alpha + x + 1}t^{\alpha - 1}e^{-t}dt \leq e^{-(\alpha + x)}\int_{\alpha + x}^{\alpha + x + 1}t^{\alpha - 1}dt,
	\end{split}
	\end{equation*}
	and
	\begin{equation*}
	\begin{split}
	\int_{\alpha + x}^{\alpha + x + 1}t^{\alpha - 1}dt = \frac{(\alpha + x + 1)^\alpha - (\alpha + x)^\alpha}{\alpha},
	\end{split}
	\end{equation*}
    immediately we have
    \begin{equation*}
    	\begin{split}
    	\frac{(\alpha + x + 1)^\alpha - (\alpha + x)^\alpha}{e^{\alpha + x + 1}\Gamma(\alpha + 1)} \leq \int_{\alpha + x}^{\alpha + x + 1}t^{\alpha - 1}e^{-t}dt \leq \frac{(\alpha + x + 1)^\alpha - (\alpha + x)^\alpha}{e^{\alpha + x}\Gamma(\alpha + 1)}.
    	\end{split}
    \end{equation*}
    \eqref{eq30} and the previous inequality together imply
	\begin{equation}
	\begin{split}
	\frac{1}{e}\cdot\frac{(\alpha + x + 1)^\alpha - (\alpha + x)^\alpha}{e^{\alpha + x}\Gamma(\alpha + 1)} \leq \bbP\left(X \geq x\right) \leq \frac{e}{e - 1}\cdot\frac{(\alpha + x + 1)^\alpha - (\alpha + x)^\alpha}{e^{\alpha + x}\Gamma(\alpha + 1)}.
	\end{split}
	\end{equation}
\end{itemize}

Now, we consider the left tail, we still discuss in two scenarios: $\alpha \geq 1$ or $\alpha < 1$.
\begin{itemize}[leftmargin = *]
	\item \bm{$\alpha \geq 1$}.
	\begin{itemize}[leftmargin = *]
		\item \bm{$0 \leq x \leq c'\alpha$}.
		By \eqref{eq123}, we know that
		\begin{equation*}
		\begin{split}
		\phi_X(-t) = \exp\left(\alpha(t - \log(1 + t))\right).
		\end{split}
		\end{equation*}
		The derivative of $t - \log(1 + t)$ satisfies
		\begin{equation*}
		\begin{split}
		\frac{2t}{3} = \frac{t}{\frac{1}{2} + 1} \leq \left[t - \log(1 + t)\right]' = \frac{t}{t + 1} \leq t, \quad \forall 0 \leq t \leq \frac{1}{2},
		\end{split}
		\end{equation*}
		which means
		\begin{equation*}
		\begin{split}
		\frac{t^2}{3} \leq t - \log(1 + t) \leq \frac{t^2}{2}, \quad \forall 0 \leq t \leq \frac{1}{2}.
		\end{split}
		\end{equation*}
		Thus
		\begin{equation*}
		\begin{split}
		\exp\left(\frac{\alpha t^2}{3}\right) \leq \phi_{X}(t) \leq \exp\left(\frac{\alpha t^2}{2}\right), \quad \forall -\frac{1}{2} \leq t \leq 0.
		\end{split}
		\end{equation*}
		By Theorem \ref{thm: general}, there exist constants $c', c, C > 0$, for all $0 \leq x \leq c'\alpha$,
		\begin{equation}\label{eq45}
		\begin{split}
		\bbP\left(X \leq -x\right) = \bbP\left(-X \geq x\right) \geq c\cdot\exp\left(-C\frac{x^2}{\alpha}\right).
		\end{split}
		\end{equation}
		\item
		\bm{$c'\alpha \leq x \leq \frac{\alpha}{\beta}$}. Denote $l = \lceil\alpha\rceil$, and suppose $W \sim \chi_1^2, W_1, \dots, W_{2l}$ are i.i.d. copies of $W$, then $\frac{\sum_{i = 1}^{2l}W_i}{2} \sim \Gamma(l, 1)$.
		\begin{equation}\label{eq39}
		\begin{split}
		\bbP\left(X \leq -x\right) =& \bbP\left(Y \leq \alpha - x\right) \geq \bbP\left(\sum_{i = 1}^{2l}W_i \leq 2(\alpha - x)\right)\\ \geq& \bbP\left(\forall 1 \leq i \leq 2l, W_i \leq \frac{2(\alpha - x)}{2l}\right)\\
		\geq& \left[\bbP\left(W \leq \frac{\alpha - x}{l}\right)\right]^{2l} = \left[\bbP\left(-\sqrt{\frac{\alpha - x}{l}} \leq Z \leq \sqrt{\frac{\alpha - x}{l}}\right)\right]^{2l},
		\end{split}
		\end{equation}
		where $Z \sim N(0, 1)$.\\
		If we can show that there exists $C_{\beta}' > 0$ that only relies on $\beta$ such that for any $c'\alpha \leq x \leq \frac{\alpha}{\beta}$,
		\begin{equation}\label{eq125}
		\begin{split}
		\bbP\left(-\sqrt{\frac{\alpha - x}{l}} \leq Z \leq \sqrt{\frac{\alpha - x}{l}}\right) \geq \exp\left(-C_{\beta}'\frac{x}{2l}\right),
		\end{split}
		\end{equation}
		immediately we have
		\begin{equation}\label{eq46}
		\begin{split}
		\bbP\left(X \leq -x\right) \geq \left[\exp\left(-C_{\beta}'\frac{x}{2l}\right)\right]^{2l} = \exp\left(-C_{\beta}'x\right) \geq \exp\left(-\frac{C_{\beta}'}{c'}\frac{x^2}{\alpha}\right).
		\end{split}
		\end{equation}
		The last inequality comes from $x \geq c'\alpha$.\\
		Therefore we only need to prove $\eqref{eq125}$.
		\begin{equation}\label{eq124}
		\begin{split}
		\bbP\left(-\sqrt{\frac{\alpha - x}{l}} \leq Z \leq \sqrt{\frac{\alpha - x}{l}}\right) =& 2\int_{0}^{\sqrt{\frac{\alpha - x}{l}}}\frac{1}{\sqrt{2\pi}}e^{-\frac{x^2}{2}}dt \geq \frac{2}{\sqrt{2\pi}}\sqrt{\frac{\alpha - x}{l}}e^{-\frac{\alpha - x}{2l}}\\ \geq& \sqrt{\frac{\alpha - x}{2l}}e^{-\frac{\alpha - x}{2l}} = \left(\frac{1}{2}\frac{\alpha - x}{l}e^{-\frac{\alpha - x}{l}}\right)^{\frac{1}{2}}.
		\end{split}
		\end{equation}
		Since $c'\alpha \leq x \leq \frac{\alpha}{\beta}$ and $1 \leq \alpha \leq l$, we know that
		$0 < \frac{(\beta - 1)x}{l} \leq \frac{\alpha - x}{l} < \frac{\alpha}{l} \leq 1.$
		Also note that $\frac{t}{e^t}$ is an increasing function of $t \in [0, 1]$,
		\begin{equation*}
		\begin{split}
		\frac{\alpha - x}{l}e^{-\frac{\alpha - x}{l}} \geq \frac{\left(\beta - 1\right)x}{l}e^{-\frac{\left(\beta - 1\right)x}{l}}.	
		\end{split}
		\end{equation*}
		\eqref{eq124} and the previous inequality tell us
		\begin{equation}\label{eq127}
		\begin{split}
		\bbP\left(-\sqrt{\frac{\alpha - x}{l}} \leq Z \leq \sqrt{\frac{\alpha - x}{l}}\right) \geq \left(\frac{1}{2}\frac{\left(\beta - 1\right)x}{l}e^{-\frac{\left(\beta - 1\right)x}{l}}\right)^{\frac{1}{2}}.
		\end{split}
		\end{equation}
		By the definition of $l$ and $\alpha \geq 1$, we have $\alpha \leq l \leq 2\alpha$. Thus
		$\frac{c'}{2} = \frac{c'\alpha}{2\alpha} \leq \frac{x}{l} \leq \frac{\alpha/\beta}{\alpha} = \frac{1}{\beta}.$
		Since $\frac{y}{e^y}$ is an increasing function on $[0, 1]$ and $(\beta - 1)\frac{x}{l} \leq \frac{\beta - 1}{\beta} < 1$,
		$\frac{(\beta - 1)\frac{c'}{2}}{e^{(\beta - 1)\frac{c'}{2}}} \leq \frac{(\beta - 1)\frac{x}{l}}{e^{(\beta - 1)\frac{x}{l}}} \leq \frac{1}{e}.$
		Therefore,
		\begin{equation*}
		\begin{split}
		0 < \frac{-\log\left(\frac{\left(\beta - 1\right)x}{2l}e^{-\left(\beta - 1\right)x/l}\right)}{x/l} \leq \frac{-\log\left((\beta - 1)c'/4\cdot e^{-(\beta - 1)c'/2}\right)}{c'/2} \triangleq C_{\beta}',
		\end{split}
		\end{equation*}
		i.e.,
		\begin{equation*}
		\begin{split}
		\frac{\left(\beta - 1\right)x}{2l}e^{-\left(\beta - 1\right)x/l} \geq \exp\left(-C_{\beta}' x/l\right).
		\end{split}
		\end{equation*}
		Combine \eqref{eq127} and the previous inequality together, we know that \eqref{eq125} is true, which means we have proved \eqref{eq46}.
	
	\end{itemize}
    \ \par
    Hence, there exists a constant $c > 0$ and $C_{\beta} > 0$ that only depends on $\beta$, for all $0 \leq x \leq \frac{\alpha}{\beta}$,
    \begin{equation*}
    \begin{split}
    \bbP\left(X \leq -x\right) \geq c\cdot\exp\left(-C_{\beta}\frac{x^2}{\alpha}\right).
    \end{split}
    \end{equation*}
	\item \bm{$\alpha < 1$}.
	For all $0 \leq x \leq \alpha$,
	\begin{equation*}
	\begin{split}
	\bbP\left(X \leq -x\right) = \bbP\left(Y \leq \alpha - x\right) = \int_{0}^{\alpha - x}\frac{1}{\Gamma(\alpha)}t^{\alpha - 1}e^{-t}dt \geq \frac{1}{e}\int_{0}^{\alpha - x}\frac{1}{\Gamma(\alpha)}t^{\alpha - 1}dt = \frac{(\alpha - x)^\alpha}{e\Gamma(\alpha + 1)}.
	\end{split}
	\end{equation*}
	And
	\begin{equation*}
	\begin{split}
	\bbP\left(X \leq -x\right) \leq \int_{0}^{\alpha - x}\frac{1}{\Gamma(\alpha)}t^{\alpha - 1}dt = \frac{(\alpha - x)^\alpha}{\Gamma(\alpha + 1)}.
	\end{split}
	\end{equation*}
\end{itemize}
In summary, we have proved our assertion.\quad $\square$

\subsection{Proof of Corollary \ref{th:chi^2}}\label{sec:proof-chi^2}
For $k \geq 2$, set $\alpha = \frac{k}{2}$, then $\frac{X}{2} \sim \Gamma(\alpha, 1)$. Apply Theorem \ref{th:Gamma tail lower bound} to $\frac{X}{2}$, we know that there exist uniform constants $C', c'>0$, and a constant $C'_{\varepsilon} > 0$ depending on $\varepsilon$ only, such that
\begin{equation}\label{eq42}
\begin{split}
\bbP\left(X \geq x\right) \geq c'\exp\left(- C'x\wedge \frac{x^2}{k}\right), \quad \forall x>0, k \geq 2;
\end{split}
\end{equation}
\begin{equation}\label{eq43}
\begin{split}
\bbP\left(X \leq - x \right) \geq c'\exp\left(-\frac{C'_\varepsilon x^2}{k}\right), \quad \forall 0< x<(1-\varepsilon)k, k \geq 2.
\end{split}
\end{equation}\\
For $k = 1$, let $X_1, X_2$ be two i.i.d. copies of $X$, then $X_1 + X_2 \sim \chi_2^2$. Set $k = 2$ in \eqref{eq42} and \eqref{eq43}, we have
\begin{equation}\label{eq40}
\begin{split}
\bbP\left(X_ 1 + X_2 \geq x\right) \geq c'\exp\left(-C'x \wedge \frac{x^2}{2}\right) \geq c'\exp\left(-\frac{C'}{2}x \wedge x^2\right), \quad \forall x > 0,
\end{split}
\end{equation}
and
\begin{equation}\label{eq41}
\begin{split}
\bbP\left(X_1 + X_2 \leq - x \right) \geq c'\exp\left(-\frac{C'_\varepsilon }{2}x^2\right), \quad \forall 0< x<2(1-\varepsilon).
\end{split}
\end{equation}
The union bound shows that
\begin{equation*}
\begin{split}
\bbP\left(X_ 1 + X_2 \geq x\right) \leq \bbP\left(\exists 1 \leq i \leq 2, X_i \geq \frac{x}{2}\right) \leq 2\bbP\left(X \geq \frac{x}{2}\right).
\end{split}
\end{equation*}
Combine \eqref{eq40} with the previous inequality, we have
\begin{equation*}
\begin{split}
\bbP\left(X \geq x\right) \geq \frac{c'}{2}\exp\left(-2C'x \wedge x^2\right), \quad \forall x > 0.
\end{split}
\end{equation*}
Similarly, \eqref{eq41} and
\begin{equation*}
\begin{split}
\bbP\left(X_ 1 + X_2 \leq -x\right) \leq \bbP\left(\exists 1 \leq i \leq 2, X_i \leq -\frac{x}{2}\right) \leq 2\bbP\left(X \leq -\frac{x}{2}\right)
\end{split}
\end{equation*}
lead us to
\begin{equation*}
\begin{split}
\bbP\left(X \leq -x\right) \geq \frac{c'}{2}\exp\left(-2C'_\varepsilon x^2\right), \quad \forall 0< x<1-\varepsilon.
\end{split}
\end{equation*}
By setting $c = \frac{c'}{2}, C = 2C', C_\varepsilon = 2C'_\varepsilon$, we know that for all $k \geq 1$, \eqref{ineq:chi-square-tail} and \eqref{eq44} hold.\quad $\square$

\subsection{Proof of Theorem \ref{thm:weighted chi^2 distribution}}\label{sec:proof-weighted-chi^2}
Without loss of generality, assume $u_1 = \|u\|_{\infty} = 1$. First, we consider the right tail. Since $\frac{x^2}{\|u\|_2^2} \asymp x$ when $x \asymp \|u\|_2^2$, it suffices to show that
\begin{equation*}
\begin{split}
c\exp\left(-Cx^2/\|u\|_2^2\right) \leq \bbP\left(X \geq x\right) \leq \exp(-\bar{c}x^2/\|u\|_2^2), \quad \forall 0 \leq x \leq c'\|u\|_2^2,
\end{split}
\end{equation*}
\begin{equation*}
\begin{split}
\tilde{c}\exp\left(-\tilde{C}x\right) \leq \bbP\left(X \geq x\right) \leq \exp\left(-\bar{c}x\right), \quad \forall x > c'\|u\|_2^2,
\end{split}
\end{equation*}
where $c' > 0$ is a constant.
\begin{itemize}[leftmargin = *]
	\item \bm{$0 \leq x \leq c'\|u\|_2^2$}. For $0 < t < \frac{1}{2}$,
	\begin{equation*}
	\begin{split}
	\phi_{Z}(t) = (1 - 2t)^{-1/2}e^{-t} = e^{-\frac{1}{2}(\log(1 - 2t) + 2t)}.
	\end{split}
	\end{equation*}
	By Taylor's expansion, we know that for $0 \leq t < \frac{1}{2}$,
	\begin{equation}\label{eq79}
	2t^2 \leq -\log(1 - 2t) - 2t = \sum_{i = 1}^{\infty}\frac{(2t)^i}{i} - 2t = \sum_{i = 2}^{\infty}\frac{(2t)^i}{i} \leq \sum_{i = 2}^{\infty}\frac{(2t)^i}{2} = \frac{2t^2}{1 - 2t}.
	\end{equation}
	Thus
	\begin{equation*}
	\begin{split}
	\exp(t^2) \leq \phi_{Z}(t) \leq \exp\left(5t^2\right), \quad \forall 0 \leq t \leq \frac{2}{5},
	\end{split}
	\end{equation*}
	which means
	\begin{equation*}
	\begin{split}
	\exp\left(\|u\|_2^2t^2\right) \leq \phi_X(t) \leq \exp\left(5\|u\|_2^2t^2\right), \quad 0 \leq t \leq \frac{2}{5\|u\|_{\infty}} = \frac{2}{5}.
	\end{split}
	\end{equation*}
	Since $\|u\|_2^2 \geq \|u\|_{\infty}^2 = 1$, set $\alpha = \|u\|_2^2, M = \frac{2}{5}$ in Theorem \ref{thm: general}, we have
	\begin{equation*}
	\begin{split}
	\bbP\left(X \geq x\right) \geq c\exp\left(-C\frac{x^2}{\|u\|_2^2}\right), \quad \forall 0 \leq x \leq c'\|u\|_2^2.
	\end{split}
	\end{equation*}
	\item \bm{$x > c'\|u\|_2^2$}.
	For $x > c'\|u\|_2^2 \geq c'$, since $\int_{s}^{\infty}\frac{1}{\sqrt{2\pi}}e^{-\frac{t^2}{2}}dt \geq \left(\frac{1}{s} - \frac{1}{s^3}\right)\frac{e^{-\frac{s^2}{2}}}{\sqrt{2\pi}}$ holds for all $s > 0$,
	\begin{equation}\label{eq78}
	\begin{split}
	\bbP\left(Z_1 \geq 2x\right) =& 2\int_{\sqrt{1 + 2x}}^{\infty}\frac{1}{\sqrt{2\pi}}e^{-\frac{t^2}{2}}dt \geq 2\left(\frac{1}{\sqrt{1 + 2x}} - \frac{1}{(\sqrt{1 + 2x})^3}\right)\frac{e^{-\frac{1 + 2x}{2}}}{\sqrt{2\pi}}\\
	=& \frac{1}{\sqrt{1 + 2x}}\left(1 - \frac{1}{1 + 2x}\right)\frac{2e^{-x}}{\sqrt{2e\pi}} \geq e^{-\frac{1 + 2x}{4}}\left(1 - \frac{1}{1 + 2c'}\right)\frac{2e^{-x}}{\sqrt{2e\pi}}\\
	=& c_0\exp\left(-\frac{3}{2}x\right),
	\end{split}
	\end{equation}
	In the second inequality, we used $e^{\frac{u}{4}} \geq 1 + \frac{u}{4} \geq \sqrt{u}$ for $u \geq 0$.\\
    Moreover, Taylor's theorem tells us
	$2t - 2t^2 \leq \log (1 + 2t) \leq 2t - \frac{4t^2}{2} + \frac{8t^3}{3} = 2t - 2t^2 + \frac{8}{3}t^3$ holds for all $t \geq 0$.
	Thus
	\begin{equation}\label{eq82}
	\begin{split}
	\frac{4t^2}{3} \leq 2t - \log(1 + 2t) \leq 2t^2, \quad 0 \leq t \leq \frac{1}{4},
	\end{split}
	\end{equation}
	which means that
	\begin{equation}\label{eq77}
	\begin{split}
	\exp\left(\frac{2}{3}t^2\right) \leq \phi_{Z}(t) \leq \exp\left(t^2\right), \quad \forall -\frac{1}{4} \leq t \leq 0.
	\end{split}
	\end{equation}
	Therefore
	\begin{equation*}
	\begin{split}
	\phi_{\sum_{i = 2}^{k}u_iZ_i}(t) \leq \exp\left(\sum_{i = 2}^{k}u_i^2t^2\right) \leq \exp\left(\|u\|_2^2t^2\right), \quad -\frac{1}{4} = -\frac{1}{4\|u\|_{\infty}} \leq t \leq 0.
	\end{split}
	\end{equation*}
	By \eqref{eq78}, the previous inequality and Theorem \ref{thm:general 2}, also notice that $\|u\|_2^2 \geq 1$, for all $x \geq c'\|u\|_2^2$,
	\begin{equation*}
	\begin{split}
	\bbP\left(X \geq x\right) \geq& \left(1 - e^{-c_2\|u\|_2^2}\right)\cdot \bbP\left(u_1Z_1 \geq 2x\right) \geq \left(1 - e^{-c_2}\right)\bbP\left(Z_1 \geq 2x\right)\\ \geq& \left(1 - e^{-c_2}\right)c_0\exp\left(-\frac{3}{2}x\right) = \tilde{c}\exp\left(-\frac{3}{2}x\right),
	\end{split}
	\end{equation*}
	where $c_0, c_2 > 0$ are constants and $\tilde{c} = \left(1 - e^{-c_2}\right)c_0$.
\end{itemize}
Next, we consider the left tail. By \eqref{eq77},
\begin{equation*}
\begin{split}
\exp\left(\frac{2}{3}\|u\|_2^2t^2\right) \leq \phi_{X}(t) \leq \exp\left(\|u\|_2^2t^2\right), \quad \forall -\frac{1}{4} \leq t \leq 0.
\end{split}
\end{equation*}
By Theorem \ref{thm: general}, there exist constants $c'', c_1, C_1> 0$, for all $0 \leq x \leq c''\|u\|_2^2$,
\begin{equation*}
\begin{split}
\bbP\left(X \leq -x\right) = \bbP\left(-X \geq x\right) \geq c_1\exp\left(-C_1\frac{x^2}{\|u\|_2^2}\right).
\end{split}
\end{equation*}

Finally,  the upper tail bound can be obtained by \cite[Lemma 1]{laurent2000adaptive}: for $x \geq 0$,
\begin{equation*}
\begin{split}
\bbP\left(X \geq 2\|u\|_2\sqrt{x} + 2\|u\|_{\infty}x\right) \leq& \exp(-x);\\
\bbP\left(X \leq -2\|u\|_2\sqrt{x}\right) \leq& \exp\left(-x\right).
\end{split}
\end{equation*}
Therefore, there exists a constant $\bar{c} > 0$ such that
\begin{equation*}
\begin{split}
\bbP\left(X \geq x\right) \leq \exp\left(-\bar{c}x^2/\|u\|_2^2\right), \quad \forall 0 \leq x \leq \frac{\|u\|_2^2}{\|u\|_{\infty}}.
\end{split}
\end{equation*}
\begin{equation*}
\begin{split}
\bbP\left(X \geq x\right) \leq \exp\left(-\bar{c}x/\|u\|_{\infty}\right), \quad \forall x > \frac{\|u\|_2^2}{\|u\|_{\infty}}.
\end{split}
\end{equation*}
\begin{equation*}
\begin{split}
\bbP\left(X \geq x\right) \leq \exp\left(-\frac{1}{4}x^2/\|u\|_2^2\right), \quad \forall x \geq 0.
\end{split}
\end{equation*}
\quad $\square$

\subsection{Proof of Theorem \ref{thm: noncentral chi-square}} \label{sec:proof-noncentral-chi^2}
The moment generating function of $X$ is
\begin{equation}
\begin{split}
\phi_{X}(t) = \frac{\exp\left(\frac{\lambda t}{1 - 2t}\right)}{(1 - 2t)^{k/2}}\exp(-(k + \lambda)t) = \exp\left(\frac{2\lambda t^2}{1 - 2t} - \frac{k}{2}(\log(1 - 2t) + 2t)\right).
\end{split}
\end{equation}
First, we consider the right tail. Since $\frac{x^2}{k + 2\lambda} \asymp x$ when $x \asymp k + 2\lambda$, we only need to prove that there exists a constant $c' > 0$,
\begin{equation*}
\begin{split}
c\exp\left(-C\frac{x^2}{k + 2\lambda}\right) \leq \bbP\left(X \geq x\right) \leq \exp\left(-\bar{c}\frac{x^2}{k + 2\lambda}\right), \quad \forall 0 \leq x \leq c'(k + 2\lambda).
\end{split}
\end{equation*}	
\begin{equation*}
\begin{split}
\widetilde{c}\exp\left(-\widetilde{C}x\right) \leq \bbP\left(X \geq x\right) \leq \exp\left(-\bar{c}x\right), \quad \forall x \geq c'(k + 2\lambda).
\end{split}
\end{equation*}
\begin{itemize}[leftmargin = *]
	\item \bm{$0 \leq x \leq c'(k + 2\lambda)$}.
	By \eqref{eq79},
	\begin{equation*}
	\begin{split}
	2t^2 \leq -\log(1 - 2t) - 2t \leq \frac{2t^2}{1 - 2t}, \quad \forall 0 \leq t < \frac{1}{2}.
	\end{split}
	\end{equation*}
	Therefore
	\begin{equation}\label{eq80}
	\begin{split}
	\exp\left((k + 2\lambda)t^2\right) \leq \phi_X(t) \leq \exp\left(5(k + 2\lambda)t^2\right), \quad \forall 0 \leq t \leq \frac{2}{5}.
	\end{split}
	\end{equation}
	By Theorem \ref{thm: general}, there exist constants $c, c', C > 0$,
	\begin{equation*}
	\begin{split}
	\bbP\left(X \geq x\right) \geq c\exp\left(-C\frac{x^2}{k + 2\lambda}\right), \quad \forall 0 \leq x \leq c'(k + 2\lambda).
	\end{split}
	\end{equation*}
	\item \bm{$x > c'(k + 2\lambda)$}.
	For $k \geq 2$, suppose $W_1 \sim \chi_1^2$ and $W_2$ is a noncentral $\chi^2$ variable with $k - 1$ degrees of freedom and noncentrality parameter $\lambda$ independent with $W_1$. Then $Z$ and $W_1 + W_2$ are identically distributed. Denote $V_1 = W_1 - 1, V_2 = W_2 - (k - 1 + \lambda)$, we have $\bbE V_1 = \bbE V_2 = 0$, $X$ and $V_1 + V_2$ are identically distributed. Note that $x \geq c'(k + 2\lambda) \geq c'$, by \eqref{eq78},
	\begin{equation}\label{eq81}
	\begin{split}
	\bbP\left(V_1 \geq 2x\right)
	\geq& \frac{2}{\sqrt{2\pi}}e^{-\frac{1 + 2x}{4}}\left(1 - \frac{1}{2x + 1}\right)e^{-\frac{1 + 2x}{2}} \geq \frac{2}{\sqrt{2\pi}}e^{-\frac{1 + 2x}{4}}\cdot\left(1 - \frac{1}{2c' + 1}\right)e^{-\frac{1 + 2x}{2}}\\
	=& c_0\exp(-C_0x).
	\end{split}
	\end{equation}
	where $c_0, C_0 > 0$ are constants. \\
	By \eqref{eq80}, we have
	\begin{equation*}
	\begin{split}
	\phi_{V_2}(t) \leq \exp\left(5(k - 1+ 2\lambda)t^2\right) \leq \exp\left(5(k + 2\lambda)t^2\right), \quad \forall 0 \leq t \leq \frac{2}{5}.
	\end{split}
	\end{equation*}
	Apply Theorem \ref{thm:general 2}, $\forall x \geq c'(k  + 2\lambda)$,
	\begin{equation}
	\begin{split}
	\bbP\left(X \geq x\right) \geq \left(1 - e^{-c_2\cdot (k  + 2\lambda)}\right)\bbP\left(V_1 \geq 2x\right) \geq \left(1 - e^{-c_2}\right)c_0e^{-C_0x} = \widetilde{c}\exp\left(- \widetilde{C}x\right).
	\end{split}
	\end{equation}
	where $c_2, \widetilde{c}$, and $\widetilde{C}$ are constants. \\
	For $k = 1$,
	\begin{equation*}
	\begin{split}
	\bbP\left(X \geq x\right) =& \bbP\left(Z \geq 1 + \lambda + x\right) \geq \int_{\sqrt{1 + \lambda + x} - \sqrt{\lambda}}^{\infty}\frac{1}{\sqrt{2\pi}}e^{-\frac{u^2}{2}}du \geq \int_{\sqrt{1 + x}}^{\infty}\frac{1}{\sqrt{2\pi}}e^{-\frac{u^2}{2}}du.
	\end{split}
	\end{equation*}
	The last inequality comes from $\sqrt{1 + \lambda + x} - \sqrt{\lambda} \leq \sqrt{1 + x}$ for all $\lambda, x \geq 0$.
	Similarly to \eqref{eq78}, we know that there exist constants $\widetilde{C}, \widetilde{c} > 0$ such that for all $x \geq c'(1 + 2\lambda) \geq c'$,
	\begin{equation*}
	\begin{split}
	\bbP\left(X \geq x\right) \geq \widetilde{c}\exp\left(-\widetilde{C}x\right).
	\end{split}
	\end{equation*}
	Thus there exist two universal constants $\widetilde{C}, \widetilde{c} > 0$, for all $k \geq 1, x \geq c'(k + 2\lambda)$,
	\begin{equation*}
		\begin{split}
		\bbP\left(X \geq x\right) \geq \widetilde{c}\exp\left(-\widetilde{C}x\right).
		\end{split}
	\end{equation*}
\end{itemize}

Now, we consider the left tail.
\begin{itemize}[leftmargin = *]
	\item \bm{$0 <  x \leq c''(k + 2\lambda)$}.  Note that
	\begin{equation*}
	\begin{split}
	\frac{4}{3}t^2 \leq \frac{2t^2}{1 - 2t} \leq 2t^2, \quad \forall -\frac{1}{4} \leq t \leq 0,
	\end{split}
	\end{equation*}
	by \eqref{eq82}, we have
	\begin{equation}
	\begin{split}
	\exp\left(\frac{2}{3}(k + 2\lambda)t^2\right) \leq \phi_{X}(t) \leq \exp\left((k + 2\lambda)t^2\right), \quad \forall -\frac{1}{4} \leq t \leq 0.
	\end{split}
	\end{equation}
	By Theorem 2, there exist constants $c_1, c'', C_1 > 0$ such that
	\begin{equation*}
	\begin{split}
	\bbP\left(X \leq -x\right) \geq c_1\exp\left(-C_1\frac{x^2}{k + 2\lambda}\right)
	\end{split}
	\end{equation*}
	holds for all $0 < x \leq c''(k + 2\lambda)$.
	\item \bm{$c''(k + 2\lambda) < x \leq \frac{k + \lambda}{\beta}$}.
	Suppose $Q_1, \dots, Q_k$ are i.i.d. non-central $\chi^2$ distributed random variables with 1 degree of freedom and noncentrality parameter $\frac{\lambda}{k}$. Denote $R_i = Q_i - \left(1 + \frac{\lambda}{k}\right)$ for $1 \leq i \leq k$, immediately we have $\bbE R_i = 0$. Moreover, $X$ and $\sum_{i = 1}^{k}R_i$ are identically distributed.
	\begin{equation}\label{eq83}
	\begin{split}
	\bbP\left(X \leq -x\right) = \bbP\left(\sum_{i = 1}^{k}R_i \leq -x\right) \geq \bbP\left(\forall 1 \leq i \leq k, R_i \leq -\frac{x}{k}\right) = \left[\bbP\left(R_1 \leq -\frac{x}{k}\right)\right]^k.
	\end{split}
	\end{equation}
	Since $e^{-\frac{u^2}{2}} \geq e^{-\frac{\left(\sqrt{\frac{k + \lambda - x}{k}} + \sqrt{\frac{\lambda}{k}}\right)^2}{2}}$ for all $-\sqrt{\frac{k + \lambda - x}{k}} - \sqrt{\frac{\lambda}{k}} \leq u \leq \sqrt{\frac{k + \lambda - x}{k}} - \sqrt{\frac{\lambda}{k}}$,
	\begin{equation}\label{eq84}
	\begin{split}
	\bbP\left(R_1 \leq -\frac{x}{k}\right) = \int_{-\sqrt{\frac{k + \lambda - x}{k}} - \sqrt{\frac{\lambda}{k}}}^{\sqrt{\frac{k + \lambda - x}{k}} - \sqrt{\frac{\lambda}{k}}}\frac{1}{\sqrt{2\pi}}e^{-\frac{u^2}{2}}du \geq \frac{2}{\sqrt{2\pi}}\sqrt{\frac{k + \lambda - x}{k}}e^{-\frac{\left(\sqrt{\frac{k + \lambda - x}{k}} + \sqrt{\frac{\lambda}{k}}\right)^2}{2}}.
	\end{split}
	\end{equation}
	Notice that $e^t \geq t + 1 \geq 2\sqrt{t}$ for all $t \geq 0$,
	for all $c''(k + 2\lambda) < x \leq \frac{k + \lambda}{\beta}$,
	we have
	\begin{equation*}
	\begin{split}
	\frac{2}{\sqrt{2\pi}}\sqrt{\frac{k + \lambda - x}{k}} \geq \frac{1}{2}\sqrt{\frac{k + \lambda - x}{k}} \geq \exp\left(-\frac{k}{k + \lambda - x}\right) \geq \exp\left(-\frac{k}{\left(1 - \frac{1}{\beta}\right)(k + \lambda)}\right).
	\end{split}
	\end{equation*}
	Also use the basic inequality $x^2 + y^2 \geq \frac{(x + y)^2}{2}$ for all $x, y \in \bbR$, for all $c''(k + 2\lambda) < x \leq \frac{k + \lambda}{\beta}$,
	\begin{equation*}
	\begin{split}
	e^{-\frac{\left(\sqrt{\frac{k + \lambda - x}{k}} + \sqrt{\frac{\lambda}{k}}\right)^2}{2}} \geq e^{-\frac{k + \lambda - x}{k} - \frac{\lambda}{k}} = \exp\left(-\frac{k + 2\lambda - x}{k}\right) \geq \exp\left(-\left(\frac{1}{c''} - 1\right)\frac{x}{k}\right),
	\end{split}
	\end{equation*}
	The previous two inequalities and \eqref{eq84} imply that
	\begin{equation}\label{eq85}
	\begin{split}
	\bbP\left(R_1 \leq -\frac{x}{k}\right) \geq \exp\left(-\frac{k}{\left(1 - \frac{1}{\beta}\right)(k + \lambda)} - \left(\frac{1}{c''} - 1\right)\frac{x}{k}\right).
	\end{split}
	\end{equation}
	\eqref{eq85} and \eqref{eq83} lead us to
	\begin{equation}
	\begin{split}
	\bbP\left(X \leq -x\right) \geq& \exp\left(-\frac{k^2}{\left(1 - \frac{1}{\beta}\right)(k + \lambda)} - \left(\frac{1}{c''} - 1\right)x\right) \geq \exp\left(-\frac{k + 2\lambda}{1 - \frac{1}{\beta}} - \left(\frac{1}{c''} - 1\right)x\right)\\
	\geq& \exp\left(-\frac{x}{\left(1 - \frac{1}{\beta}\right)c''} - \left(\frac{1}{c''} - 1\right)x\right) = \exp\left(-C_2(\beta)x\right),
	\end{split}
	\end{equation}
	where $C_3(\beta) = \frac{1}{\left(1 - \frac{1}{\beta}\right)c''} + \left(\frac{1}{c''} - 1\right)$.\\
	Since $x \leq \frac{1}{c''}\frac{x^2}{k + 2\lambda}$ for $c''(k + 2\lambda) < x \leq \frac{k + \lambda}{\beta}$, let $C_4(\beta) = \frac{1}{c''}C_3(\beta)$, we have
	\begin{equation}
	\begin{split}
	\bbP\left(X \leq -x\right) \geq \exp\left(-C_4(\beta)\frac{x^2}{k + 2\lambda}\right), \quad \forall c''(k + 2\lambda) < x \leq \frac{k + \lambda}{\beta}.
	\end{split}
	\end{equation}
\end{itemize}
Set $C_{\beta} = \max\{C_4(\beta), C_1\}$, for all $0 < x \leq \frac{k + \lambda}{\beta}$,
\begin{equation*}
\begin{split}
\bbP\left(X \leq -x\right) \geq c_1\exp\left(-C_{\beta}\frac{x^2}{k + 2\lambda}\right).
\end{split}
\end{equation*}

Finally, by Lemma 8.1 in \cite{birge2001alternative},
\begin{equation*}
\begin{split}
\bbP\left(X \geq 2\sqrt{(k + 2\lambda)x} + 2x\right) \leq \exp(-x),
\end{split}
\end{equation*}
and
\begin{equation*}
\begin{split}
\bbP\left(X \leq - 2\sqrt{(k + 2\lambda)x}\right) \leq \exp\left(-x\right).
\end{split}
\end{equation*}
Thus there exists a constant $\bar{c} > 0$,
\begin{equation*}
\begin{split}
\bbP\left(X \geq x\right) \leq \exp\left(-\bar{c}\frac{x^2}{k + 2\lambda} \vee x\right), \quad \forall x \geq 0,
\end{split}
\end{equation*}
and
\begin{equation*}
\begin{split}
\bbP\left(X \leq -x\right) \leq \exp\left(-\frac{1}{4}\frac{x^2}{k + 2\lambda}\right), \quad \forall 0 \leq x \leq k + \lambda.
\end{split}
\end{equation*}
\quad $\square$

\subsection{Proof of Theorem \ref{thm:beta}}\label{sec:proof-beta}
Suppose $R_1 \sim$ $\Gamma(\alpha, 1)$, $R_2 \sim$ $\Gamma(\beta, 1)$, $R_1$ and $R_2$ are independent. We first consider the right tail.
\begin{itemize}[leftmargin = *]
	\item \textbf{Upper Bound}.
	Then $\forall x \in (0, \frac{\beta}{\alpha + \beta})$,
	\begin{equation}\label{eq62}
	\begin{split}
	\bbP\left(Z \geq \frac{\alpha}{\alpha + \beta} + x\right) = \bbP\left(\frac{R_1}{R_1 + R_2} \geq \frac{\alpha}{\alpha + \beta} + x\right).
	\end{split}
	\end{equation}
	By \eqref{ineq:previous-gamma}, 
	$\forall z > 0, a > 1$, $Y \sim$ $\Gamma(a, 1)$, we have
	\begin{equation}\label{gamma tail bound}
	\begin{split}
	\bbP\left(Y > a + \sqrt{2az} + z\right) &\leq e^{-z},\\
	\bbP\left(Y < a - \sqrt{2az}\right) &\leq e^{-z}.
	\end{split}
	\end{equation}
	$\forall x \in (0, \frac{\beta}{\alpha + \beta})$, if $\exists t \in (0, \infty)$ such that
	\begin{equation}\label{eq54}
	\begin{split}
	\frac{\alpha + \sqrt{2\alpha t} + t}{\alpha + \beta + \sqrt{2\alpha t} - \sqrt{2\beta t} + t} = \frac{\alpha}{\alpha + \beta} + x,
	\end{split}
	\end{equation}
	then immediately we have
	\begin{equation}\label{eq59}
	\begin{split}
	\bbP\left(Z \geq \frac{\alpha}{\alpha + \beta} + x\right) \leq \bbP\left(R_1 > \alpha + \sqrt{2\alpha t} + t\right) + \bbP\left(R_2 < \beta - \sqrt{2\beta t}\right) \leq 2e^{-t}.
	\end{split}
	\end{equation}
	Now, we prove that there exists $t > 0$ satisfying \eqref{eq54}. Actually, since for all $t > 0$,
	\begin{equation*}
	\begin{split}
	\alpha + \beta + \sqrt{2\alpha t} - \sqrt{2\beta t} + t > \beta + t - \sqrt{2\beta t} \geq 2\sqrt{\beta t}  - \sqrt{2\beta t} > 0,
	\end{split}
	\end{equation*}
	\eqref{eq54} is equivalent to
	\begin{equation*}
	\begin{split}
	A_xt + B_x\sqrt{t} - D_x = 0,
	\end{split}
	\end{equation*}
	where
	\begin{equation*}
	\begin{split}
	A_x = \beta - (\alpha + \beta)x, \quad B_x = \beta\sqrt{2\alpha} + \alpha\sqrt{2\beta} + (\alpha + \beta)\left(\sqrt{2\beta} - \sqrt{2\alpha}\right)x, \quad D_x = (\alpha + \beta)^2x.
	\end{split}
	\end{equation*}
	Since $A_x > 0$ and $D_x > 0$ as $0 < x < \frac{\beta}{\alpha + \beta}$, we know that $A_xu^2 + B_xu - D_x = 0$ has a positive solution $u$. Let $t = u^2$, then $x$ satisfies \eqref{eq54}.\\
	Moreover, $B_x$ is a linear function of $x$ and
	\begin{equation*}
	\begin{split}
	B_0 = \beta\sqrt{2\alpha} + \alpha\sqrt{2\beta} > 0, \quad B_{\frac{\beta}{\alpha + \beta}} = (\alpha + \beta)\sqrt{2\beta} > 0,
	\end{split}
	\end{equation*}
	therefore, for all $0 < x < \frac{\beta}{\alpha + \beta}$, $B_x > 0$.
	For convenience, let $A, B, D$ denote $A_x, B_x, D_x$, separately.
	\begin{equation}\label{eq55}
	\begin{split}
	t = u^2 = \left(\frac{-B + \sqrt{B^2 + 4AD}}{2A}\right)^2 = \left(\frac{2D}{B + \sqrt{B^2 + 4AD}}\right)^2 \asymp \min\left\{\frac{D^2}{B^2}, \frac{D}{A}\right\}.
	\end{split}
	\end{equation}
	Now we discuss in two scenarios: $\beta > \alpha$ and $\alpha \geq \beta$.
	\begin{itemize}[leftmargin = *]
		\item[(1)]
		\bm{$\beta > \alpha$}. For $0 < x <\frac{\beta}{\alpha + \beta}$,
		\begin{equation}\label{eq57}
		\begin{split}
		0 = A_{\frac{\beta}{\alpha + \beta}}< A = A_x < A_0 = \beta,
		\end{split}
		\end{equation}
		and
		\begin{equation}\label{eq56}
		\begin{split}
		\frac{D}{B} \asymp \min\left\{\frac{(\alpha + \beta)^2x}{\beta\sqrt{2\alpha} + \alpha\sqrt{2\beta}}, \frac{(\alpha + \beta)^2x}{(\alpha + \beta)\left(\sqrt{2\beta} - \sqrt{2\alpha}\right)x}\right\} \asymp \min\left\{\frac{\beta x}{\sqrt{\alpha}}, \frac{\beta\sqrt{\beta}}{\beta - \alpha}\right\}.
		\end{split}
		\end{equation}
		Combine \eqref{eq55}, \eqref{eq57} and \eqref{eq56} together,  we know that for all $\beta > \alpha$, $0 < x < \frac{\beta}{\alpha + \beta}$,
		\begin{equation*}
		\begin{split}
		t \succeq \min\left\{\frac{\beta^2 x^2}{\alpha}, \frac{\beta^3}{\left(\beta - \alpha\right)^2}, \frac{\left(\alpha + \beta\right)^2x}{\beta}\right\} \asymp \min\left\{\frac{\beta^2 x^2}{\alpha}, \frac{\beta^3}{\left(\beta - \alpha\right)^2}, \beta x\right\} = \min\{\frac{\beta^2 x^2}{\alpha}, \beta x\}.
		\end{split}
		\end{equation*}
		The last equation holds since $\beta x < \beta < \frac{\beta^3}{\left(\beta - \alpha\right)^2}$ holds for all $0 < x < \frac{\beta}{\alpha + \beta}$. \\
		Combine \eqref{eq59} with the previous inequality together, there exists a constant $c > 0$, for all $0 < x < \frac{\beta}{\alpha + \beta}$,
		\begin{equation}
		\begin{split}
		\bbP\left(Z \leq \frac{\alpha}{\alpha + \beta} + x\right) \leq 2\exp\left(-c \min\left\{\frac{\beta^2 x^2}{\alpha}, \beta x\right\}\right).
		\end{split}
		\end{equation}
		\item[(2)]\bm{$\alpha \geq \beta$}. For $0 < x < \frac{\beta}{\alpha + \beta}$,
		\begin{equation*}
		\begin{split}
		\sqrt{2}\cdot\alpha\sqrt{\beta} \leq (\alpha + \beta)\sqrt{2\beta} = B_{\frac{\beta}{\alpha + \beta}} \leq B = B_x \leq  B_0 = \beta\sqrt{2\alpha} + \alpha\sqrt{2\beta} \leq 2\sqrt{2}\cdot\alpha\sqrt{\beta}.
		\end{split}
		\end{equation*}
		Thus
		\begin{equation}\label{eq58}
		\begin{split}
		\frac{D}{B} \asymp \frac{(\alpha + \beta)^2x}{\alpha\sqrt{\beta}} \asymp \frac{\alpha x}{\sqrt{\beta}}.
		\end{split}
		\end{equation}
		By \eqref{eq55}, \eqref{eq57} and \eqref{eq58}, for $0 < x < \frac{\beta}{\alpha + \beta}$,
		\begin{equation*}
		\begin{split}
		t \succeq \min\left\{\left(\frac{\alpha x}{\sqrt{\beta}}\right)^2, \frac{(\alpha + \beta)^2x}{\beta}\right\} \asymp \min\left\{\frac{\alpha^2x^2}{\beta}, \frac{\alpha^2 x}{\beta}\right\} = \frac{\alpha^2x^2}{\beta}.
		\end{split}
		\end{equation*}
		\eqref{eq59} and the previous inequality imply that there eixsts a constant $c > 0$, for all $0 < x < \frac{\beta}{\alpha + \beta}$,
		\begin{equation*}
		\begin{split}
		\bbP\left(Z \geq \frac{\alpha}{\alpha + \beta} + x\right) \leq 2\exp\left(-c\frac{\alpha^2x^2}{\beta}\right).
		\end{split}
		\end{equation*}
		
	\end{itemize}
	In summary, there exists a constant $c > 0$, for all $0 < x < \frac{\beta}{\alpha + \beta}$,
	\begin{equation}
	\begin{split}
	\bbP\left(Z \geq \frac{\alpha}{\alpha + \beta} + x\right) \leq \left\{\begin{array}{ll}
	2\exp\left(-c\frac{\beta^2x^2}{\alpha} \wedge \beta x\right), & \text{if } \beta > \alpha;\\
	2\exp\left(-c\frac{\alpha^2x^2}{\beta}\right), & \text{if } \alpha \geq \beta.
	\end{array}\right.
	\end{split}
	\end{equation}
	\item \textbf{Lower Bound}.
	By Theorem \ref{th:Gamma tail lower bound}, for all $\eta > 1$, there exist constants $\bar{c}, \bar{C}, \bar{C}_{\eta}> \bar{C} > 0$, for all $y \geq 0$,
	\begin{equation}\label{eq61}
	\begin{split}
	\bbP\left(R_1 \geq \alpha + y\right) \geq \bar{c}\cdot\exp\left(-\bar{C}\frac{y^2}{\alpha} \wedge y\right),
	\end{split}
	\end{equation}
	and for all $0 \leq y \leq \frac{\beta}{\eta}$,
	\begin{equation}\label{eq63}
	\begin{split}
	\bbP\left(R_2 \leq \beta - y\right) \geq \bar{c}\cdot\exp\left(-\bar{C}_{\eta}\frac{y^2}{\beta}\right).
	\end{split}
	\end{equation}
	Fix $y = (\alpha + \beta)x$, then $0 \leq y \leq (\alpha + \beta)\frac{\beta}{\eta(\alpha + \beta)} = \frac{\beta}{\eta}$,
	\begin{equation*}
	\begin{split}
	\frac{\alpha + y}{\alpha + y + \beta - y} = \frac{\alpha}{\alpha + \beta} + x
	\end{split}
	\end{equation*}
	By \eqref{eq62}, \eqref{eq61}, \eqref{eq63} and the previous inequality, for all $0 < x \leq \frac{\beta}{\eta(\alpha + \beta)}$,
	\begin{equation}\label{eq64}
	\begin{split}
	&\bbP\left(Z \geq \frac{\alpha}{\alpha + \beta} + x\right) \geq \bbP\left(R_1 \geq \alpha + y, R_2 \leq \beta - y\right) = \bbP\left(R_1 \geq \alpha + y\right)\cdot\bbP\left(R_2 \leq \beta - y\right)\\
	\geq& \bar{c}\exp\left(-\bar{C}\frac{y^2}{\alpha} \wedge y\right)\cdot\bar{c}\exp\left(-\bar{C}_{\eta}\frac{y^2}{\beta}\right)
	\geq \bar{c}\exp\left(-\bar{C}_{\eta}\frac{y^2}{\alpha} \wedge y\right)\cdot\bar{c}\exp\left(-\bar{C}_{\eta}\frac{y^2}{\beta}\right).
	\end{split}
	\end{equation}
	We discuss in two scenarios: $\beta > \alpha$ or $\alpha \geq \beta$.
	\begin{itemize}[leftmargin = *]
		\item[(1)] \bm{$\beta > \alpha$}. Immediately we have$\frac{y^2}{\beta} \leq \frac{y^2}{\alpha}$. Since $0 \leq y < \beta$, we know that $y \geq \frac{y^2}{\beta}$, thus
		\begin{equation*}
		\begin{split}
		\frac{y^2}{\beta} \leq \frac{y^2}{\alpha} \wedge y.
		\end{split}
		\end{equation*}
		Combine \eqref{eq64} and the previous inequality, we conclude that
		\begin{equation}
		\begin{split}
		\bbP\left(Z \geq \frac{\alpha}{\alpha + \beta} + x\right) \geq& \bar{c}\exp\left(-\bar{C}_{\eta}\frac{y^2}{\alpha} \wedge y\right)\cdot\bar{c}\exp\left(-\bar{C}_{\eta}\frac{y^2}{\alpha} \wedge y\right)\\ =& \bar{c}^2\exp\left(-2\bar{C}_{\eta}\frac{y^2}{\alpha} \wedge y\right)
		\geq \bar{c}^2\exp\left(-8\bar{C}_{\eta}\frac{\beta^2x^2}{\alpha} \wedge \beta x\right).
		\end{split}
		\end{equation}
		\item[(2)] \bm{$\alpha \geq \beta$}. Since $0 \leq y < \beta$, we have
		\begin{equation*}
		\begin{split}
		\frac{y^2}{\alpha} \leq \frac{y^2}{\beta} < y.
		\end{split}
		\end{equation*}
		\eqref{eq64} and the previous inequality show that
		\begin{equation}
		\begin{split}
		\bbP\left(Z \geq \frac{\alpha}{\alpha + \beta} + x\right) \geq& \bar{c}\exp\left(-\bar{C}_{\eta}\frac{y^2}{\beta}\right)\cdot\bar{c}\exp\left(-\bar{C}_{\eta}\frac{y^2}{\beta}\right)\\
		=& \bar{c}^2\exp\left(-2\bar{C}_{\eta}\frac{y^2}{\beta}\right)
		\geq \bar{c}^2\exp\left(-8\bar{C}_{\eta}\frac{\alpha^2x^2}{\beta}\right).
		\end{split}
		\end{equation}
	\end{itemize}
	Thus there exists a universal constant $c > 0$ and a constant $C_{\eta} > 0$ that only depends on $\eta$, for all $0 < x \leq \frac{\beta}{\eta(\alpha + \beta)}$,
	\begin{equation}
	\begin{split}
	\bbP\left(Z \geq \frac{\alpha}{\alpha + \beta} + x\right) \geq \left\{\begin{array}{ll}
	c\exp\left(-C_{\eta}\left(\frac{\beta^2x^2}{\alpha} \wedge \beta x\right)\right), & \text{if } \beta > \alpha;\\
	c\exp\left(-C_{\eta}\frac{\alpha^2x^2}{\beta}\right), & \text{if } \alpha \geq \beta.
	\end{array}\right.
	\end{split}
	\end{equation}
\end{itemize}
Next, we consider the left tail.
\begin{equation*}
\begin{split}
\bbP\left(Z \leq \frac{\alpha}{\alpha + \beta} - x\right) = \bbP\left(1 - Z \geq \frac{\beta}{\alpha + \beta} + x\right).
\end{split}
\end{equation*}
Note that if $Z \sim \text{Beta}(\alpha, \beta)$, then $1 - Z \sim \text{Beta}(\beta, \alpha)$. Directly applying the results for the right tail, we reach the conclusion for the left tail.\quad $\square$

\subsection{Proof of Theorem \ref{th:Sum of i.i.d. Bernoulli tail lower bound}}\label{sec:proof-binomial}

First, the moment generating function of $X$ is
\begin{equation*}
\begin{split}
\phi_X(t) = \mathbb{E}\exp(tX) = \frac{\left(1 - p + pe^t\right)^k}{e^{kpt}} = \exp\left(k\log\left(1 - p + pe^t\right) - kpt\right).
\end{split}
\end{equation*}
In order to apply Lemma \ref{th:reverse-chernoff}, we aim to choose $\theta$ and $t'$ that satisfy
\begin{equation*}
\begin{split}
& t\theta = \argmin_{u}\left\{k\log\left(1 - p + pe^u\right) - kpu - u\delta x\right\}\\
\text{and}\quad &  t - t' = \argmin_{u}\left\{k\log\left(1 - p + pe^u\right) - kpu - ux\right\}.
\end{split}
\end{equation*}
By calculation, we can see the previous equations are equivalent to
\begin{equation*}
\begin{split}
t\theta = \log \left(\frac{(1 - p)(p + \frac{\delta x}{k})}{p(1 - p - \frac{\delta x}{k})}\right), \quad t - t' = \log \left(\frac{(1 - p)(p + \frac{x}{k})}{p(1 - p - \frac{x}{k})}\right).
\end{split}
\end{equation*}
Set $t = \log \left(\frac{(1 - p)(p + \frac{\delta' x}{k})}{p(1 - p - \frac{\delta' x}{k})}\right)$ for some $1 < \delta' < \delta$ to be specified later. For $0 < x < \frac{k(1 - p)}{\delta}$, since
\begin{equation}\label{ineq:p-q}
\begin{split}
(1 - p)q > (1 - q)p, \text{ for any }0 < p < q < 1,
\end{split}
\end{equation}
we set $q = p + \frac{x}{k}, p + \frac{\delta' x}{k}$, $p + \frac{\delta x}{k}$, and can conclude that $t - t' > 0, t\theta > 0$, and $t > 0$, respectively. By \eqref{ineq:p-q}, we also know that $t\theta > t > t - t'$. Thus, $\theta > 1$ and $t > t' > 0$. Next, we evaluate that
\begin{equation*}
\begin{split}
& \phi_X(t)\exp(-t\delta x) = \exp\left(k\log\left(\frac{1 - p}{1 - p - \frac{\delta' x}{k}}\right) - (kp + \delta x)\log\left(\frac{(1 - p)(p + \frac{\delta'x}{k})}{(1 - p - \frac{\delta' x}{k})p}\right)\right),\\
& \phi_X(t\theta)\exp(-t\theta\delta x) = \exp\left(k\log\left(\frac{1 - p}{1 - p - \frac{\delta x}{k}}\right) - (kp + \delta x)\log\frac{(1 - p)(p + \frac{\delta x}{k})}{(1 - p - \frac{\delta x}{k})p}\right),\\
& \phi_X(t - t')\exp(-(t\delta - t')x) = \exp\Bigg(k\log\left(\frac{1 - p}{1 - p - \frac{x}{k}}\right) - (kp + x)\log\left(\frac{(1 - p)(p + \frac{x}{k})}{(1 - p - \frac{x}{k})p}\right) \\
& \qquad\qquad\qquad\qquad\qquad\qquad\qquad - (\delta - 1)x \log \left(\frac{(1 - p)(p + \frac{\delta' x}{k})}{p(1 - p - \frac{\delta' x}{k})}\right)\Bigg).
\end{split}
\end{equation*}
\begin{equation*}
\begin{split}
& \log\left(\phi_X(t)\exp(-t\delta x)\right) - \log\left(\phi_X(t\theta)\exp(-t\theta\delta x)\right)\\
=& k\log\left(\frac{1 - p - \frac{\delta x}{k}}{1 - p - \frac{\delta' x}{k}}\right) - (kp + \delta x)\log\left(\frac{(p + \frac{\delta'x}{k})(1 - p - \frac{\delta x}{k})}{(p + \frac{\delta x}{k})(1 - p - \frac{\delta' x}{k})}\right)\\
=& k\left(p+\frac{\delta x}{k}\right) \log\left(\frac{p+\frac{\delta x}{k}}{p+\frac{\delta' x}{k}}\right) + k\left(1 - p-\frac{\delta x}{k}\right)\log\left(\frac{1 - p-\frac{\delta x}{k}}{1 - p-\frac{\delta' x}{k}}\right)\\
=& k\cdot h_{p + \frac{\delta'x}{k}}\left(p + \frac{\delta x}{k}\right).
\end{split}
\end{equation*}
Similarly,
\begin{equation*}
\begin{split}
\log(\phi_X(t)\exp(-t\delta x)) - \log(\phi_X(t - t')\exp(-(t\delta - t')x)) = kh_{p + \frac{\delta'x}{k}}\left(p + \frac{x}{k}\right).
\end{split}
\end{equation*}
By Lemma \ref{th:reverse-chernoff},
\begin{equation}\label{eq8}
\begin{split}
& \bbP(X \geq x) \\
\geq & \phi_X(t)\exp(-t\delta x)\cdot \left[1 - \exp\left\{-kh_{p + \frac{\delta'x}{k}}\left(p + \frac{\delta x}{k}\right)\right\} - \exp\left\{-kh_{p + \frac{\delta'x}{k}}\left(p + \frac{x}{k}\right)\right\}\right]\\
\geq & \phi_X(t\theta)\exp(-t\theta\delta x)\cdot \left[1 - \exp\left\{-kh_{p + \frac{\delta'x}{k}}\left(p + \frac{\delta x}{k}\right)\right\} - \exp\left\{-kh_{p + \frac{\delta'x}{k}}\left(p + \frac{x}{k}\right)\right\}\right]\\
= & \exp\left(-k\cdot h_{p}\left(p + \frac{\delta x}{k}\right)\right)\\
& \cdot \left[1 - \exp\left\{-kh_{p + \frac{\delta'x}{k}}\left(p + \frac{\delta x}{k}\right)\right\} - \exp\left\{-kh_{p + \frac{\delta'x}{k}}\left(p + \frac{x}{k}\right)\right\}\right].
\end{split}
\end{equation}
\eqref{eq8} will be used as a central technical tool for the rest of the analysis. Next, we consider the case that $k$ is larger than a constant $C_0(\beta)$ that only depends on $\beta$ and specifically discuss in four scenarios. We will study later the scenario that $k$ is no more than $C_0(\beta)$ based on discussions in two scenarios.

We first assume $k \geq C_0(\beta)$, where $C_0(\beta)$ is to be determined later in the analysis.
\begin{enumerate}[leftmargin=*]
	\item \bm{$p \leq \frac{1}{2k}$}. If $x \leq \frac{1}{2}$,
	\begin{equation*}
	\begin{split}
	\bbP\left(X \geq x\right) = \bbP\left(Z_1 + Z_2 + \dots + Z_k \geq 1\right) = 1 - \bbP\left(Z_1 + Z_2 + \dots + Z_k = 0\right) = 1 - (1 - p)^k.
	\end{split}
	\end{equation*}
	Here, $Z_1, \dots, Z_d$ are i.i.d. Bernoulli random variables with parameter $p$, i.e., $\bbP\left(Z_i = 1\right) = p$ and $\bbP\left(Z_i = 0\right) = 1 - p$.
	
	Next we assume $x\geq 1/2$. Let $\tau, \mu>0$. Note that
	\begin{equation}\label{eq4}
	\begin{split}
	\frac{\partial h_{p + \frac{\mu x}{k}}\left(p + \frac{\tau x}{k}\right)}{\partial \tau} = \frac{x}{k}\log\left(\frac{p + \frac{\tau x}{k}}{p + \frac{\mu x}{k}}\cdot \frac{1 - p - \frac{\mu x}{k}}{1 - p - \frac{\tau x}{k}}\right)
	\end{split}
	\end{equation}
	is an increasing function of $\tau$, immediately we know that $h_{p + \frac{\mu x}{k}}\left(p + \frac{\tau x}{k}\right)$ is a convex function of $\tau$, if $\frac{1}{2} \leq x < \frac{1 - p}{\tau \vee \mu}$.
	Since $p\leq 1/(2k)$, we have $pk/x \leq 1$ for all $x\geq 1/2$. Then for any $\tau > \mu$,
	\begin{equation*}
	\begin{split}
	\frac{p + \frac{\tau x}{k}}{p + \frac{\mu x}{k}} = \frac{\tau+\frac{pk}{x}}{\mu+\frac{pk}{x}} \geq \frac{\tau + 1}{\mu + 1}
	\end{split}
	\end{equation*}
	\begin{equation}\label{eq5}
	\begin{split}
	\text{and consequently} \quad \frac{\partial h_{p + \frac{\mu x}{k}}\left(p + \frac{\tau x}{k}\right)}{\partial \tau} \geq \frac{x}{k}\log\left(\frac{\tau + 1}{\mu + 1}\right),
	\end{split}
	\end{equation}
	for all $\frac{1}{2} \leq x < \frac{k(1 - p)}{\tau \vee \mu}$.
	Now we fix $\delta' = \frac{1 + \delta}{2}$. By convexity,
	\begin{equation*}
	\begin{split}
	h_{p + \frac{\delta'x}{k}}\left(p + \frac{\delta x}{k}\right) \geq & h_{p + \frac{\delta'x}{k}}\left(p + \frac{(3\delta +1) x}{4k}\right) + \left(\delta - \frac{(3\delta+1)}{4}\right) \cdot \frac{\partial h_{p + \frac{\delta' x}{k}}\left(p + \frac{\tau x}{k}\right)}{\partial \tau}\Big|_{\tau = \frac{3\delta+1}{4}}\\
	\geq & h_{p + \frac{\delta'x}{k}}\left(p + \frac{(3\delta + 1)x}{4k}\right) + \frac{\delta - 1}{4}\cdot\log\left(\frac{3\delta + 5}{2\delta + 6}\right)\cdot\frac{x}{k} \\
	\geq & \frac{\delta - 1}{4}\cdot\log\left(\frac{3\delta + 5}{2\delta + 6}\right)\cdot\frac{x}{k}.
	\end{split}
	\end{equation*}
	Therefore, there exists $C_1(\delta)$ that only depends on $\delta$, such that for all $C_1(\delta) \leq x  < \frac{k(1 - p)}{\delta}$,
	\begin{equation}\label{ineq:1}
	\begin{split}
	h_{p + \frac{\delta'x}{k}}\left(p + \frac{\delta x}{k}\right) \geq \frac{1}{k}.
	\end{split}
	\end{equation}
	Similarly to \eqref{eq5}, for any $\tau < \delta'$ and $\frac{1}{2} \leq x < \frac{k(1 - p)}{\delta'}$, we have
	\begin{equation*}
	\begin{split}
	\frac{\partial h_{p + \frac{\delta'x}{k}}\left(p + \frac{\tau x}{k}\right)}{\partial \tau} \leq \frac{x}{k}\log\left(\frac{\tau + 1}{\delta' + 1}\right).    	
	\end{split}
	\end{equation*}
	Thus,
	\begin{equation*}
	\begin{split}
	h_{p + \frac{\delta'x}{k}}\left(p + \frac{x}{k}\right) \geq & h_{p + \frac{\delta'x}{k}}\left(p + \frac{(\delta + 3)x}{4k}\right) - \left(\frac{(\delta+3)}{4} - 1\right) \cdot \frac{\partial h_{p + \frac{\delta' x}{k}}\left(p + \frac{\tau x}{k}\right)}{\partial \tau}\Big|_{\tau = \frac{\delta+3}{4}}\\
	\geq & h_{p + \frac{\delta'x}{k}}\left(p + \frac{(\delta + 3)x}{4k}\right) - \frac{(\delta - 1)x}{4k}\log\left(\frac{\frac{\delta + 3}{4} + 1}{\frac{\delta + 1}{2} + 1}\right) \\
	\geq & \frac{\delta - 1}{4}\cdot \log\left(\frac{2\delta + 6}{\delta + 7}\right)\cdot\frac{x}{k}.
	\end{split}
	\end{equation*}
	Therefore, there exists another constant that only relies on $\delta$: $C_2(\delta) > C_1(\delta)$, such that for all $C_2(\delta) \leq x < \frac{k(1 - p)}{\delta}$,
	\begin{equation}\label{ineq:2}
	\begin{split}
	h_{p + \frac{\delta'x}{k}}\left(p + \frac{x}{k}\right) \geq \frac{1}{k}.	
	\end{split}
	\end{equation}
	Next, assume $0<\eta<1$. By \eqref{eq4} and Taylor's Theorem, for any $b > a, 0 < x \leq \frac{(1 - \eta)(1 - p)k}{b}$, we have
	\begin{equation}\label{eq6}
	\begin{split}
	& \frac{\partial h_{p + \mu\frac{x}{k}}\left(p + \frac{\tau x}{k}\right)}{\partial \tau} \bigg|^{\tau = b}_{\tau = a} = \frac{\partial h_{p + \mu\frac{x}{k}}\left(p + \frac{b x}{k}\right)}{\partial b} - \frac{\partial h_{p + \mu\frac{x}{k}}\left(p + \frac{a x}{k}\right)}{\partial a} \\
	= & \frac{x}{k}\log\left(\frac{p + \frac{bx}{k}}{p + \frac{a x}{k}}\frac{1 - p - \frac{a x}{k}}{1 - p - \frac{b x}{k}}\right) \leq \frac{x}{k}\left(\frac{b}{a} \cdot \frac{1-p}{1-p-\frac{b}{k}\cdot \frac{(1-\eta)(1-p)k}{b}}\right) \\
	= & \frac{x}{k}\log\left(\frac{b}{a\cdot \eta}\right).
	\end{split}
	\end{equation}
	By convexity, for $0 < x \leq \frac{k(1 - p)(1 - \eta)}{\delta}$,
	\begin{equation}\label{eq13}
	\begin{split}
	h_p\left(p + \frac{\delta x}{k}\right) \leq h_p\left(p\right) + \delta\frac{\partial h_{p}\left(p + \frac{\tau x}{k}\right)}{\partial \tau}\bigg|_{\tau = \delta} = \delta\frac{\partial h_{p}\left(p + \frac{\tau x}{k}\right)}{\partial \tau}\bigg|_{\tau = \delta}.
	\end{split}
	\end{equation}
	\begin{equation}\label{eq14}
	\begin{split}
	h_p\left(p + \frac{x}{k}\right) \geq h_p\left(p + \frac{x}{2k}\right) + \frac{1}{2}\frac{\partial h_{p}\left(p + \frac{\tau x}{k}\right)}{\partial \tau}\bigg|_{\tau = \frac{1}{2}} \geq \frac{1}{2}\frac{\partial h_{p}\left(p + \frac{\tau x}{k}\right)}{\partial \tau}\bigg|_{\tau = \frac{1}{2}}.
	\end{split}
	\end{equation}
	By (\ref{eq5}) and (\ref{eq6}),
	\begin{equation*}
	\begin{split}
	\frac{\frac{\partial h_{p}\left(p + \frac{\tau x}{k}\right)}{\partial \tau}\bigg|_{\tau = \delta}}{\frac{\partial h_{p}\left(p + \frac{\tau x}{k}\right)}{\partial \tau}\bigg|_{\tau = \frac{1}{2}}} = 1 + \frac{\frac{\partial h_{p}\left(p + \frac{\tau x}{k}\right)}{\partial \tau}\bigg|_{\tau = \frac{1}{2}}^{\tau = \delta}}{\frac{\partial h_{p}\left(p + \frac{\tau x}{k}\right)}{\partial \tau}\bigg|_{\tau = \frac{1}{2}}} \leq 1 + \frac{\frac{x}{k}\cdot\log \left(\frac{2\delta}{\eta}\right)}{\frac{x}{k}\cdot\log\left(\frac{3}{2}\right)} \leq 1 +  \frac{\log \left(\frac{2\delta}{\eta}\right)}{\log\left(\frac{3}{2}\right)}.
	\end{split}
	\end{equation*}
	By \eqref{eq13}, \eqref{eq14}, and the previous inequality, there exists another constant $C'(\eta, \delta)$ that only depends on $\eta$ and $\delta$, such that for $0 \leq x \leq \frac{k(1 - p)(1 - \eta)}{\delta}$,
	\begin{equation}\label{eq7}
	\begin{split}
	h_p\left(p + \frac{\delta x}{k}\right) \leq C'(\eta, \delta)\cdot h_p\left(p + \frac{x}{k}\right).
	\end{split}
	\end{equation}
	
	Now we set $\eta = 1 - \sqrt{\frac{1}{\beta}}, \delta = \sqrt{\beta}$. Then, 
	Inequalities \eqref{eq8}, \eqref{ineq:1}, and \eqref{ineq:2} together imply
	\begin{equation*}
	\begin{split}
	\bbP\left(X \geq x\right) \geq \exp\left(-k\cdot h_{p}\left(p + \frac{x}{k}\right)\right)\cdot(1 - 2e^{-1}),
	\end{split}
	\end{equation*}
	for any $x$ satisfying $\bar C_2(\beta) = C_2(\delta) \leq x \leq \frac{k(1 - p)}{\beta}$.
	
	When $k \geq \frac{\beta\cdot \bar C_2(\beta)}{(1 - p)}$,
	set $\mu = 0, x = \frac{1}{2}$ in (\ref{eq5}), for all $0 < \tau \leq 1$,
	\begin{equation*}
	\begin{split}
	\frac{\partial h_p\left(p + \frac{\tau}{2k}\right)}{\partial \tau} \geq \frac{1}{2k}\log\left(\tau + 1\right).
	\end{split}
	\end{equation*}
	Therefore, by the property of convex function, we have
	\begin{equation*}
	\begin{split}
	h_p\left(p + \frac{1}{2k}\right) \geq&  h_p\left(p + \frac{1}{4k}\right) + \left(1 - \frac{1}{2}\right)\frac{\partial h_p\left(p + \frac{\tau}{2k}\right)}{\partial \tau}\bigg|_{\tau = \frac{1}{2}}\\ \geq& h_p\left(p + \frac{1}{4k}\right) + \frac{1}{4k}\log\left(\frac{1}{2} + 1\right)\\
	\geq& \frac{\log(\frac{3}{2})}{4k}.
	\end{split}
	\end{equation*}
	\begin{equation*}
	\begin{split}
	h_p\left(p + \frac{\bar{C}_2\left(\beta\right)}{k}\right) \leq h_p(p) + 2\bar{C}_2(\beta)\cdot\frac{\partial h_{p}\left(p + \frac{\tau}{2k}\right)}{\partial \tau}\bigg|_{\tau = 2\bar{C}_2(\beta)} = 2\bar{C}_2(\beta)\cdot\frac{\partial h_{p}\left(p + \frac{\tau}{2k}\right)}{\partial \tau}\bigg|_{\tau = 2\bar{C}_2(\beta)}.
	\end{split}
	\end{equation*}
	Also, set $x = \frac{1}{2}, \mu = 0, b = 2\bar{C}_2(\beta), a = \frac{1}{2}, \eta = 1 - \frac{1}{\beta} \leq 1 - \frac{bx}{(1 - p)k}$ in \eqref{eq6},
	\begin{equation*}
	\begin{split}
	\frac{\partial h_{p}\left(p + \frac{\tau}{2k}\right)}{\partial \tau}\bigg|_{\tau = \frac{1}{2}}^{\tau = 2\bar{C}_2(\beta)} \leq \frac{1}{2k}\log\left(\frac{4\beta\cdot\bar{C}_2(\beta)}{\beta - 1}\right).
	\end{split}
	\end{equation*}
	Similarly to \eqref{eq7}, there exists a constant $C_3(\beta) > 1$, such that for all $\frac{1}{2} \leq x \leq \bar C_2(\beta)$, we have
	\begin{equation*}
	\begin{split}
	h_p\left(p + \frac{\bar C_2(\beta)}{k}\right) \leq C_3(\beta)\cdot h_p\left(p + \frac{1}{2k}\right) \leq C_3(\beta)\cdot h_p\left(p + \frac{x}{k}\right).
	\end{split}
	\end{equation*}
	Thus for all $\frac{1}{2} \leq x \leq \bar C_2(\beta)$ and $C_{\beta} = C'(\beta)\cdot C_3(\beta) > C'(\beta)$,
	\begin{equation*}
	\begin{split}
	\bbP\left(X \geq x\right) \geq \bbP\left(X \geq \bar C_2(\beta)\right) \geq \exp\left(-C_{\beta}k\cdot h_{p}\left(p + \frac{x}{k}\right)\right)\cdot(1 - 2e^{-1}).
	\end{split}
	\end{equation*}
	
	In summary, for all $k \geq C_0(\beta):= 2\beta\cdot \bar{C}_2(\beta) \geq \frac{\beta\cdot \bar{C}_2(\beta)}{1-p}$, (which only relies on $\beta$) and $\frac{1}{2} \leq x \leq \frac{k(1 - p)}{\beta}$, there exists a constant $C_\beta>0$ that only depends on $\beta$, such that
	\begin{equation*}
	\begin{split}
	\bbP\left(X \geq x\right) \geq \exp\left(-C_{\beta}k\cdot h_{p}\left(p + \frac{x}{k}\right)\right)\cdot(1 - 2e^{-1}).
	\end{split}
	\end{equation*}
	\item
	\bm{$\frac{1}{2k} \leq p \leq \frac{1}{2}$}. By (\ref{eq4}), for any $\tau > \mu \geq 1$, we have
	\begin{equation*}
	\begin{split}
	\frac{\partial h_{p + \frac{\mu x}{k}}\left(p + \frac{\tau x}{k}\right)}{\partial \tau} = \frac{x}{k}\log\left(1 + \frac{ \frac{(\tau - \mu) x}{k}}{(p + \frac{\mu x}{k})(1 - p - \frac{\tau x}{k})}\right) \geq \frac{x}{k}\log\left(1 + \frac{(\tau - \mu)x}{k(p + \frac{\mu x}{k})}\right).
	\end{split}
	\end{equation*}
	Since $(\log(1 + t))' = \frac{1}{1+t}$ is a decreasing function of $t \geq 0$, by Taylor's Theorem,
	\begin{equation}\label{ineq:log(1+t)>}
	\begin{split}
	\log(1+t) \geq \log(1) + t\cdot (\log(1+t))' = \frac{t}{1+t}, \quad \forall t\geq 0.
	\end{split}
	\end{equation}
	Then,
	\begin{equation*}
	\begin{split}
	& \log\left(1 + \frac{\tau - \mu}{\mu}\right) \geq \frac{\frac{\tau-\mu}{\mu}}{\frac{\tau-\mu}{\mu} + 1} = \frac{\tau - \mu}{\mu}\cdot \frac{\mu}{\tau}\\
	\text{and}\quad & \log\left(1 + z\right) \geq z\cdot\frac{\mu}{\tau} \quad \forall 0\leq z \leq \frac{\tau-\mu}{\mu}.
	\end{split}
	\end{equation*}
	Since $\frac{(\tau-\mu)x}{k(p+\frac{\mu x}{k})} \leq \frac{(\tau-\mu)x}{\mu x}\leq \frac{\tau-\mu}{\mu}$, we have for all $0< x < \frac{k(1 - p)}{\tau}$,
	\begin{equation*}
	\begin{split}
	\frac{\partial h_{p + \frac{\mu x}{k}}\left(p + \frac{\tau x}{k}\right)}{\partial \tau} \geq \frac{x}{k}\frac{(\tau - \mu)x}{k(p + \frac{\mu x}{k})}\cdot\frac{\mu}{\tau} \geq \frac{(\tau - \mu)}{2\tau k}\min\left\{\mu\frac{x^2}{pk}, x\right\}.
	\end{split}
	\end{equation*}
	Since $pk \geq \frac{1}{2}$ and $\mu \geq 1$, for all $C > 1$ and $C\sqrt{pk} \leq x < \frac{k(1 - p)}{\tau}$, we have
	\begin{equation*}
	\begin{split}
	\frac{\partial h_{p + \frac{\mu x}{k}}\left(p + \frac{\tau x}{k}\right)}{\partial \tau}  \geq \frac{(\tau - \mu)C}{2\sqrt{2}\tau k}.
	\end{split}
	\end{equation*}
	We still fix $\delta' = \frac{1 + \delta}{2}$, by convexity, for all $C > 1$ and $C\sqrt{pk} \leq x < \frac{k(1 - p)}{\delta}$,
	\begin{equation}
	\begin{split}
	h_{p + \frac{\delta'x}{k}}\left(p + \frac{\delta x}{k}\right) \geq& h_{p + \frac{\delta'x}{k}}\left(p + \frac{(3\delta + 1)x}{4k}\right) + \left(\delta - \frac{3\delta + 1}{4}\right)\frac{\partial h_{p + \frac{\delta' x}{k}}\left(p + \frac{\tau x}{k}\right)}{\partial \tau}\bigg|_{\tau = \frac{3\delta + 1}{4}}\\ \geq& h_{p + \frac{\delta'x}{k}}\left(p + \frac{(3\delta + 1)x}{4k}\right) + \frac{\delta - 1}{4}\cdot\frac{\frac{\delta - 1}{4}C}{2\sqrt{2}\frac{3\delta + 1}{4}k}\\ \geq& \frac{\delta - 1}{4}\cdot\frac{\frac{\delta - 1}{4}C}{2\sqrt{2}\frac{3\delta + 1}{4}k}.
	\end{split}
	\end{equation}
	Thus, there exists a constant $C_1(\delta) > 1$, such that for all $C_1(\delta)\sqrt{pk} \leq x < \frac{k(1 - p)}{\delta}$,
	\begin{equation}\label{eq18}
	\begin{split}
	h_{p + \frac{\delta'x}{k}}\left(p + \frac{\delta x}{k}\right) \geq \frac{1}{k}.
	\end{split}
	\end{equation}
	Similarly to the derivation of \eqref{ineq:2} and \eqref{eq18}, there exists another constant $C_2(\delta) \geq C_1(\delta)$, such that for all $C_2(\delta)\sqrt{pk} \leq x < \frac{k(1 - p)}{\delta}$,
	\begin{equation}\label{eq19}
	\begin{split}
	h_{p + \frac{\delta'x}{k}}\left(p + \frac{x}{k}\right) \geq \frac{1}{k}.
	\end{split}
	\end{equation}
	For any $b>a>0$ and $0 \leq x \leq \frac{k(1 - p)(1 - \eta)}{b}$, since $p + \frac{ax}{k} \geq p$ and $1 - p - \frac{bx}{k} \geq (1 - p)\eta \geq \frac{\eta}{2}$, we have
	\begin{equation*}
	\begin{split}
		\frac{\partial h_{p + \mu\frac{x}{k}}\left(p + \frac{\tau x}{k}\right)}{\partial \tau} \bigg|^{\tau = b}_{\tau = a} \overset{\eqref{eq4}}{=}& \frac{x}{k}\log\left(1 + \frac{(b - a)\frac{x}{k}}{\left(p + \frac{a x}{k}\right)\left(1 - p - \frac{b x}{k}\right)}\right)\\ \leq& \frac{x}{k}\log\left(1 + \frac{(b - a)\frac{x}{k}}{p(1 - p)\eta}\right).
	\end{split}
	\end{equation*}
	\begin{equation*}
	\begin{split}
	\frac{\partial h_{p}\left(p + \frac{\tau x}{k}\right)}{\partial \tau} = \frac{x}{k}\log\left(1 + \frac{\tau\frac{x}{k}}{p(1 - p - \frac{\tau x}{k})}\right) \geq \frac{x}{k}\log\left(1 + \frac{\tau\frac{x}{k}}{p(1 - p)}\right).
	\end{split}
	\end{equation*}
	Note that $\log(1 + ax) \leq a\log(1 + x)$ holds for all $x \geq 0$ and $a \geq 1$,
	\begin{equation}\label{eq15}
	\begin{split}
	\frac{\frac{\partial h_{p}\left(p + \frac{\tau x}{k}\right)}{\partial \tau}\bigg|_{\tau = \delta}}{\frac{\partial h_{p}\left(p + \frac{\tau x}{k}\right)}{\partial \tau}\bigg|_{\tau = \frac{1}{2}}} = 1 + \frac{\frac{\partial h_{p}\left(p + \frac{\tau x}{k}\right)}{\partial \tau}\bigg|_{\tau = \frac{1}{2}}^{\tau = \delta}}{\frac{\partial h_{p}\left(p + \frac{\tau x}{k}\right)}{\partial \tau}\bigg|_{\tau = \frac{1}{2}}} \leq 1 + \frac{\frac{x}{k}\cdot\log\left(1 + \frac{(\delta - \frac{1}{2})\frac{x}{k}}{p(1 - p)\eta}\right)}{\frac{x}{k}\cdot\log\left(1 + \frac{\frac{x}{k}}{2p(1 - p)}\right)} \leq 1 +  \frac{2\delta - 1}{\eta}.
	\end{split}
	\end{equation}
	By (\ref{eq13}), (\ref{eq14}) and (\ref{eq15}), there exists $C(\eta, \delta)$ that only depends on $\eta$ and $\delta$, such that for $0 \leq x \leq \frac{k(1 - p)(1 - \eta)}{\delta}$,
	\begin{equation}\label{eq17}
	\begin{split}
	h_p\left(p + \frac{\delta x}{k}\right) \leq C(\eta, \delta)\cdot h_p\left(p + \frac{x}{k}\right).
	\end{split}
	\end{equation}
	Now, we set $\delta = \sqrt{\beta}, \eta = 1 - \sqrt{\frac{1}{\beta}}$.
	\begin{itemize}[leftmargin=*]
		\item The previous discussions on \eqref{eq8}, \eqref{eq18}, \eqref{eq19}, and \eqref{eq17} imply that there exists $\bar{C}_2(\beta) = C_2(\sqrt{\beta})$ and $C(\beta)$ that only relies on $\beta$, such that for all 
		$\bar{C}_2(\beta)\sqrt{pk} \leq x \leq \frac{k(1 - p)}{\beta}$,
		\begin{equation*}
		\begin{split}
		\bbP\left(X \geq x\right) \geq \exp\left(-C(\beta)k\cdot h_{p}\left(p + \frac{x}{k}\right)\right)\cdot(1 - 2e^{-1}).
		\end{split}
		\end{equation*}
		\item Noting that $\log(1 + x) \leq x$ for all $x > -1$, for $x \leq \bar C_2(\beta)\sqrt{pk} \leq \frac{k(1 - p)}{\beta}, 0 \leq \tau \leq 1$, we have
		\begin{equation*}
		\begin{split}
		\frac{\partial h_p(p + \frac{\tau x}{k})}{\partial \tau} = \frac{x}{k}\log\left(1 + \frac{\frac{\tau x}{k}}{p(1 - p - \frac{\tau x}{k})}\right) \leq \frac{\tau x^2}{pk^2(1 - p - \frac{\tau x}{k})} \leq \frac{\tau x^2}{pk^2\cdot (1 - p)(1 - \frac{1}{\beta})}.
		\end{split}
		\end{equation*}
		Since $\frac{\partial h_p(p + \frac{\tau x}{k})}{\partial \tau}$ is an increasing function for $0\leq \tau \leq 1$, we have
		\begin{equation}\label{ineq:h_p(p+x/k)<=}
		\begin{split}
		h_p\left(p + \frac{x}{k}\right) \leq h_p(p) + \frac{\partial h_p(p + \frac{\tau x}{k})}{\partial \tau}\bigg|_{\tau = 1} \leq \frac{2x^2}{pk^2\cdot (1 - \frac{1}{\beta})},
		\end{split}
		\end{equation}
		which means
		\begin{equation*}
		\begin{split}
		h_p\left(p + \frac{\bar C_2(\beta)\sqrt{pk}}{k}\right) \leq \frac{2\bar C_2^2(\beta)}{k\left(1 - \frac{1}{\beta}\right)}.
		\end{split}
		\end{equation*}
		Consequently, for any $x\leq \bar{C}_2(\beta)\sqrt{pk} \leq \frac{k(1-\beta)}{\beta}$,
		\begin{equation*}
		\begin{split}
		\bbP\left(X \geq x\right) \geq & \bbP\left(X \geq \bar C_2(\beta)\sqrt{pk}\right) \geq \exp\left(-C(\beta)\frac{2 \bar C_2^2(\beta)}{1 - \frac{1}{\beta}}\right)\cdot (1 - 2e^{-1})\\
		\geq & \left[\exp\left(-C(\beta)\frac{2 \bar C_2^2(\beta)}{1 - \frac{1}{\beta}}\right)(1-2e^{-1})\right]\exp\left(-C(\beta) k\cdot h_p\left(p+\frac{x}{k}\right)\right).
		\end{split}
		\end{equation*}
	\end{itemize}
	Denote $C_0(\beta) = 2\bar{C}_2^2(\beta)\cdot\beta^2$. In summary of the previous two bullet points, there exist $c_{\beta}, C_{\beta} > 0$ that only relies on $\beta$, such that for all $k \geq C_0(\beta) \geq \frac{\bar{C}_2^2(\beta)\cdot\beta^2 p}{(1 - p)^2}$ and $0 \leq x \leq \frac{k(1 - p)}{\beta}$, we have
	\begin{equation*}
	\begin{split}
	\bbP\left(X \geq x\right) \geq c_{\beta}\exp\left(-C_{\beta}h_p\left(p + \frac{x}{k}\right)\right).
	\end{split}
	\end{equation*}
	\item \bm{$\frac{1}{2} \leq p \leq 1 - \frac{\widetilde{C}_{\beta}}{k}$}. Here, $\tilde{C}_\beta$ is a to-be-specified constant that only relies on $\beta$.
	
	By \eqref{eq4}, for any $0 < \eta < 1$, $\tau > \mu \geq 0$, and $x \leq \frac{(1 - \eta)(1 - p)k}{\tau}$, we have
	\begin{equation*}
	\begin{split}
	\frac{\partial h_{p + \frac{\mu x}{k}}\left(p + \frac{\tau x}{k}\right)}{\partial \tau} = \frac{x}{k}\log\left(1 + \frac{ \frac{(\tau - \mu) x}{k}}{(p + \frac{\mu x}{k})(1 - p - \frac{\tau x}{k})}\right) \geq \frac{x}{k}\log\left(1 + \frac{(\tau - \mu)x}{k(1 - p)}\right).
	\end{split}
	\end{equation*}
	Since $\frac{(\tau - \mu)x}{k(1 - p)} \leq \frac{\tau - \mu}{\tau}$ and $\frac{\log(1 + x)}{x}$ is a decreasing function for $x > 0$,
	\begin{equation*}
	\begin{split}
	\log\left(1 + \frac{(\tau - \mu)x}{k(1 - p)}\right) \geq \frac{(\tau - \mu)x}{k(1 - p)}\cdot \frac{\log\left(1 + \frac{\tau - \mu}{\tau}\right)}{\frac{\tau - \mu}{\tau}} \overset{\eqref{ineq:log(1+t)>}}{\geq} \frac{(\tau - \mu)x}{k(1 - p)}\cdot \frac{1}{\frac{\tau - \mu}{\tau} + 1}.
	\end{split}
	\end{equation*}
	Fix $\delta' = \frac{1 + \delta}{2}$, similar to the discussions in the case that $\frac{1}{2k} \leq p \leq \frac{1}{2}$, we can find a constant $C_1(\delta) > 0$ that only depends on $\delta$, such that for all $k \geq \frac{\left(\delta C_1(\delta)\right)^2}{(1 - \eta)^2(1 - p)}$ and $C_1(\delta)\sqrt{(1 - p)k} \leq x \leq \frac{(1 - \eta)(1 - p)k}{\delta}$, we have
	\begin{equation}\label{eq16}
	\begin{split}
	h_{p + \frac{\delta'x}{k}}\left(p + \frac{\delta x}{k}\right) \geq \frac{1}{k} \quad \text{and} \quad h_{p + \frac{\delta'x}{k}}\left(p + \frac{x}{k}\right) \geq \frac{1}{k}.
	\end{split}
	\end{equation}
	Let $0 < \eta < 1$ be a fixed and to-be-specified value. For any $x > 0, b > a$, and $x \leq \frac{(1 - \eta)(1 - p)k}{b}$, we have $p + \frac{ax}{k} \geq p \geq \frac{1}{2}$ and $1 - p \geq 1 - p - \frac{bx}{k} \geq (1 - p)\eta$. By \eqref{eq4} and Taylor's Theorem, we have
	\begin{equation}\label{eq10}
    \begin{split}
    	\frac{\partial h_{p + \mu\frac{x}{k}}\left(p + \frac{\tau x}{k}\right)}{\partial \tau} \bigg|^{\tau = b}_{\tau = a} \overset{\eqref{eq4}}{=}& \frac{x}{k}\log\left(1 + \frac{(b - a)\frac{x}{k}}{\left(p + \frac{a x}{k}\right)\left(1 - p - \frac{b x}{k}\right)}\right)\\ \leq& \frac{x}{k}\log\left(1 + \frac{2(b - a)\frac{x}{k}}{\eta(1 - p)}\right),
    \end{split}
	\end{equation}
	\begin{equation}\label{eq11}
	\begin{split}
	\frac{\partial h_{p}\left(p + \frac{\tau x}{k}\right)}{\partial \tau} \bigg|_{\tau = b} \overset{\eqref{eq4}}{=} \frac{x}{k}\log\left(1 +\frac{\frac{bx}{k}}{p\left(1 - p - \frac{b x}{k}\right)}\right) \geq \frac{x}{k}\log\left(1 +\frac{\frac{bx}{k}}{1 - p}\right).
	\end{split}
	\end{equation}
	Noticing that $\frac{\log(1 + u_1)}{u_1} \leq \frac{\log(1 + u_2)}{u_2}$ for all $u_1 \geq u_2\geq 0$, and $\frac{2(\delta-1/2)\frac{x}{k}}{\eta(1 - p)} \geq \frac{\frac{x}{2k}}{1 - p}$, we have
	\begin{equation*}
	\begin{split}
	\frac{\frac{\partial h_{p}\left(p + \frac{\tau x}{k}\right)}{\partial \tau}\bigg|_{\tau = \delta}}{\frac{\partial h_{p}\left(p + \frac{\tau x}{k}\right)}{\partial \tau}\bigg|_{\tau = \frac{1}{2}}} = 1 + \frac{\frac{\partial h_{p}\left(p + \frac{\tau x}{k}\right)}{\partial \tau}\bigg|_{\tau = \frac{1}{2}}^{\tau = \delta}}{\frac{\partial h_{p}\left(p + \frac{\tau x}{k}\right)}{\partial \tau}\bigg|_{\tau = \frac{1}{2}}} \overset{\eqref{eq10}\eqref{eq11}}{\leq} 1+\frac{\log\left(1 + \frac{2(\delta-1/2)\frac{x}{k}}{\eta(1 - p)}\right)}{\log\left(1 +\frac{\frac{x}{2k}}{1 - p}\right)} \leq 1 + \frac{4\delta - 2}{\eta}.
	\end{split}
	\end{equation*}
	Using (\ref{eq13}) and (\ref{eq14}) again, we can find a constant $C'(\alpha, \delta)$ that only depends on $\eta$ and $\delta$, such that for $x \leq \frac{k(1 - p)(1 - \eta)}{\delta}$,
	\begin{equation}\label{eq12}
	\begin{split}
	h_p\left(p + \frac{\delta x}{k}\right) \leq C'(\eta, \delta)\cdot h_p\left(p + \frac{x}{k}\right).
	\end{split}
	\end{equation}
	Now we specify $\delta = \sqrt{\beta}$, $\eta = 1 - \frac{1}{\sqrt{\beta}}$, $\bar C_1(\beta) = C_1(\delta)$, and $\widetilde{C}_{\beta} = \beta\cdot\bar C_1(\beta)$, which all only depend on $\beta$. For any $k \geq 2\widetilde C_\beta$, we have $1 - \frac{\widetilde{C}_{\beta}}{k} \geq \frac{1}{2}$ and $\frac{(1 - p)k}{\beta} \geq \sqrt{(1 - p)k}\cdot\frac{\widetilde C_\beta}{\beta} = \bar{C}_1(\beta)\sqrt{(1 - p)k}$.
	\begin{itemize}[leftmargin=*]
		\item For all $\bar{C}_1(\beta)\sqrt{(1 - p)k} \leq x \leq \frac{(1 - p)k}{\beta}$, \eqref{eq16}, \eqref{eq12} and \eqref{eq8} lead us to
		\begin{equation*}
		\begin{split}
		\bbP\left(X \geq x\right) \geq \exp\left(-C_{\beta} k \cdot h_p\left(p + \frac{x}{k}\right)\right)\cdot (1 - 2e^{-1}).
		\end{split}
		\end{equation*}
		\item For $x \leq \bar{C}_1(\beta)\sqrt{(1 - p)k} \leq \frac{k(1 - p)}{\beta}$ and $0 \leq \tau \leq 1$,
		\begin{equation*}
		\begin{split}
		\frac{\partial h_p(p + \frac{\tau x}{k})}{\partial \tau} = \frac{x}{k}\log\left(1 + \frac{\frac{\tau x}{k}}{p(1 - p - \frac{\tau x}{k})}\right) \leq \frac{\tau\frac{x^2}{k^2}}{p(1 - p - \frac{\tau x}{k})} \leq \frac{2\tau\frac{x^2}{k^2}}{(1 - p)(1 - \frac{1}{\beta})}.
		\end{split}
		\end{equation*}
		Therefore,
		\begin{equation*}
		\begin{split}
		h_p\left(p + \frac{x}{k}\right) \overset{\eqref{ineq:h_p(p+x/k)<=}}{\leq} & h_p(p) + \frac{\partial h_p(p + \frac{\tau x}{k})}{\partial \tau}\bigg|_{\tau = 1} \leq \frac{2x^2}{k^2\cdot (1 - p)(1 - \frac{1}{\beta})}\\
		\leq & \frac{2\bar{C}_1(\beta)^2}{k(1-\frac{1}{\beta})}, \quad \forall 0\leq x \leq \bar{C}_1(\beta)\sqrt{(1-p)k}\leq \frac{k(1-p)}{\beta}.
		\end{split}
		\end{equation*}
		Consequently, for any $0\leq x \leq \bar{C}_1(\beta)\sqrt{(1-p)k}\leq \frac{k(1-p)}{\beta}$,
		\begin{equation*}
		\begin{split}
		\bbP\left(X \geq x\right) \geq & \bbP\left(X \geq \bar{C}_1(\beta)\sqrt{(1 - p)k}\right) \geq \exp\left(-C_{\beta}\cdot \frac{2}{1 - \frac{1}{\beta}}\right)\cdot (1 - 2e^{-1})\\
		\geq & \exp\left(-C_\beta k \cdot h_p(p+\frac{x}{k})\right)\cdot\left[\exp\left(-C_{\beta}\cdot \frac{2}{1 - \frac{1}{\beta}}\right)\cdot (1 - 2e^{-1})\right].
		\end{split}
		\end{equation*}
		for some uniform constant $C_\beta>0$ that only depends on $\beta$.
	\end{itemize}
	In summary from the previous two bullet points, there exist $c_{\beta}, C_{\beta} > 0$, which only rely on $\beta$, such that for all $k \geq C_0(\beta) := 2\widetilde{C}_{\beta}$, $0 \leq x \leq \frac{k(1 - p)}{\beta}$, we have
	\begin{equation*}
	\begin{split}
	\bbP\left(X\geq x\right) \geq c_{\beta}\exp\left(-C_{\beta}k\cdot h_p\left(p + \frac{x}{k}\right)\right).
	\end{split}
	\end{equation*}
	\item
	\bm{$1 - \frac{\widetilde{C}_{\beta}}{k} \leq p \leq 1$}. If $k \geq \widetilde{C}_{\beta} + 1$ and $x \leq k(1 - p)$,
	\begin{equation*}
	\begin{split}
	\bbP\left(X \geq x\right) \geq \bbP\left(Z_1 + \cdots + Z_k = k\right) = \left[\bbP\left(Z_1 = 1\right)\right]^k = p^k \geq \left(1 - \frac{\widetilde C_{\beta}}{k}\right)^k.
	\end{split}
	\end{equation*}
	By arithmetic mean-geometric mean inequality, for all integer $m \geq \widetilde{C}_{\beta} + 1$,
	\begin{equation*}
	\begin{split}
	\left(1 - \frac{\widetilde C_{\beta}}{m}\right)^m = 1 \cdot \left(1 - \frac{\widetilde C_{\beta}}{m}\right)^m \leq \left(\frac{1 + m\cdot \left(1 - \frac{\widetilde C_{\beta}}{m}\right)}{m + 1}\right)^{m + 1} \leq \left(1 - \frac{\widetilde C_{\beta}}{m + 1}\right)^{m + 1}.
	\end{split}
	\end{equation*}
	Thus, for all $k \geq \widetilde{C}_{\beta} + 1$,
	\begin{equation*}
	\begin{split}
	\bbP\left(X \geq x\right) \geq \left(1 - \frac{\widetilde C_{\beta}}{k}\right)^k \geq \left(1 - \frac{\widetilde C_{\beta}}{k-1}\right)^{k-1} \geq \cdots \geq \left(1 - \frac{\widetilde{C}_{\beta}}{\lceil \widetilde{C}_{\beta}\rceil + 1}\right)^{\left\lceil\widetilde{C}_{\beta}\right\rceil + 1}.
	\end{split}
	\end{equation*}
	Therefore, there exists $c_{\beta} = \left(1 - \frac{\widetilde{C}_{\beta}}{\lceil\widetilde{C}_{\beta}\rceil + 1}\right)^{\left\lceil\widetilde{C}_{\beta}\right\rceil + 1} > 0$ that only relies on $\beta$, such that for all $k \geq \widetilde{C}_{\beta} + 1$ and $0 \leq x \leq (1 - p)k$,
	\begin{equation*}
	\begin{split}
	\bbP\left(X \geq x\right) \geq c_{\beta} \geq c_{\beta}\exp\left(-h_p\left(p + \frac{x}{k}\right)\right).
	\end{split}
	\end{equation*}
\end{enumerate}
Finally, it remains to consider the case that $k \leq C_0(\beta)$, where $C_0(\beta)$ is any constant that only relies on $\beta$. We discuss in two scenarios.
\begin{enumerate}[leftmargin=*]
	\item For any $0 < p \leq \frac{1}{2}$ and $1 \leq pk + x \leq k$,
	\begin{equation*}
	\begin{split}
	\bbP\left(X \geq x\right) = & \sum_{i = \lceil kp + x\rceil}^{k}\binom{k}{i}p^i(1 - p)^{k - i} \geq p^{\lceil pk + x\rceil}(1 - p)^{k - \lceil pk + x\rceil}\\
	\geq & p^{\lceil pk + x\rceil}\cdot \frac{1}{2^{k - \lceil pk + x\rceil}}.
	\end{split}
	\end{equation*}
	Since $pk + x \geq 1$, $\lceil pk + x\rceil \leq 2\left(pk + x\right)$, we have
	\begin{equation}\label{ineq:k-constant-1}
	\begin{split}
	\bbP\left(X \geq x\right) \geq \frac{1}{2^{C_0(\beta)}}p^{2\left(pk + x\right)}.
	\end{split}
	\end{equation}
	Also for any $0 \vee (1 - pk) < x \leq \frac{(1 - p)k}{\beta}$,
	\begin{equation*}
	\begin{split}
		\exp\left(-2k\cdot h_p\left(p + \frac{x}{k}\right)\right) =& \left(\frac{p}{p + \frac{x}{k}}\right)^{2pk + 2x}\left(\frac{1 - p}{1 - p - \frac{x}{k}}\right)^{2k(1 - p) - 2x}\\ \leq& \left(\frac{p}{p + \frac{x}{k}}\right)^{2pk + 2x}\left(\frac{1}{1 - \frac{1}{\beta}}\right)^{2k}.
	\end{split}
	\end{equation*}
	Since $pk + x \geq 1$, we have $(p + kx)^{p + kx} \geq 1$; we also have $pk+x\leq k \leq C_0(\beta)$. Then,
	\begin{equation}\label{ineq:k-constant-2}
	\begin{split}
	& \exp\left(-2k\cdot h_p\left(p + \frac{x}{k}\right)\right) \leq \left(\frac{kp}{pk + x}\right)^{2pk + 2x} \left(\frac{1}{1 - \frac{1}{\beta}}\right)^{2k} \\
	\leq & (C_0(\beta) p)^{2pk+2x} \left(\frac{1}{1-\frac{1}{\beta}}\right)^{2C_0(\beta)} \leq p^{2pk + 2x}\left(\frac{C_0(\beta)}{1 - \frac{1}{\beta}}\right)^{2C_0(\beta)}.
	\end{split}
	\end{equation}
	Combining \eqref{ineq:k-constant-1} and \eqref{ineq:k-constant-2}, we conclude that there exist $C_{\beta} > 0$ and $c_{\beta} > 0$ that only rely on $\beta$, such that for any $k \leq C_0(\beta)$ and $0 \vee (1 - pk) < x \leq \frac{(1 - p)k}{\beta}$, we have
	\begin{equation}
	\begin{split}
	\bbP\left(X \geq x\right) \geq c_{\beta}\cdot\exp\left(-C_{\beta}k\cdot h_p\left(p + \frac{x}{k}\right)\right).
	\end{split}
	\end{equation}
	
	\item
	Similarly to the derivation of \eqref{ineq:k-constant-1} in the previous scenario, for any $\frac{1}{2} \leq p < 1, pk \leq pk + x \leq k-1$, we have
	\begin{equation*}
	\begin{split}
	\bbP\left(X \geq x\right)  \geq \frac{1}{2^{C_0(\beta)}}(1 - p)^{k(1 - p) - x} \geq \frac{1}{2^{C_0(\beta)}}(1 - p)^{k(1 - p)}.
	\end{split}
	\end{equation*}
	Noting that the function $u^u$ ($0<u\leq 1$) takes the minimum at $u=1/e$, we have
	\begin{equation*}
	\begin{split}
	\frac{1}{2^{C_0(\beta)}} (1-p)^{k(1-p)} = \frac{1}{2^{C_0(\beta)}} \left((1-p)^{1-p}\right)^{k} \geq \frac{1}{2^{C_0(\beta)}}(1/e)^{k/e} \geq  \frac{1}{2^{C_0(\beta)}}\left(\frac{1}{e}\right)^{\frac{C_0(\beta)}{e}},
	\end{split}
	\end{equation*}
	which can be lower bounded by a positive constant depending on $\beta$.
	Similarly as the proof of \eqref{ineq:k-constant-2}, for any $0 < x \leq \frac{(1 - p)k}{\beta}$, we have
	\begin{equation*}
	\begin{split}
	\exp\left(-k\cdot h_p\left(p + \frac{x}{k}\right)\right) \leq \left(\frac{1}{1 - \frac{1}{\beta}}\right)^{C_0(\beta)},
	\end{split}
	\end{equation*}
	which can be upper bounded by a constant depending on $\beta$.
	Therefore there exists $c_{\beta} > 0$ that only relies on $\beta$, such that for $\frac{1}{2} \leq p < 1$, $k \leq C_0(\beta)$ and
	$0 \leq x \leq \frac{(1 - p)k}{\beta}$, we have
	\begin{equation*}
	\begin{split}
	\bbP\left(X \geq x\right) \geq c_{\beta}\cdot\exp\left(-k\cdot h_p\left(p + \frac{x}{k}\right)\right).
	\end{split}
	\end{equation*}
\end{enumerate}

In summary, we have finished the proof of this theorem.\quad $\square$

\subsection{Proof of Theorem \ref{cor:Poisson tail lower bound}}\label{sec:proof-Poisson}
We first consider the right tail.
We discuss in two scenarios: $\lambda + x \leq 1$ or $\lambda + x > 1$.
\begin{itemize}[leftmargin=*]
	\item
	$\bm{\lambda + x \leq 1}.$
	\begin{equation*}
	\begin{split}
	\bbP\left(X \geq x\right) = \bbP\left(Y \geq \lambda + x\right) = \bbP\left(Y \geq 1\right) = 1 - \bbP(Y = 1) = 1 - e^{-\lambda}.
	\end{split}
	\end{equation*}
	\item
	$\bm{\lambda + x > 1}.$
	Suppose $Y_n \sim \text{Bi}(n, \frac{\lambda}{n})$ and $X_n = Y_n - \lambda$. Then $Y_n \stackrel{d.}{\rightarrow} Y$, and consequently, $X_n \stackrel{d.}{\rightarrow} X$, i.e.,
	\begin{equation}\label{eq20}
	\begin{split}
	\bbP\left(X \geq x\right) = \lim_{n \to \infty}\bbP\left(X_n \geq x\right).
	\end{split}
	\end{equation}
	For all $n \geq \max\{2\lambda, 4x\}$, we have
	\begin{equation*}
	\begin{split}
	\frac{n(1 - \frac{\lambda}{n})}{2} \geq \frac{4x\cdot(1 - \frac{\lambda}{2\lambda})}{2} = x.
	\end{split}
	\end{equation*}
	Note that $n\cdot\left(\frac{\lambda}{n}\right) + x = \lambda + x > 1$, by Theorem \ref{th:Sum of i.i.d. Bernoulli tail lower bound} ($\beta = 2, p = \frac{\lambda}{n}$), there exist constants $c, C > 0$, for all $n \geq \max\{2\lambda, 4x\}$,
	\begin{equation*}
	\begin{split}
	\bbP\left(X_n \geq x\right) \geq c\cdot\exp\left(-Cn\cdot h_{\frac{\lambda}{n}}\left(\frac{\lambda}{n} + \frac{x}{n}\right)\right).
	\end{split}
	\end{equation*}
	By \eqref{eq20} and the previous inequality, there exist constants $c, C > 0$ such that
	\begin{equation}\label{eq22}
	\begin{split}
	\bbP\left(X \geq x\right) \geq \varliminf_{n \to \infty}c\cdot\exp\left(-Cn\cdot h_{\frac{\lambda}{n}}\left(\frac{x + \lambda}{n}\right)\right) = c\cdot\exp\left(-C\left(\varlimsup_{n \to \infty}n\cdot h_{\frac{\lambda}{n}}\left(\frac{x + \lambda}{n}\right)\right)\right).
	\end{split}
	\end{equation}
	For all $x \geq 0, n \geq \lambda + x$,
	\begin{equation}\label{eq131}
	\begin{split}
	n\cdot h_{\frac{\lambda}{n}}\left(\frac{x + \lambda}{n}\right) =& n\left(\frac{x + \lambda}{n}\log\left(\frac{\frac{x + \lambda}{n}}{\frac{\lambda}{n}}\right) + \left(1 - \frac{x + \lambda}{n}\right)\log\left(\frac{1 - \frac{x + \lambda}{n}}{1 - \frac{\lambda}{n}}\right)\right)\\
	=& \left(x + \lambda\right)\log\left(\frac{x + \lambda}{\lambda}\right) - x\cdot\left(\frac{1 - \frac{x + \lambda}{n}}{\frac{x}{n}}\right)\log\left(1 + \frac{\frac{x}{n}}{1 - \frac{x + \lambda}{n}}\right),
	\end{split}
	\end{equation}
	Noting that $\lim_{y \to 0}\frac{\log(1 + y)}{y} = 1$ and $\frac{\frac{x}{n}}{1 - \frac{x + \lambda}{n}} \to 0$ as $n \to \infty$, we conclude that
	\begin{equation}\label{eq21}
	\begin{split}
	\lim_{n \to \infty}n\cdot h_{\frac{\lambda}{n}}\left(\frac{x + \lambda}{n}\right) =& \left(x + \lambda\right)\log\left(\frac{x + \lambda}{\lambda}\right) - x\lim_{n \to \infty}\left(\frac{1 - \frac{x + \lambda}{n}}{\frac{x}{n}}\right)\log\left(1 + \frac{\frac{x}{n}}{1 - \frac{x + \lambda}{n}}\right)\\
	=&  \lambda\left[\left(1 + \frac{x}{\lambda}\right)\log\left(1 + \frac{x}{\lambda}\right) - \frac{x}{\lambda}\right]
	= \frac{x^2}{2\lambda}\psi_{Benn}(x/\lambda).
	\end{split}
	\end{equation}
	\eqref{eq22} and \eqref{eq21} lead us to
	\begin{equation*}
	\begin{split}
	\bbP\left(X \geq x\right) \geq c\cdot\exp\left(-C\cdot\frac{x^2}{2\lambda}\psi_{Benn}(x/\lambda)\right).
	\end{split}
	\end{equation*}
\end{itemize}
Then we consider the left tail. For all $0 \leq x \leq \frac{\lambda}{\beta}$, since $X_n \stackrel{d.}{\rightarrow} X$,
\begin{equation}\label{eq23}
\begin{split}
\bbP\left(X \leq -x\right) = \lim_{n \to \infty}\bbP\left(X_n \leq -x\right).
\end{split}
\end{equation}
For all $n \geq \lambda + 1$, we have $x + n\left(1 - \frac{\lambda}{n}\right) \geq n - \lambda = 1$, by Theorem \ref{th:Sum of i.i.d. Bernoulli tail lower bound}, there exist two constants $c_{\beta}, C_{\beta} > 0$ depending only on $\beta$ such that for all $0 \leq x \leq \frac{n\cdot\left(\frac{\lambda}{n}\right)}{\beta}$,
\begin{equation*}
\begin{split}
\bbP\left(X_n \leq -x\right) \geq c_\beta \exp\left(-C_\beta n\cdot h_{\frac{\lambda}{n}}\left(\frac{\lambda}{n} - \frac{x}{n}\right)\right).
\end{split}
\end{equation*}
\eqref{eq23} and the previous inequality lead us to
\begin{equation}
\begin{split}
\bbP\left(X \leq -x\right) \geq c_{\beta}\cdot\exp\left(-C_{\beta}\left(\varliminf_{n \to \infty}n\cdot h_{\frac{\lambda}{n}}\left(\frac{\lambda - x}{n}\right)\right)\right).
\end{split}
\end{equation}
Similarly to \eqref{eq21}, for all $x < \lambda$ and $n \geq \lambda + 1$,
\begin{equation}\label{eq24}
\begin{split}
\lim_{n \to \infty}n\cdot h_{\frac{\lambda}{n}}\left(\frac{\lambda - x}{n}\right) = \frac{x^2}{2\lambda}\psi_{Benn}(-x/\lambda).
\end{split}
\end{equation}
\eqref{eq23} and \eqref{eq24} imply that for all $0 \leq x \leq \frac{\lambda}{\beta}$,
\begin{equation*}
\begin{split}
\bbP\left(X \leq -x\right) \geq c_{\beta}\cdot\exp\left(-C_{\beta}\cdot\frac{x^2}{2\lambda}\psi_{Benn}(-x/\lambda)\right).
\end{split}
\end{equation*}

\subsection{Proof of Corollary \ref{thm:Irwin-Hall Distribution}}\label{sec:Irwin-Hall}
Suppose $U \sim U[0, 1]$, by convexity of $e^{tx}$,
\begin{equation}\label{eq50}
\begin{split}
\phi_{X}(t) = \left[\bbE e^{t(U - \frac{1}{2})}\right]^k \leq \left[\bbE\left((1 - U)e^{-\frac{t}{2}} + Ue^{\frac{t}{2}}\right)\right]^k = \left(\frac{1}{2}\left(e^{-\frac{t}{2}} + e^{\frac{t}{2}}\right)\right)^k.
\end{split}
\end{equation}
Denote $X_1 = \sum_{i = 1}^kV_i$, where $V_i$ are i.i.d. random variables with $\bbP\left(V_i = \frac{1}{2}\right) = \bbP\left(V_i = -\frac{1}{2}\right) = \frac{1}{2}$, then the RHS of \eqref{eq50} is exactly the moment generating function of $X_1$.
For all $0 \leq \tau \leq 1$ and $0 \leq x \leq \frac{k}{2}$,
\begin{equation*}
\begin{split}
\frac{\partial h_{\frac{1}{2}}\left(\frac{1}{2} + \frac{\tau x}{k}\right)}{\partial \tau} = \frac{x}{k}\log\left(1 + \frac{\frac{\tau x}{k}}{\frac{1}{2}\left(\frac{1}{2} - \frac{\tau x}{k}\right)}\right) \geq \frac{x}{k}\log\left(1 + \frac{4\tau x}{k}\right) \geq \frac{x}{k}\cdot\frac{1}{2}\frac{4\tau x}{k} = \frac{2\tau x^2}{k^2}.
\end{split}
\end{equation*}
The second inequality holds since $\inf_{0 \leq t\leq 2}\frac{\log(1 + t)}{t} = \frac{\log 3}{2} \geq \frac{1}{2}$.
\begin{equation}\label{eq128}
\begin{split}
h_{\frac{1}{2}}\left(\frac{1}{2} + \frac{x}{k}\right) = h_{\frac{1}{2}}\left(\frac{1}{2}\right) + \int_{0}^{1}\frac{\partial h_{\frac{1}{2}}\left(\frac{1}{2} + \frac{\tau x}{k}\right)}{\partial \tau}d\tau \geq 0 + \int_{0}^1\frac{2\tau x^2}{k^2}d\tau = \frac{x^2}{k^2}.
\end{split}
\end{equation}
By the Chernoff Bound, \cite{boucheron2013concentration} page 22-23, and \eqref{eq128}, we know that for any $0 \leq x \leq \frac{k}{2}$,
\begin{equation*}
\begin{split}
\bbP\left(X \geq x\right) \leq \inf_{t\geq 0}\frac{\phi_{X}(t)}{e^{tx}} \leq  \inf_{t\geq 0}\frac{\phi_{X_1}(t)}{e^{tx}} \leq \exp\left(-k\cdot h_{\frac{1}{2}}\left(\frac{1}{2} + \frac{x}{k}\right)\right) \leq \exp\left(-\frac{x^2}{k}\right).
\end{split}
\end{equation*}
Next, we consider the lower bound. The arithmetic mean-geometric mean inequality tells us
\begin{equation}\label{eq51}
\begin{split}
\phi_{X}(t) =& \left(\int_{0}^{1}e^{t(u - \frac{1}{2})}du\right)^k = \left(\int_{0}^{\frac{1}{2}}\frac{1}{2}\left(e^{tu} + e^{t(\frac{1}{2} - u)}\right)du + \int_{-\frac{1}{2}}^{0}\frac{1}{2}\left(e^{tu} + e^{t(-\frac{1}{2} - u)}\right)du\right)^k\\
\geq& \left(\frac{1}{2}\left(e^{-\frac{t}{4}} + e^{\frac{t}{4}}\right)\right)^k.
\end{split}
\end{equation}
Combine \eqref{eq50} and \eqref{eq51} together, we know that for all $t \geq 0$,
\begin{equation}
\begin{split}
\left(\frac{1}{2}\left(e^{-\frac{t}{4}} + e^{\frac{t}{4}}\right)\right)^k \leq \phi_{X}(t) \leq \left(\frac{1}{2}\left(e^{-\frac{t}{2}} + e^{\frac{t}{2}}\right)\right)^k.
\end{split}
\end{equation}
By L'Hopital's rule,
$\lim_{t \to 0}\frac{\log\left(\frac{1}{2}\left(e^{-\frac{t}{4}} + e^{\frac{t}{4}}\right)\right)}{t^2} = \frac{1}{32}.$
Moreover, for all $t > 0$, $\frac{\log\left(\frac{1}{2}\left(e^{-\frac{t}{4}} + e^{\frac{t}{4}}\right)\right)}{t^2}$ is positive and continuous. Therefore
\begin{equation*}
\begin{split}
c_1 = \inf_{0 \leq t \leq 1}\frac{\log\left(\frac{1}{2}\left(e^{-\frac{t}{4}} + e^{\frac{t}{4}}\right)\right)}{t^2} > 0.
\end{split}
\end{equation*}
Similarly,
\begin{equation*}
\begin{split}
C_1 = \sup_{0 \leq t \leq 1}\frac{\log\left(\frac{1}{2}\left(e^{-\frac{t}{2}} + e^{\frac{t}{2}}\right)\right)}{t^2} < \infty.
\end{split}
\end{equation*}
For all $0 \leq t \leq 1$, we have
\begin{equation*}
\begin{split}
e^{c_1kt^2} \leq \phi_{X}(t) \leq e^{C_1kt^2}.
\end{split}
\end{equation*}
By Theorem \ref{thm: general}, there exist constants $c, c', C > 0$, for all $0 \leq x \leq c'k$,
\begin{equation*}
\begin{split}
\bbP\left(X \geq x\right) \geq c\cdot e^{-C\frac{x^2}{k}},
\end{split}
\end{equation*}
which has finished the proof of this theorem.

\subsection{Proof of Theorem \ref{th:extreme-value-gaussian}}\label{sec:extreme-value-gaussian}
Let $Y_i = \frac{X_i}{\sqrt{\alpha\|u\|_2^2}}$, $i = 1, \dots, k$, immediately we have
\begin{equation*}
\begin{split}
c_2\exp\left(c_1t^2\right) \leq \phi_{Y_i}(t) \leq C_2\exp\left(C_1t^2\right), \quad \forall t \in \bbR.
\end{split}
\end{equation*}
By Theorem \ref{thm: general}, there exist two constants $\bar{C} > 1 > \bar{c} > 0$ such that
\begin{equation}\label{eq108}
\begin{split}
\bar{c}\exp\left(-\bar{C}x^2\right) \leq \bbP\left(Y_i \geq x\right), \bbP\left(Y_i \leq -x\right) \leq \bar{C}\exp\left(-\bar{c}x^2\right), \quad \forall x \geq 0, 1 \leq i \leq k.
\end{split}
\end{equation}
By the formula for expected value,
\begin{equation}\label{eq102}
\begin{split}
\bbE\left[\sup_{1 \leq i \leq k}Y_i\right] =& \bbE\left[\left(\sup_{1 \leq i \leq k}Y_i\right) \vee 0\right] + \bbE\left[\left(\sup_{1 \leq i \leq k}Y_i\right) \wedge 0\right]\\ =& \int_{0}^{\infty}\bbP\left(\left(\sup_{1 \leq i \leq k}Y_i\right) \vee 0 > t\right)dt - \int_{0}^{\infty}\bbP\left(\left(\sup_{1 \leq i \leq k}Y_i\right) \wedge 0 < -t\right)dt\\
=& \int_{0}^{\infty}\left[1 - \bbP\left(\sup_{1 \leq i \leq k}Y_i \leq t\right) - \bbP\left(\sup_{1 \leq i \leq k}Y_i < -t\right)\right]dt.
\end{split}
\end{equation}
First, we consider the upper bound. For all $t \geq \sqrt{\frac{\log\bar{C}}{\bar{c}}}$, we have $\bar{C}\exp\left(-\bar{c}t^2\right) \leq 1$. Thus
\begin{equation}\label{eq104}
\begin{split}
&1 - \bbP\left(\sup_{1 \leq i \leq k}Y_i \leq t\right) - \bbP\left(\sup_{1 \leq i \leq k}Y_i < -t\right)
= 1 - \prod_{i = 1}^k\left[1 - \bbP\left(Y_i > t\right)\right] - \prod_{i = 1}^k\bbP\left(Y_i < -t\right)\\
\leq& 1 - \prod_{i = 1}^{k}\left[1 - \bar{C}\exp\left(-\bar{c}t^2\right)\right] - 0 = 1 - \prod_{i = 1}^{k}\left[1 - \bar{C}\exp\left(-\bar{c}t^2\right)\right] \leq k\bar{C}\exp\left(-\bar{c}t^2\right).
\end{split}
\end{equation}
The last inequality holds since
$\prod_{i = 1}^{k}\left[1 - \bar{C}\exp\left(-\bar{c}t^2\right)\right] \geq 1 - \sum_{i = 1}^{k}\bar{C}\exp\left(-\bar{c}t^2\right) = 1 - k\bar{C}\exp\left(-\bar{c}t^2\right).$
Combine \eqref{eq102} and \eqref{eq104} together, we have
\begin{equation*}
\begin{split}
\bbE\left[\sup_{1 \leq i \leq k}Y_i\right] \leq& \sqrt{\frac{\log\left(\bar{C}k\right)}{\bar{c}}} + \int_{t = \sqrt{\frac{\log\left(\bar{C}k\right)}{\bar{c}}}}^{\infty}k\bar{C}\exp\left(-\bar{c}t^2\right)dt\\
=& \sqrt{\frac{\log\left(\bar{C}k\right)}{\bar{c}}} + \sqrt{\frac{1}{\bar{c}}}\int_{s = 0}^{\infty}\exp\left(-s^2 - 2s\sqrt{\log\left(\bar{C}k\right)}\right)ds\\
\leq& \sqrt{\frac{\log\left(\bar{C}k\right)}{\bar{c}}} + \sqrt{\frac{1}{\bar{c}}}\int_{s = 0}^{\infty}\exp\left(- 2s\sqrt{\log\left(\bar{C}k\right)}\right)ds\\
=& \sqrt{\frac{\log\left(\bar{C}k\right)}{\bar{c}}} + \sqrt{\frac{1}{\bar{c}}}\frac{1}{2\sqrt{\log\left(\bar{C}k\right)}}
\leq 2\sqrt{\frac{\log\left(\bar{C}k\right)}{\bar{c}}}.
\end{split}
\end{equation*}
Thus
\begin{equation*}
\begin{split}
\bbE \left[\sup_{1 \leq i \leq k}X_i\right] = \sqrt{\alpha}\bbE\left[\sup_{1 \leq i \leq k}Y_i\right] \leq 2\sqrt{\frac{\log\left(\bar{C}k\right)}{\bar{c}}\alpha} \leq 2\sqrt{\frac{\bar{C}\log k}{\bar{c}}\alpha}.
\end{split}
\end{equation*}
Next, we consider the lower bound. For all $t \geq 0$,
\begin{equation}\label{eq105}
\begin{split}
&1 - \bbP\left(\sup_{1 \leq i \leq k}Y_i \leq t\right)
= 1 - \prod_{i = 1}^k\left[1 - \bbP\left(Y_i > t\right)\right]
\geq 1 - \prod_{i = 1}^k\left[1 - \bar{c}\exp\left(-\bar{C}t^2\right)\right]\\ \geq& 1 - \prod_{i = 1}^k\exp\left(- \bar{c}\exp\left(-\bar{C}t^2\right)\right)
\geq 1 - \exp\left(- k\bar{c}\exp\left(-\bar{C}t^2\right)\right).
\end{split}
\end{equation}
In the second inequality, we used $0 \leq 1 - x \leq e^{-x}$ for all $x \in [0, 1]$.\\
In addition,
\begin{equation}\label{eq107}
\begin{split}
\bbP\left(\sup_{1 \leq i \leq k}Y_i < -t\right) = \prod_{i = 1}^k\bbP\left(Y_i < -t\right) \leq \prod_{i = 1}^k\bar{C}\exp\left(-\bar{c}t^2\right) = \bar{C}^k\exp\left(-\bar{c}kt^2\right).
\end{split}
\end{equation}
For all $t \geq \sqrt{\frac{2\log \bar{C}}{\bar{c}}}$,
\begin{equation*}
\begin{split}
\bbP\left(\sup_{1 \leq i \leq k}Y_i < -t\right) \leq \left(\exp\left(\frac{\bar{c}}{2}t^2\right)\right)^k\exp\left(-\bar{c}kt^2\right) = \exp\left(-\frac{\bar{c}}{2}kt^2\right).
\end{split}
\end{equation*}
By \eqref{eq105}, for all $0 \leq t \leq \sqrt{\frac{\log k}{\bar{C}}}$,
\begin{equation}\label{eq106}
\begin{split}
1 - \bbP\left(\sup_{1 \leq i \leq k}Y_i \leq t\right) \geq 1 - \exp\left(- k\bar{c}\exp\left(-\bar{C}\left(\sqrt{\frac{\log k}{\bar{C}}}\right)^2\right)\right) = 1 - e^{-\bar{c}}.
\end{split}
\end{equation}
Combine \eqref{eq102}, \eqref{eq107} and \eqref{eq106} together, and also notice that\\ $1 - \bbP\left(\sup_{1 \leq i \leq k}X_i \leq t\right) \geq 0, \bbP\left(\sup_{1 \leq i \leq k}Y_i < -t\right) \leq 1$ for all $t \geq 0$,
\begin{equation*}
\begin{split}
\bbE\left[\sup_{1 \leq i \leq k}Y_i\right] \geq& \int_{0}^{\sqrt{\frac{\log k}{\bar{C}}}}(1 - e^{-\bar{c}})dt - \int_{0}^{\sqrt{\frac{2\log \bar{C}}{\bar{c}}}}1dt - \int_{\sqrt{\frac{2\log \bar{C}}{\bar{c}}}}^{\infty}\exp\left(-\frac{\bar{c}}{2}kt^2\right)dt\\
\geq & \sqrt{\frac{\log k}{\bar{C}}}\left(1 - e^{-\bar{c}}\right) - \sqrt{\frac{2\log \bar{C}}{\bar{c}}} - \sqrt{\frac{\pi}{2\bar{c}k}}.
\end{split}
\end{equation*}
Hence there exists a constant $\widetilde{C} > 0$, for any $k \geq \widetilde{C}$,
\begin{equation}\label{eq109}
\begin{split}
\bbE\left[\sup_{1 \leq i \leq k}X_i\right] = \sqrt{\alpha\|u\|_2^2}\bbE\left[\sup_{1 \leq i \leq k}Y_i\right] \geq \frac{1 - e^{-\bar{c}}}{2}\sqrt{\frac{\alpha\|u\|_2^2\log k}{\bar{C}}}.
\end{split}
\end{equation}
Finally, for $2 \leq k \leq \widetilde{C}$, by \eqref{eq108},
\begin{equation*}
\begin{split}
\bbE\left[\sup_{1 \leq i \leq k}Y_i\right] \geq& \bbE \max\{Y_1, Y_2\} = \bbE Y_1 + \bbE\left[\max\{Y_1, Y_2\} - Y_1\right] = \bbE\max\{0, Y_2 - Y_1\}\\
\geq & \bbE 1_{\{Y_2 - Y_1 \geq 1\}} \geq \bbP\left(Y_1 \leq 0, Y_2 \geq 1\right) = \bbP\left(Y_1 \leq 0\right)\bbP\left(Y_2 \geq 1\right) \geq \bar{c}^2e^{-\bar{C}}.
\end{split}
\end{equation*}
Thus for all any $2 \leq k \leq \widetilde{C}$,
\begin{equation}\label{eq110}
\begin{split}
\bbE\left[\sup_{1 \leq i \leq k}X_i\right] = \sqrt{\alpha}\bbE\left[\sup_{1 \leq i \leq k}Y_i\right] \geq \frac{\bar{c}^2e^{-\bar{C}}}{\sqrt{\log \widetilde{C}}}\sqrt{\alpha\log k}.
\end{split}
\end{equation}
\eqref{eq109} and \eqref{eq110} show that there exists a constant $c > 0$, for any $k \geq 2$,
\begin{equation}
\begin{split}
\bbE\left[\sup_{1 \leq i \leq k}X_i\right] \geq c\sqrt{\alpha\|u\|_2^2\log k}.
\end{split}
\end{equation}

\subsection{Proof of Theorem \ref{th:extreme-value-exponential}}\label{sec:extreme-value-exponential}
By Theorem \ref{cor: sub-exponential}, there exist constants $\bar{C} > 1 > \bar{c}  > 0$ such that for all $x \geq 0$,
\begin{equation}\label{eq118}
\begin{split}
\bar{c}\exp\left(-\bar{C}\left(\frac{x^2}{\alpha\|u\|_2^2} \wedge \frac{Mx}{\|u\|_{\infty}}\right)\right) \leq \bbP\left(X_i \geq x\right) \leq \exp\left(-\bar{c}\left(\frac{x^2}{\alpha\|u\|_2^2} \wedge \frac{Mx}{\|u\|_{\infty}}\right)\right).
\end{split}
\end{equation}
Similarly to \eqref{eq102}, we have
\begin{equation}\label{eq111}
\begin{split}
\bbE \left[\sup_{1 \leq i \leq k}X_i\right] = \int_{0}^{\infty}\left[1 - \bbP\left(\sup_{1 \leq i \leq k}X_i \leq t\right)\right]dt - \int_{0}^{\infty}\bbP\left(\sup_{1 \leq i \leq k}X_i < -t\right)dt.
\end{split}
\end{equation}
First, we consider the upper bound. For all $t \geq 2\left(\sqrt{\frac{\alpha\|u\|_2^2\log k}{\bar{c}}} \vee \frac{\|u\|_{\infty}\log k}{\bar{c}M}\right) \geq \sqrt{\frac{\alpha\|u\|_2^2}{\bar{c}}} \vee \frac{\|u\|_{\infty}}{\bar{c}M}$,
\begin{equation*}
\begin{split}
&1 - \bbP\left(\sup_{1 \leq i \leq k}X_i \leq t\right) = 1 - \prod_{i = 1}^k\left[1 - \bbP\left(X_i > t\right)\right] \leq 1 - \prod_{i = 1}^k\left[1 - \exp\left(-\bar{c}(\frac{t^2}{\alpha\|u\|_2^2} \wedge \frac{Mt}{\|u\|_{\infty}})\right)\right]\\
\leq& 1 - \left[1 - \sum_{i = 1}^k\exp\left(-\bar{c}\left(\frac{t^2}{\alpha\|u\|_2^2} \wedge \frac{Mt}{\|u\|_{\infty}}\right)\right)\right]
\leq k\left[\exp\left(-\bar{c}\frac{t^2}{\alpha\|u\|_2^2}\right) + \exp\left(-\bar{c}\frac{Mt}{\|u\|_{\infty}}\right)\right].
\end{split}
\end{equation*}
The second inequality comes from Bernoulli's inequality.\\
By \eqref{eq111} and the previous inequality,
\begin{equation}\label{eq112}
\begin{split}
\bbE \left[\sup_{1 \leq i \leq k}X_i\right] \leq& \int_{0}^{\infty}\left[1 - \bbP\left(\sup_{1 \leq i \leq k}X_i \leq t\right)\right]dt\\
\leq& \int_{0}^{2\left(\sqrt{\frac{\alpha\|u\|_2^2\log k}{\bar{c}}} \vee \frac{\|u\|_{\infty}\log k}{\bar{c}M}\right)}1 dt + \int_{2\left(\sqrt{\frac{\alpha\|u\|_2^2\log k}{\bar{c}}} \vee \frac{\|u\|_{\infty}\log k}{\bar{c}M}\right)}^{\infty}k\exp\left(-\bar{c}\frac{t^2}{\alpha\|u\|_2^2}\right)dt\\
&+ \int_{2\left(\sqrt{\frac{\alpha\|u\|_2^2\log k}{\bar{c}}} \vee \frac{\|u\|_{\infty}\log k}{\bar{c}M}\right)}^{\infty}k\exp\left(-\bar{c}\frac{Mt}{\|u\|_{\infty}}\right)dt\\
\leq& 2\left(\sqrt{\frac{\alpha\|u\|_2^2\log k}{\bar{c}}} \vee \frac{\|u\|_{\infty}\log k}{\bar{c}M}\right) + \sqrt{\frac{\alpha\|u\|_2^2}{\bar{c}}}\int_{2\sqrt{\log k}}^{\infty}k\exp\left(-u^2\right)du\\
&+ \frac{\|u\|_{\infty}}{\bar{c}M}\int_{2\log k}^{\infty}k\exp\left(-u\right)du.
\end{split}
\end{equation}
Moreover, we have
\begin{equation}\label{eq113}
\begin{split}
\int_{2\sqrt{\log k}}^{\infty}k\exp\left(-u^2\right)du \leq& \int_{\sqrt{\log k}}^{\infty}k\exp\left(-u^2\right)du = \int_{0}^{\infty}\exp\left(-s^2 - 2s\sqrt{\log k}\right)ds\\
\leq& \int_{0}^{\infty}\exp\left( - 2s\sqrt{\log k}\right)ds = \frac{1}{2\sqrt{\log k}},
\end{split}
\end{equation}
and
\begin{equation}\label{eq114}
\begin{split}
\int_{2\log k}^{\infty}k\exp\left(-u\right)du = k\exp\left(-2\log k\right) = \frac{1}{k}.
\end{split}
\end{equation}
\eqref{eq112}, \eqref{eq113} and \eqref{eq114} together imply
\begin{equation}
\begin{split}
\bbE\left[\sup_{1 \leq i \leq k}X_i\right] \leq& 2\left(\sqrt{\frac{\alpha\|u\|_2^2\log k}{\bar{c}}} \vee \frac{\|u\|_{\infty}\log k}{\bar{c}M}\right) + \sqrt{\frac{\alpha\|u\|_2^2}{\bar{c}}}\frac{1}{2\sqrt{\log k}} + \frac{\|u\|_{\infty}}{\bar{c}M}\frac{1}{k}\\
\leq& 4\left(\sqrt{\frac{\alpha\|u\|_2^2\log k}{\bar{c}}} \vee \frac{\|u\|_{\infty}\log k}{\bar{c}M}\right) \leq \frac{4}{\bar{c}}\left[\sqrt{\alpha\|u\|_2^2\log k} \vee \frac{\|u\|_{\infty}\log k}{M}\right].
\end{split}
\end{equation}
Next, we consider the lower bound. For any $0 \leq t \leq \sqrt{\frac{\alpha\|u\|_2^2\log k}{\bar{C}}} \vee \frac{\|u\|_{\infty}\log k}{\bar{C}M}$,
\begin{equation*}
\begin{split}
1 - \bbP\left(\sup_{1 \leq i \leq k}X_i \leq t\right) =& 1 - \prod_{i = 1}^k\left[1 - \bbP\left(X_i > t\right)\right] \geq 1 - \prod_{i = 1}^k\left[1 - \bar{c}\exp\left(-\bar{C}(\frac{t^2}{\alpha\|u\|_2^2} \wedge \frac{Mt}{\|u\|_{\infty}})\right)\right]\\
\geq& 1 - \prod_{i = 1}^k\exp\left(-\bar{c}\exp\left(-\bar{C}\left(\frac{t^2}{\alpha\|u\|_2^2} \wedge \frac{Mt}{\|u\|_{\infty}}\right)\right)\right)\\
=& 1 - \exp\left(-\bar{c}k\exp\left(-\bar{C}\left(\frac{t^2}{\alpha\|u\|_2^2} \wedge \frac{Mt}{\|u\|_{\infty}}\right)\right)\right) \geq 1 - e^{-\bar{c}}.
\end{split}
\end{equation*}
Thus
\begin{equation}\label{eq116}
\begin{split}
\int_{0}^{\infty}\left[1 - \bbP\left(\sup_{1 \leq i \leq k}X_i \leq t\right)\right]dt \geq& \int_{0}^{\sqrt{\frac{\alpha\|u\|_2^2\log k}{\bar{C}}} \vee \frac{\|u\|_{\infty}\log k}{\bar{C}M}}\left[1 - \bbP\left(\sup_{1 \leq i \leq k}X_i \leq t\right)\right]dt\\
\geq& \left(1 - e^{-\bar{c}}\right)\left[\sqrt{\frac{\alpha\|u\|_2^2\log k}{\bar{C}}} \vee \frac{\|u\|_{\infty}\log k}{\bar{C}M}\right].
\end{split}
\end{equation}
Note that for all $-M \leq t \leq 0$, $\phi_{Z_{ij}}(t) \leq \exp\left(C_1\alpha t^2\right)$,
by \eqref{eq121},
\begin{equation*}
\begin{split}
\bbP\left(X_i \leq -x\right) = \bbP\left(-X_i \geq x\right) \leq \left\{\begin{array}{l}
e^{-\frac{x^2}{4C_1\alpha\|u\|_2^2}}, \quad 0 \leq x \leq \frac{2MC_1\alpha\|u\|_2^2}{\|u\|_{\infty}},\\
e^{-\frac{Mx}{2\|u\|_{\infty}}}, \quad x > \frac{2MC_1\alpha\|u\|_2^2}{\|u\|_{\infty}},
\end{array}\right.
\end{split}
\end{equation*}
which means that for all $t \geq 0$,
\begin{equation*}
	\begin{split}
	\bbP\left(X_i \leq -t\right) \leq \exp\left(-\frac{t^2}{4C_1\alpha\|u\|_2^2} \wedge \frac{Mt}{2\|u\|_{\infty}}\right) \leq \exp\left(-\frac{1}{4C_1}\left(\frac{t^2}{\alpha\|u\|_2^2} \wedge \frac{Mt}{\|u\|_{\infty}}\right)\right).
	\end{split}
\end{equation*}
Therefore
\begin{equation*}
\begin{split}
\int_{0}^{\infty}\bbP\left(\sup_{1 \leq i \leq k}X_i < -t\right)dt =& \int_{0}^{\infty}\prod_{i = 1}^k\bbP\left(X_i < -t\right)dt\\ \leq& \int_{0}^{\infty}\prod_{i = 1}^k\exp\left(-\frac{1}{4C_1}\left(\frac{t^2}{\alpha\|u\|_2^2} \wedge \frac{Mt}{\|u\|_{\infty}}\right)\right)dt\\
\leq&  \int_{0}^{\infty}\exp\left(-\frac{k}{4C_1}\frac{t^2}{\alpha\|u\|_2^2}\right)dt +  \int_{0}^{\infty}\exp\left(-\frac{k}{4C_1} \frac{Mt}{\|u\|_{\infty}}\right)dt\\
=& \frac{\sqrt{\pi}}{2}\sqrt{\frac{\alpha\|u\|_2^2}{k}\cdot 4C_1} + 4C_1\frac{\|u\|_{\infty}}{Mk}.
\end{split}
\end{equation*}
By \eqref{eq111}, \eqref{eq116} and the previous inequality, there exists a constant $\widetilde{C} > 0$ such that for all $k \geq \widetilde{C}$,
\begin{equation}\label{eq120}
\begin{split}
\bbE\left[\sup_{1 \leq i \leq k}X_i\right] \geq& \frac{1 - e^{-\bar{c}}}{2}\left[\sqrt{\frac{\alpha\|u\|_2^2\log k}{\bar{C}}} \vee \frac{\|u\|_{\infty}\log k}{\bar{C}M}\right]\\
\geq& \frac{1 - e^{-\bar{c}}}{2\bar{C}}\left[\sqrt{\alpha\|u\|_2^2\log k} \vee \frac{\|u\|_{\infty}\log k}{M}\right].
\end{split}
\end{equation}
Now, we consider the case $2 \leq k \leq \widetilde{C}$. Notice that
\begin{equation*}
\begin{split}
\bbE\left[X_1 \vee X_2\right] =& \bbE X_1 + \bbE\left[0 \vee (X_2 - X_1)\right]\\ \geq& 0 + \frac{1}{2}\left(\sqrt{\frac{\alpha\|u\|_2^2}{\bar{C}}} \vee \frac{\|u\|_{\infty}}{\bar{C}M}\right)\bbP\left(X_2 - X_1 \geq \frac{1}{2}\left(\sqrt{\frac{\alpha\|u\|_2^2}{\bar{C}}} \vee \frac{\|u\|_{\infty}}{\bar{C}M}\right)\right).
\end{split}
\end{equation*}
By \eqref{eq118},
\begin{equation*}
\begin{split}
&\bbP\left(X_2 - X_1 \geq \frac{1}{2}\left(\sqrt{\frac{\alpha\|u\|_2^2}{\bar{C}}} \vee \frac{\|u\|_{\infty}}{\bar{C}M}\right)\right)\\\geq& \bbP\left(X_2 \geq \sqrt{\frac{\alpha\|u\|_2^2}{\bar{C}}} \vee \frac{\|u\|_{\infty}}{\bar{C}M}\right)\cdot\bbP\left(X_1 \leq \frac{1}{2}\sqrt{\frac{\alpha\|u\|_2^2}{\bar{C}}} \vee \frac{\|u\|_{\infty}}{\bar{C}M}\right)\\
\geq& \bar{c}e^{-1}\cdot\left(1 - e^{-\frac{\bar{c}}{4\bar{C}}}\right).    	
\end{split}
\end{equation*}
Thus
\begin{equation*}
\begin{split}
\bbE\left[X_1 \vee X_2\right] \geq \frac{1}{2}\bar{c}e^{-1}\cdot\left(1 - e^{-\frac{\bar{c}}{4\bar{C}}}\right)\left(\sqrt{\frac{\alpha\|u\|_2^2}{\bar{C}}} \vee \frac{\|u\|_{\infty}}{\bar{C}M}\right).
\end{split}
\end{equation*}
Therefore, there exists a constants $\widetilde{c} > 0$ such that for all $2 \leq k \leq \widetilde{C}$,
\begin{equation}\label{eq119}
\begin{split}
\bbE\left[\sup_{1 \leq i \leq k}X_i\right] \geq& \bbE\left[X_1 \vee X_2\right] \geq \widetilde{c}\left[\sqrt{\alpha\|u\|_2^2\log \widetilde{C}} \vee \frac{\|u\|_{\infty}\log \widetilde{C}}{M}\right]\\ \geq& \widetilde{c}\left[\sqrt{\alpha\|u\|_2^2\log k} \vee \frac{\|u\|_{\infty}\log k}{M}\right].
\end{split}
\end{equation}
Combine \eqref{eq119} and \eqref{eq120} together, we have reached the conclusion.\quad $\square$

\subsection{Proof of Theorem \ref{thm: signal}}\label{sec:proof-signal}
Since $\lambda > \mu$, a natural idea to separate the signal and noise is to set a cut-off value $\theta$, and classify $Y_i$ to be a signal if $Y_i \geq \theta$. Such a scheme is indeed optimal in minimizing the Hamming distance misclassification rate due to the following lemma.
\begin{Lemma}\label{lm: classification}
	$\widetilde{Z}_i = 1_{\{Y_i > \widetilde{\theta}\}}$ is the best identification strategy.
\end{Lemma}
{\noindent\bf Proof of Lemma \ref{lm: classification}.}
\begin{equation*}
\begin{split}
\bbE kM(\hat{Z}) =& \bbE\sum_{i = 1}^{k}\left|\hat{Z}_i - Z_i\right| = \sum_{i = 1}^{k}\sum_{j = 0}^{\infty}\left[\hat{Z}_i(j)\bbP\left(Y_i = j, Z_i = 0\right) + (1 - \hat{Z}_i(j))\bbP\left(Y_i = j, Z_i = 1\right)\right]\\
\geq& \sum_{i = 1}^{k}\sum_{j = 0}^{\infty}\min\{\bbP\left(Y_i = j, Z_i = 0\right), \bbP\left(Y_i = j, Z_i = 1\right)\}.
\end{split}
\end{equation*}
"=" holds if and only if
\begin{equation}\label{eq99}
\begin{split}
\hat{Z}_i(j) = \left\{\begin{array}{ll}
0, & \bbP\left(Y_i = j, Z_i = 0\right) > \bbP\left(Y_i = j, Z_i = 1\right);\\
1, & \bbP\left(Y_i = j, Z_i = 0\right) < \bbP\left(Y_i = j, Z_i = 1\right).
\end{array}\right.
\end{split}
\end{equation}
For all $1 \leq i \leq k$,
\begin{equation*}
\begin{split}
\frac{\bbP\left(Y_i = j, Z_i = 0\right)}{\bbP\left(Y_i = j, Z_i = 1\right)} = \frac{\left(1 - \varepsilon\right)\frac{\mu^j}{j!}e^{-\mu}}{\varepsilon\frac{\lambda^j}{j!}e^{-\lambda}} = \frac{1 - \varepsilon}{\varepsilon}\left(\frac{\mu}{\lambda}\right)^j e^{\lambda - \mu}.
\end{split}
\end{equation*}
We can see if $j > \widetilde\theta = \frac{\log\left(\frac{1 - \varepsilon}{\varepsilon}\right) + \lambda - \mu}{\log\left(\frac{\lambda}{\mu}\right)}$, $\bbP\left(Y_i = j, Z_i = 0\right) < \bbP\left(Y_i = j, Z_i = 1\right)$; if $j < \widetilde\theta = \frac{\log\left(\frac{1 - \varepsilon}{\varepsilon}\right) + \lambda - \mu}{\log\left(\frac{\lambda}{\mu}\right)}$, $\bbP\left(Y_i = j, Z_i = 0\right) > \bbP\left(Y_i = j, Z_i = 1\right)$. Thus \eqref{eq99} is equivalent to
\begin{equation*}
\begin{split}
\hat{Z}_i = \left\{\begin{array}{ll}
0, & Y_i < \widetilde{\theta};\\
1, & Y_i > \widetilde{\theta},
\end{array}\right.
\end{split}
\end{equation*}
for all $1 \leq i \leq k$,
which means that our assertion is true. \qed

We go back to the proof of Theorem \ref{thm: signal}.
Notice that $(Y_i, Z_i, \widetilde{Z}_i)$ are i.i.d. and
\begin{equation*}
\begin{split}
\bbE M(\widetilde{Z}_i) = \frac{1}{k}\bbE \sum_{i = 1}^{k}\left|\widetilde{Z}_i - Z_i\right| = \bbE\left|\widetilde{Z}_1 - Z_1\right|.
\end{split}
\end{equation*}
By Lemma \ref{lm: classification},
\begin{equation}\label{eq97}
\begin{split}
\bbE M(\hat{Z}_i) \geq \bbE M(\widetilde{Z}_i) = \left[(1 - \varepsilon)\bbP_{\mu}\left(Y_1 > \widetilde\theta\right) + \varepsilon\bbP_{\lambda}\left(Y_1 \leq \widetilde\theta\right)\right].
\end{split}
\end{equation}
\begin{itemize}[leftmargin = *]
	\item \bm{$\varepsilon^-\leq \varepsilon \leq \varepsilon^+$}.
	Here
	\begin{equation*}
	\begin{split}
	\mu = \frac{\log\left(\frac{1 - \varepsilon^+}{\varepsilon^+}\right) + \lambda - \mu}{\log\left(\frac{\lambda}{\mu}\right)} \leq \widetilde{\theta} = \frac{\log\left(\frac{1 - \varepsilon}{\varepsilon}\right) + \lambda - \mu}{\log\left(\frac{\lambda}{\mu}\right)} \leq  \frac{\log\left(\frac{1 - \varepsilon^-}{\varepsilon^-}\right) + \lambda - \mu}{\log\left(\frac{\lambda}{\mu}\right)} = \lambda.
	\end{split}
	\end{equation*}
	By \eqref{ineq:previous-poisson1} \eqref{ineq:previous-poisson2}, 
	\begin{equation}\label{eq94}
	\begin{split}
	&(1 - \varepsilon)\bbP_{\mu}\left(Y_1 > \widetilde\theta\right) + \varepsilon\bbP_{\lambda}\left(Y_1 \leq \widetilde\theta\right)\\ \leq& (1 - \varepsilon)\exp\left(-\frac{(\widetilde\theta - \mu)^2}{2\mu}\psi_{Benn}\left(\frac{\widetilde\theta - \mu}{\mu}\right)\right) + \varepsilon\exp\left(-\frac{(\lambda - \widetilde\theta)^2}{2\lambda}\psi_{Benn}\left(\frac{\widetilde\theta - \lambda}{\lambda}\right)\right)\\
	=& \exp\left(-g(\widetilde{\theta})\right).
	\end{split}
	\end{equation}
	Note that $\lambda \leq C_2\mu$, we know that $\lambda - \theta \leq \lambda - \mu \leq \lambda - \frac{\lambda}{C_2} = \lambda\left(1 - \frac{1}{C_2}\right)$.	
	By Theorem \ref{cor:Poisson tail lower bound}, there exist two constants $\bar c  > 0, C > 1$ such that
	\begin{equation}\label{eq95}
	\begin{split}
	&(1 - \varepsilon)\bbP_{\mu}\left(Y_1 > \widetilde\theta\right) + \varepsilon\bbP_{\lambda}\left(Y_1 \leq \widetilde\theta\right)\\ \geq& (1 - \varepsilon)\bar c\exp\left(-C\frac{(\widetilde\theta - \mu)^2}{2\mu}\psi_{Benn}\left(\frac{\widetilde\theta - \mu}{\mu}\right)\right) + \varepsilon \bar c\exp\left(-C\frac{(\lambda - \widetilde\theta)^2}{2\lambda}\psi_{Benn}\left(\frac{\widetilde\theta - \lambda}{\lambda}\right)\right)\\
	\geq&\bar c\left[(1 - \varepsilon)^C\exp\left(-C\frac{(\widetilde\theta - \mu)^2}{2\mu}\psi_{Benn}(\frac{\widetilde\theta - \mu}{\mu})\right) + \varepsilon^C\exp\left(-C\frac{(\lambda - \widetilde\theta)^2}{2\lambda}\psi_{Benn}(\frac{\widetilde\theta - \lambda}{\lambda})\right)\right]\\
	\geq&\bar c\cdot 2^{1 - C}\left[(1 - \varepsilon)\exp\left(\frac{(\widetilde\theta - \mu)^2}{2\mu}\psi_{Benn}(\frac{\widetilde\theta - \mu}{\mu})\right) + \varepsilon\exp\left(-\frac{(\lambda - \widetilde\theta)^2}{2\lambda}\psi_{Benn}(\frac{\widetilde\theta - \lambda}{\lambda})\right)\right]^C\\
	=& \bar c\cdot 2^{1 - C}\exp\left(-Cg(\widetilde{\theta})\right).
	\end{split}
	\end{equation}
	In the last inequality, we used $x^C + y^C \geq 2^{1 - C}(x + y)^C$ for all $x, y \geq 0, C \geq 1$.
	Combining \eqref{eq97}, \eqref{eq94} and \eqref{eq95} together, there exist two constants $C, c > 0$, for  any classification $\hat{Z} \in \{0, 1\}^k$ based on $\{Y_i\}_{i = 1}^k$,
	\begin{equation*}
	\begin{split}
	 \bbE M(\hat{Z}) \geq c\exp\left(-Cg(\widetilde{\theta})\right),
	\end{split}
	\end{equation*}
	and
	\begin{equation*}
	\begin{split}
	c\exp\left(-Cg(\widetilde{\theta})\right) \leq \bbE M(\hat{Z}) \leq \exp\left(-g(\widetilde{\theta})\right).
	\end{split}
	\end{equation*}
	
	\item \bm{$\varepsilon^+ <\varepsilon < 1$}.
	\begin{equation*}
	\begin{split}
	\widetilde{\theta} < \frac{\log\left(\frac{1 - \varepsilon^+}{\varepsilon^+}\right) + \lambda - \mu}{\log\left(\frac{\lambda}{\mu}\right)} = \mu.
	\end{split}
	\end{equation*}
	For the lower bound, by Theorem \ref{cor:Poisson tail lower bound}, there exists a constant $c > 0$,
	\begin{equation*}
	\begin{split}
	(1 - \varepsilon)\bbP_{\mu}\left(Y_1 > \widetilde\theta\right) + \varepsilon\bbP_{\lambda}\left(Y_1 \leq \widetilde\theta\right) \geq (1 - \varepsilon)\bbP_{\mu}\left(Y_1 > \widetilde\theta\right) \geq (1 - \varepsilon)\bbP_{\mu}\left(Y_1 - \mu \geq 0\right) \geq (1 - \varepsilon)c.
	\end{split}
	\end{equation*}
	\eqref{eq97} and the previous inequality show that
	\begin{equation}
	\begin{split}
	\bbE M(\widetilde{Z}_i) = \left[(1 - \varepsilon)\bbP_{\mu}\left(Y_1 > \widetilde\theta\right) + \varepsilon\bbP_{\lambda}\left(Y_1 \leq \widetilde\theta\right)\right] \geq c(1 - \varepsilon).
	\end{split}
	\end{equation}
	holds for any classification $\hat{Z} \in \{0, 1\}^k$ based on $\{Y_i\}_{i = 1}^k$.\\
	For the upper bound, by \eqref{eq97},
	\begin{equation}
	\begin{split}
	\bbE M(\widetilde{Z}_i) \leq \bbE\left|1 - Z_1\right| = \bbP\left(Z_1 = 0\right) = 1 - \varepsilon.
	\end{split}
	\end{equation}
	\item \bm{$0 < \varepsilon < \varepsilon^+$}. For the lower bound, note that
	\begin{equation*}
	\begin{split}
	\widetilde{\theta} > \frac{\log\left(\frac{1 - \varepsilon^-}{\varepsilon^-}\right) + \lambda - \mu}{\log\left(\frac{\lambda}{\mu}\right)} = \lambda.
	\end{split}
	\end{equation*}
	By Theorem 11 (set $\beta = 2$), there exists a constant $c > 0$,
	\begin{equation*}
	\begin{split}
	(1 - \varepsilon)\bbP_{\mu}\left(Y_1 > \widetilde\theta\right) + \varepsilon\bbP_{\lambda}\left(Y_1 \leq \widetilde\theta\right) \geq \varepsilon\bbP_{\lambda}\left(Y_1 \leq \widetilde\theta\right) \geq \varepsilon\bbP_{\lambda}\left(Y_1 - \lambda \leq 0\right) \geq c\varepsilon.
	\end{split}
	\end{equation*}
	\eqref{eq97} and the previous inequality imply that for any classification $\hat{Z} \in \{0, 1\}^k$ based on $\{Y_i\}_{i = 1}^k$,
	\begin{equation}
	\begin{split}
	\bbE M(\hat{Z}) \geq \bbE M(\widetilde{Z}) = (1 - \varepsilon)\bbP_{\mu}\left(Y_1 > \widetilde\theta\right) + \varepsilon\bbP_{\lambda}\left(Y_1 \leq \widetilde\theta\right) \geq c\varepsilon.
	\end{split}
	\end{equation}
	For the upper bound, by \eqref{eq97},
	\begin{equation}
	\begin{split}
	\bbE M(\widetilde{Z}_i) \leq \bbE\left|0 - Z_1\right| = \bbP\left(Z_1 = 1\right) = \varepsilon.
	\end{split}
	\end{equation}
\end{itemize}
In summary, we have arrived at the conclusion.\quad $\square$

\subsection{Proof of Lemma \ref{lm:imply}}\label{sec:proof-imply}
If $\bbP\left(|Z_i| \geq x\right)\leq C_3\exp(-c_3x^2)$ holds for all $x \geq 0$, 
\begin{equation*}
\begin{split}
\bbE |Z_i|^p =& \int_{0}^{\infty}\bbP\left(|Z_i|^p > t\right)dt = \int_{0}^{\infty}\bbP\left(|Z_i| > s\right)ps^{p - 1}ds \leq p\int_{0}^{\infty}C_3e^{-c_3s^2}s^{p - 1}ds\\ =& pC_3c_3^{-\frac{p}{2}}\frac{1}{2}\int_{0}^{\infty}e^{-u}u^{\frac{p}{2}-1}du = \frac{pC_3}{2}c_3^{-\frac{p}{2}}\Gamma\left(\frac{p }{2}\right) \leq \frac{pC_3}{2}c_3^{-\frac{p}{2}}\left(\frac{p}{2}\right)^{\frac{p}{2}} = \frac{pC_3}{2}\left(\frac{p}{2c_3}\right)^{\frac{p}{2}}.
\end{split}
\end{equation*}
Thus for $p \geq 1$,
\begin{equation}\label{eq88}
\begin{split}
\left(\bbE |Z_i|^p\right)^{1/p} \leq \left(\frac{pC_3}{2}\left(\frac{p}{2c_3}\right)^{\frac{p}{2}}\right)^{\frac{1}{p}} = p^{\frac{1}{p}}\left(\frac{C_3}{2}\right)^{1/p}\left(\frac{p}{2c_3}\right)^{\frac{1}{2}} \leq 2\max\left\{1, \frac{C_3}{2}\right\}\left(\frac{p}{2c_3}\right)^{\frac{1}{2}} = C_4\sqrt{p}.
\end{split}
\end{equation}
Here $C_4 = 2\max\{1, C_3/2\}\left(\frac{1}{2C_3}\right)^{\frac{1}{2}}$. The second inequality follows from $p^{\frac{1}{p}} \leq 2$ and\\ $(C_3/2)^{1/p} \leq \max\{1, C_3/2\}$ for $p \geq 1$.\\
By Taylor's expansion,
\begin{equation}\label{eq89}
\begin{split}
\bbE e^{tZ_i} = 1  + t\bbE Z_i + \sum_{p = 2}^\infty \frac{t^p\bbE Z_i^p}{p!} \leq 1  + \sum_{p = 2}^\infty \frac{t^p\bbE |Z_i|^p}{p!} \leq 1 + \sum_{p = 2}^\infty\frac{t^p\left(C_4\sqrt{p}\right)^p}{p!} \leq 1 + \sum_{p = 2}^\infty\left(\frac{C_5t}{\sqrt{p}}\right)^p.
\end{split}
\end{equation}
Here $C_5 = eC_4$. The first inequality follows from $\bbE Z_i = 0$; in the second one we use \eqref{eq88}; the third one holds since $p! \geq \left(p/e\right)^p$.
\begin{itemize}[leftmargin = *]
	\item \bm{$0 \leq t \leq \frac{1}{2C_5}$}.
	By \eqref{eq89},
	\begin{equation*}
	\begin{split}
	\bbE e^{tZ_i} \leq 1 + \sum_{p = 2}^{\infty}\left(C_5t\right)^p = 1 + \frac{C_5^2t^2}{1 - C_5t} \leq 1 + 2C_5^2t^2 \leq e^{2C_5^2t^2}.
	\end{split}
	\end{equation*}
	\item
	\bm{$t > \frac{1}{2C_5}$}.
	Set $c_4 = \frac{1}{4eC_4^2}$, by \eqref{eq88} and $p! \geq (p/e)^p$,
	\begin{equation*}
	\begin{split}
	\bbE e^{c_4Z_i^2} = \sum_{p = 0}^{\infty}\frac{c_4^p \bbE Z_i^{2p}}{p!} \leq 1 + \sum_{p = 1}^{\infty}\frac{c_4^p\left(C_4\sqrt{2p}\right)^{2p}}{\left(\frac{p}{e}\right)^p} = \sum_{p = 0}^{\infty}\left(2ec_4C_4^2\right)^p = \sum_{p = 0}^{\infty}\left(\frac{1}{2}\right)^p = 2.
	\end{split}
	\end{equation*}
	Thus
	\begin{equation*}
	\begin{split}
	\bbE e^{tZ_i} \leq \bbE e^{c_4Z_i^2 + \frac{t^2}{4c_4}} \leq 2e^{\frac{t^2}{4c_4}} \leq e^{4C_5^2t^2}e^{\frac{t^2}{4c_4}} = e^{\left(4C_5^2 + eC_4^2\right)t^2}.
	\end{split}
	\end{equation*}
	The first inequality holds since $c_4Z_i^2 + \frac{t^2}{4c_4} \geq tZ_i$; the third inequality is due to $t > \frac{1}{2C_5}$.
\end{itemize}
Therefore
\begin{equation*}
\begin{split}
\phi_{Z_i}(t) \leq \exp\left(C_1t^2\right), \quad \forall t \geq 0.
\end{split}
\end{equation*}
If $\bbP\left(Z_i \geq x\right) \geq c_2\exp\left(-C_2x^2\right)$ for all $x \geq 0$, for any $t\geq 0$,
\begin{equation*}
\begin{split}
\mathbb{E}\exp(tZ_i) \geq& \mathbb{E}\exp(tZ_i) 1_{\{Z_i \geq \frac{t}{2C_2}\}} \geq \exp\left(\frac{t^2}{2C_2}\right)\cdot \bbP\left(Z_i \geq \frac{t}{2C_2}\right)\\ \geq& \exp\left(\frac{t^2}{2C_2}\right)\cdot c_2\exp\left(-C_2 \left(\frac{t}{2C_2}\right)^2\right) = c_2\exp\left(\frac{t^2}{4C_2}\right).
\end{split}
\end{equation*}
If $c_2 \geq 1$, then immediately we have
\begin{equation}\label{eq70}
\begin{split}
\bbE\exp\left(tZ_i\right) \geq \exp\left(\frac{t^2}{4C_2}\right), \quad \forall t \geq 0.
\end{split}
\end{equation}
If $c_2 < 1$, then for all $t \geq \sqrt{8C_2\log\left(\frac{1}{c_2}\right)}$,
\begin{equation}\label{eq69}
\begin{split}
\bbE\exp\left(tZ_i\right) \geq c_2\exp\left(\frac{t^2}{8C_2}\right)\cdot\exp\left(\frac{t^2}{8C_2}\right) \geq \exp\left(\frac{t^2}{8C_2}\right).
\end{split}
\end{equation}
For any $0 \leq t \leq \sqrt{8C_2\log\left(\frac{1}{c_2}\right)}$, by Taylor's Theorem,
\begin{equation*}
\begin{split}
e^{tu} = 1 + tu + e^{\xi}\frac{(tu)^2}{2},
\end{split}
\end{equation*}
where $\xi = \xi_{tu}$ is a real number between $0$ and $tu$.\\ For $u \geq 1, t \geq 0$,
\begin{equation*}
\begin{split}
e^{\xi}\frac{(tu)^2}{2} \geq \frac{(tu)^2}{2} \geq \frac{t^2}{2}.
\end{split}
\end{equation*}
Note that $\bbE X = 0$, for all $t \geq 0$, we have
\begin{equation}\label{eq67}
\begin{split}
\bbE e^{tX} =& \bbE\left(1 + tX + e^{\xi}\frac{(tX)^2}{2}\right) = 1 + \bbE\left(e^{\xi}\frac{(tX)^2}{2}\right) \geq 1 + \bbE\left(e^{\xi}\frac{(tX)^2}{2}1_{\{X \geq 1\}}\right)\\
\geq& 1 + \bbE \left(\frac{t^2}{2}1_{\{X \geq 1\}}\right) = 1 + \frac{t^2}{2}\bbP\left(X \geq 1\right) \geq 1 + c_2e^{-C_2}\frac{t^2}{2}.
\end{split}
\end{equation}
In the last step, we used the condition $\bbP\left(Z_i \geq x\right) \geq c_2\exp\left(-C_2x^2\right)$ for all $x \geq 0$.
Since $\frac{\log(1 + x)}{x}$ is a decreasing function of $x > 0$,
\begin{equation*}
\begin{split}
\bar c = \inf_{0 \leq t \leq c_2e^{-C_2}\cdot4C_2\log\left(\frac{1}{c_2}\right)}\frac{\log(1 + x)}{x} = \frac{\log\left(1 + 4c_2e^{-C_2}\cdot C_2\log\left(\frac{1}{c_2}\right)\right)}{4c_2e^{-C_2}\cdot C_2\log\left(\frac{1}{c_2}\right)}> 0,
\end{split}
\end{equation*}
Then for all $0 \leq t \leq \sqrt{8C_2\log\left(\frac{1}{c_2}\right)}$,
\begin{equation}
\begin{split}
\log\left(1 + c_2e^{-C_2}\frac{t^2}{2}\right) \geq \bar c\cdot c_2e^{-C_2}\frac{t^2}{2}.
\end{split}
\end{equation}
By \eqref{eq67} and the previous inequality, we know that for all $0 \leq t \leq \sqrt{8C_2\log\left(\frac{1}{c_2}\right)}$,
\begin{equation}\label{eq68}
\begin{split}
\bbE e^{tX} \geq \exp\left(\frac{\bar{c}\cdot c_2e^{-C_2}}{2}t^2\right).
\end{split}
\end{equation}
\eqref{eq70}, \eqref{eq69}, and \eqref{eq68} together imply that
\begin{equation}
\begin{split}
\bbE\left(tZ_i\right) \geq \exp(c_1t^2), \quad \forall t \geq 0,
\end{split}
\end{equation}
where $c_1 = \frac{1}{4C_2}$ if $c_2 \geq 1$ and $c_1 = \min\{\frac{1}{8C_2}, \frac{\bar{c}\cdot c_2e^{-C_2}}{2}\}$ if $0< c_2 < 1$.
Thus if Statement 2 (Equation \ref{ineq:sub-gaussian-state2}) holds, Statement 1 (Equation \ref{ineq:sub-gaussian-state1}) also holds.\quad $\square$


\end{document}